\documentclass[12pt]{article}

\usepackage{epsfig}

\usepackage{amssymb}

\usepackage{amssymb}

\newtheorem{theorem}{Theorem}[section]
\newtheorem{lemma}[theorem]{Lemma}
\newtheorem{corollary}[theorem]{Corollary}
\newtheorem{proposition}[theorem]{Proposition}

\newtheorem{problem}[theorem]{Problem}

\newtheorem{definition}[theorem]{Definition}
\newtheorem{example}[theorem]{Example}
\newtheorem{examples}[theorem]{Examples}
\newtheorem{question}[theorem]{Question}
\newtheorem{conjecture}[theorem]{Conjecture}

\newenvironment{proof}{{\it Proof:\/}}{$\Box$\vskip 0.08in}
\begin{document}
\centerline{ {\large \bf From 3-moves to Lagrangian tangles 
and cubic skein modules }}
\centerline{ \bf J\'ozef H. Przytycki }
\ \\
\centerline{{\bf Abstract}} 
{\footnotesize
We present an expanded version of four talks describing 
recent developments in Knot Theory to which the author 
contributed\footnote{Containing several results 
that are not yet published elsewhere.}. We discuss several 
open problems in classical Knot Theory and we develop techniques  
that allow us to study them:  Lagrangian tangles, skein modules 
and Burnside groups. The method of Burnside groups of links was 
discovered and developed only half a year after the last talk was 
delivered in Kananaskis\footnote{The first three talks were delivered 
at International Workshop on Graphs -- Operads -- Logic;
Cuautitl\'an, Mexico, March 12-16, 2001 and the fourth talk 
``Symplectic Structures on Colored Tangles" at the workshop
New Techniques in Topological Quantum Field Theory;
 Calgary/Kananaskis, August 22-27, 2001.} 
}.

\section{Open problems in Knot Theory that everyone can try 
to solve}

When did Knot Theory start? Was it in 1794 when 
C.~F. Gauss\footnote{Gauss'
notebooks contain several drawings related to knots, for example
a braid with complex coordinate description (see \cite{Ep-1, P-14}) or
the mysterious ``framed tangle" which is published here for the first
time, see Fig.1.2. \cite{Ga}.} (1777-1855) 
copied figures of knots from a book written in English (Fig.1.1)?
\\
\ \\
\centerline{\psfig{figure=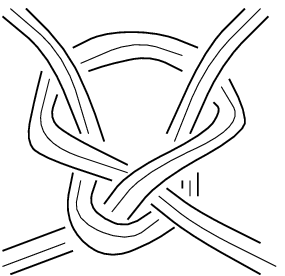,height=5.1cm}} 
\centerline{Fig. 1.1; \ Gauss' meshing knot from 1794}

\centerline{\psfig{figure=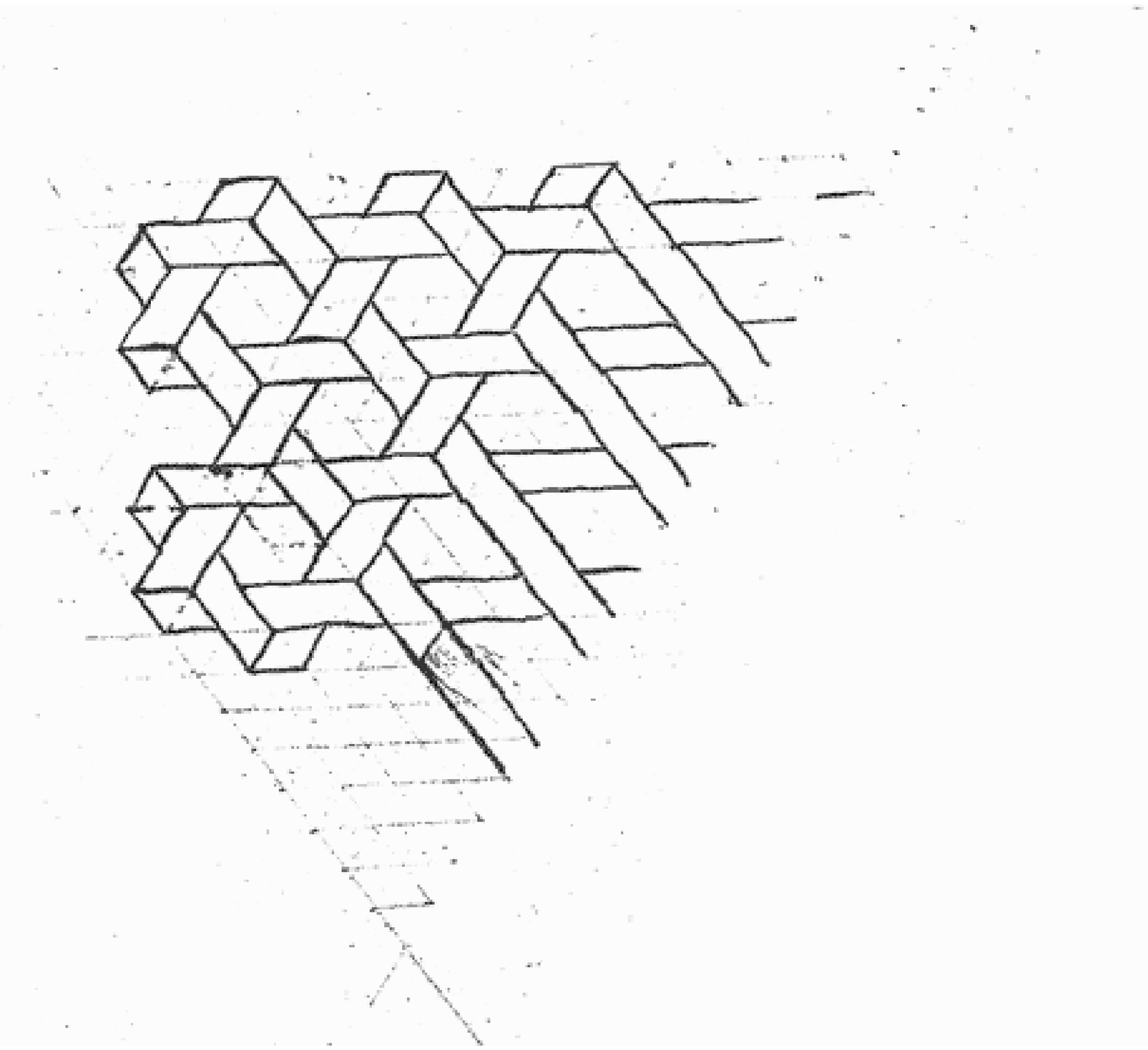,height=10.6cm}}
\centerline{Fig. 1.2; \ Framed tangle from Gauss'notebook \cite{Ga}}
\ \\

Or was it before that, in 1771,
when A-T. Vandermonde (1735-1796) considered knots and braids 
as a part of Leibniz's {\it analysis situs}?

\centerline{\psfig{figure=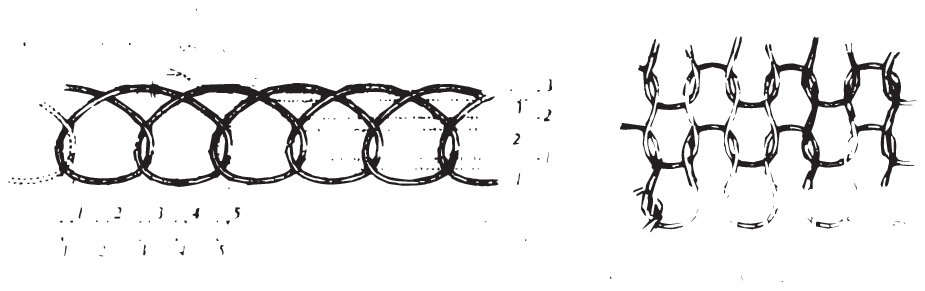,height=3.5cm}}
\centerline{Fig. 1.3; \ Vandermonde drawings of 1771}
\ \\

Perhaps engravings by Leonardo da Vinci\footnote{Giorgio Vasari writes 
in \cite{Va}: ``[Leonardo da Vinci] spent much time in making a regular
design of a series of knots so that the cord may be traced
from one end to the other, the whole filling a round space.
There is a fine engraving of this most difficult design,
and in the middle are the words: Leonardus Vinci Academia."} (1452-1519)
 \cite{Mac}  and woodcuts by Albrecht 
D\"urer\footnote{``Another great artist with whose works D\"urer now became
acquainted was Leonardo da Vinci. It does not seem likely that the
two artists ever met, but he may have been brought into relation
with him through Luca Pacioli, the author of the book De Divina
Proportione, which appeared at Venice in 1509, and an intimate friend
of the great Leonardo. D\"urer would naturally be deeply interested in
the proportion theories of Leonardo and Pacioli. He was certainly
acquainted with some engravings of Leonardo's school, representing a
curious circle of concentric scrollwork on a black ground, one of them
entitled Accademia Leonardi Vinci; for he himself executed six woodcuts
in imitation, the Six Knots, as he calls them himself. D\"urer was
amused by and interested in all scientific or mathematical problems..." 
From: http://www.cwru.edu/edocs/7/258.pdf, compare 
\cite{Dur-2}.} (1471-1528) \cite{Dur-1,Ha} should also be taken into
account, Fig.1.4.\\
\centerline{\psfig{figure=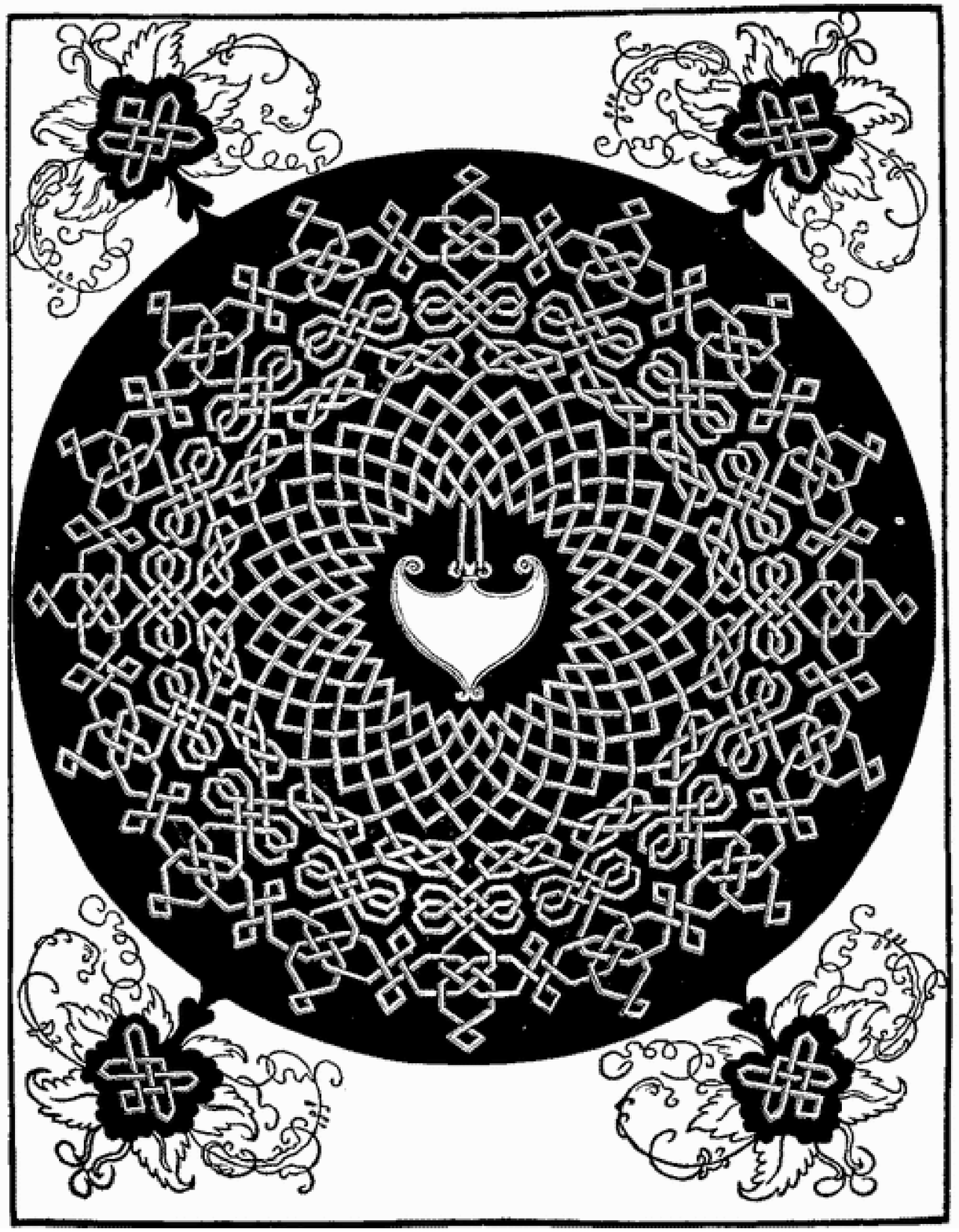,height=10.4cm}}
\centerline{Fig. 1.4; \ A knot by D\"urer \cite{Ku}; c. 1505-1507}

One can go back in time even further to ancient Greece
where surgeons considered sling knots, Fig.1.5  \cite{Da,P-14}.
\\
\ \\
\centerline{\psfig{figure=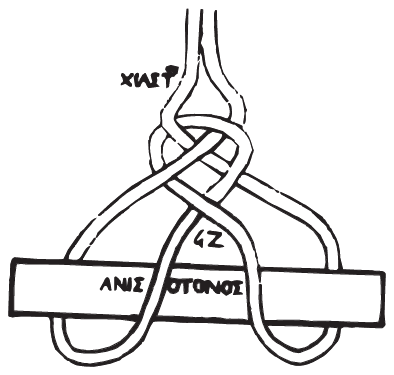,height=5.6cm}}\\
\ \\
\centerline{Fig. 1.5; \ A sling knot of Heraclas}

Moreover, we can appreciate ancient stamps and seals with knots
and links as their motifs.
The oldest examples that I am aware of are from  the pre-Hellenic Greece.
Excavations at Lerna
by the American School of Classical Studies under the direction of
Professor J.~L. Caskey (1952-1958) discovered two rich
deposits of clay seal-impressions.  The second deposit dated 
from about 2200 BC
contains several impressions of knots and links\footnote{The early Bronze
Age in Greece is divided,
as in Crete and the Cyclades, into three phases. The second phase lasted
from 2500 to 2200 BC, and was marked by a considerable increase in
prosperity. There were palaces at Lerna, and Tiryns, and probably elsewhere,
in contact with the Second City of Troy. The end
of this phase (in the Peloponnese)
was brought about by invasion and mass burnings.
The invaders are thought to be the first speakers of
the Greek language to arrive in Greece.} \cite{Hig,Hea,Wie}
 (see Fig.1.6). \\
\ \\
\ \\

\centerline{\psfig{figure=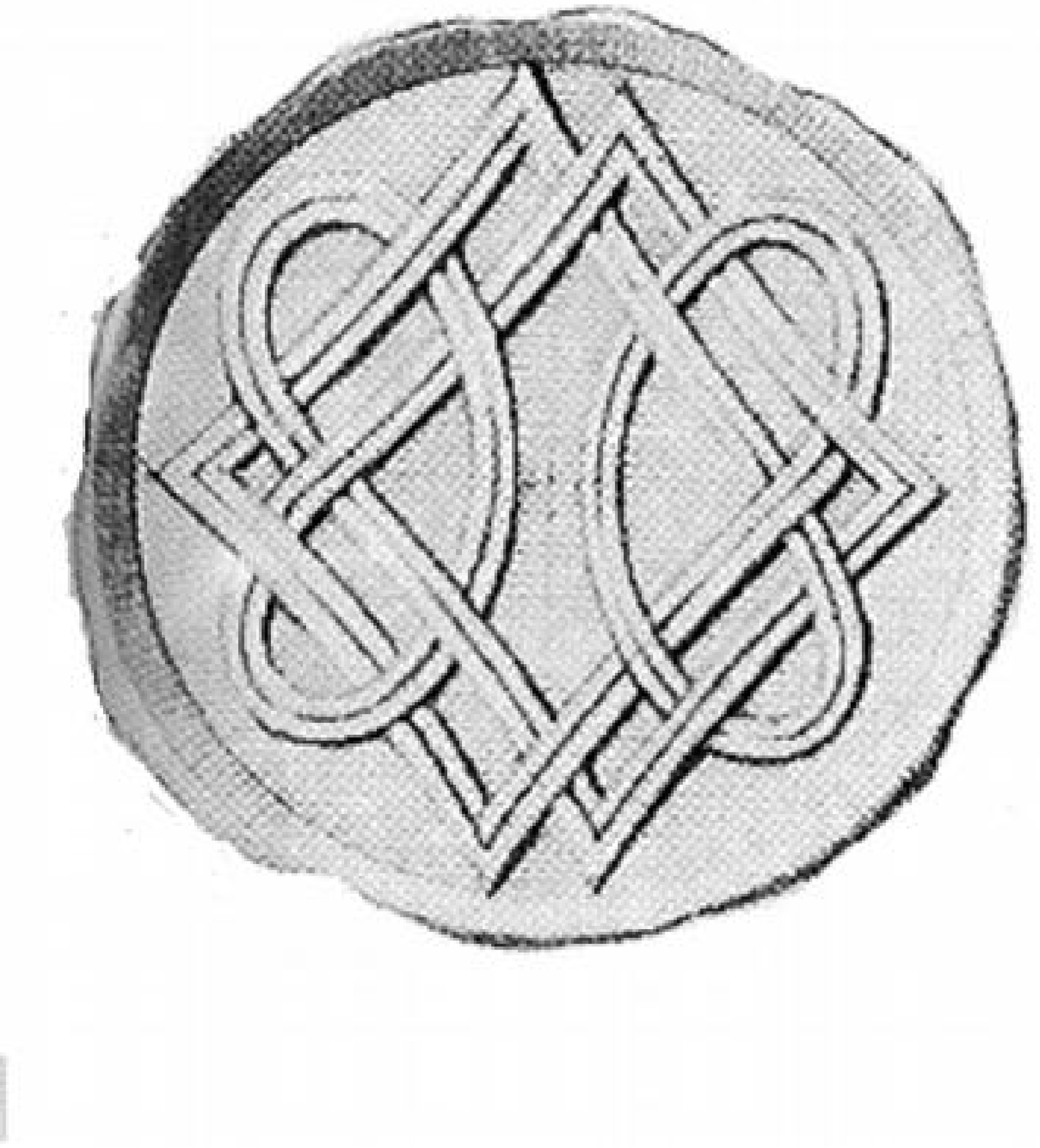,height=10.5cm}}
\centerline{Fig. 1.6; \ A seal-impression from the House of
the Tiles in Lerna \cite{Hig}.}
\ \\

As we see Knot Theory has a long history but despite this, or maybe
because of this, one still can find inspiring elementary open problems.
These problems are not just interesting puzzles but they lead
to an interesting theory.

 In this section our presentation is absolutely elementary.
Links are circles embedded in our space, $R^3$, 
up to topological deformation,
that is, two links are equivalent if one can be deformed into the other 
in space without cutting and pasting. We represent links using
their plane diagrams. 

First we introduce the concept of an $n$ move on a link. 
\begin{definition}\label{1.1}
An \textup{$n$-move} on a link is a local transformation of the 
link illustrated in Figure 1.7.
\end{definition}
\centerline{\psfig{figure=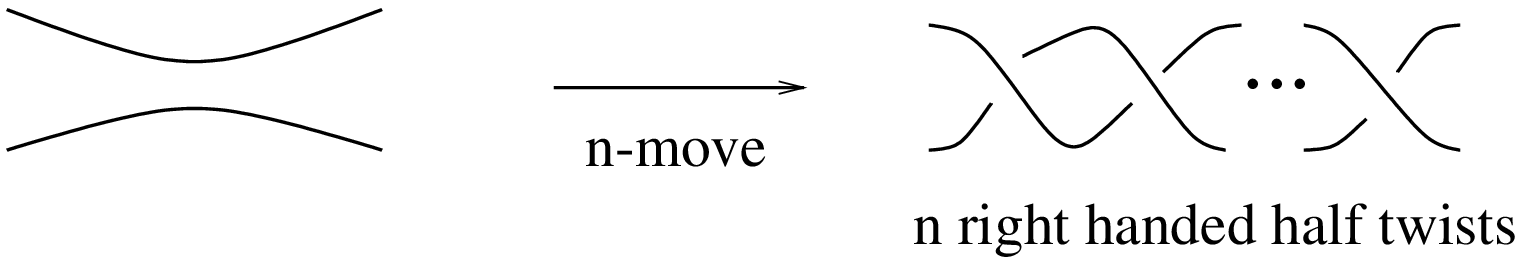,height=1.6cm}}
\centerline{Fig. 1.7; \ $n$-move}
\ \\

In our convention, the part of the link outside of the disk 
in which the move takes place, remains unchanged. 
One should stress that an $n$-move can change the topology of the link.
For example 
\parbox{2.9cm}{\psfig{figure=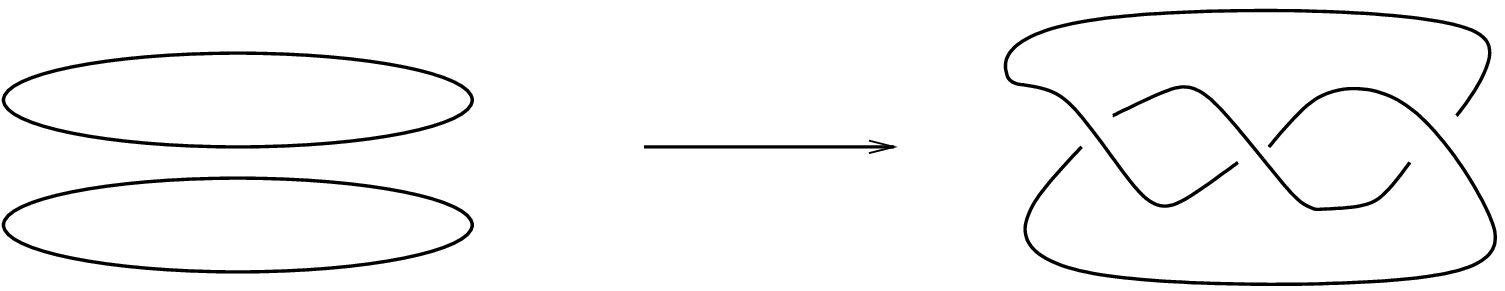,height=0.5cm}} 
illustrates a 3-move.
\begin{definition}\label{1.2}
We say that two links, $L_1$ and  $L_2$,  are 
\textup{$n$-move equivalent} if one can 
obtain $L_2$ from $L_1$  by a finite number of $n$-moves 
and $(-n)$-moves (inverses of $n$-moves). 
\end{definition}

If we work with diagrams of links then the topological deformation
 of links is captured by Reidemeister moves, 
that is, two diagrams represent the same link in space if and 
only if one can obtain one of them from the other 
by a sequence of Reidemeister moves:\ \\
\ \\
\centerline{\psfig{figure=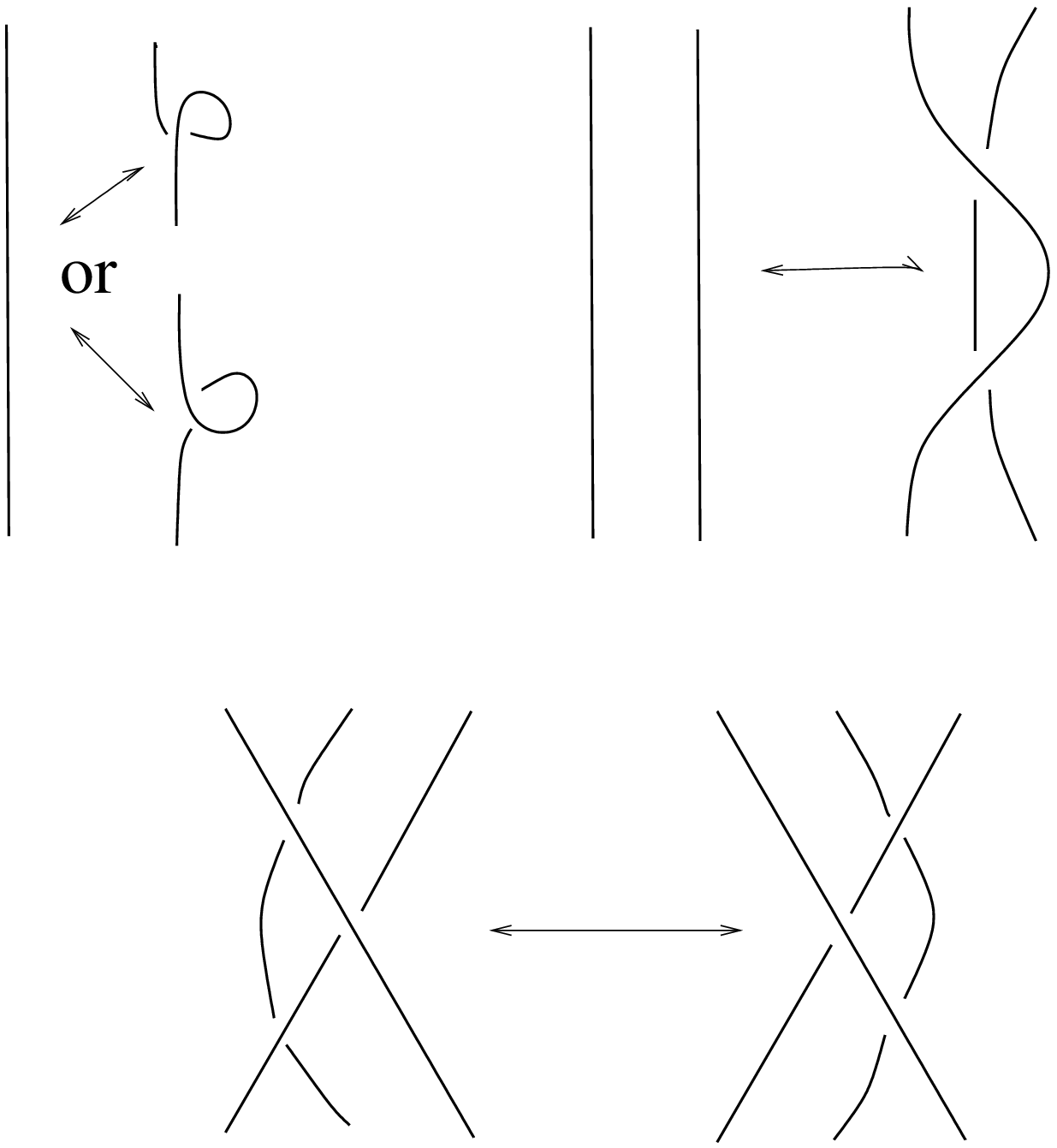,height=5.8cm}}
\centerline{Fig. 1.8; \ Reidemeister moves}
\ \\

Thus, we say that two diagrams, $D_1$ and $D_2$, are $n$-move equivalent
if one can be  obtained from the other by a sequence of
 $n$-moves, their inverses and Reidemeister moves. 
To illustrate this, we show that the move
\parbox{2.3cm}{\psfig{figure=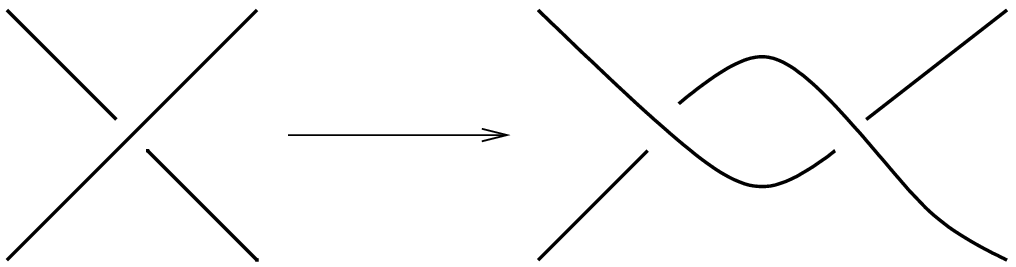,height=0.6cm}}    
is the result of an application of a 3-move followed by the second 
Reidemeister move (Fig.1.9).
\\
\ \\
\centerline{\psfig{figure=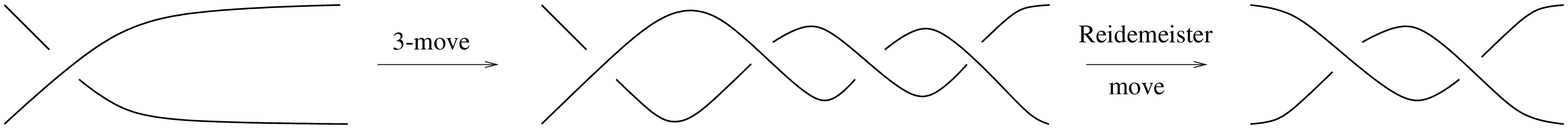,height=1.2cm}}
\centerline{Fig. 1.9}

\begin{conjecture} [Montesinos-Nakanishi]
\label{Conjecture 1.3}\ \\
 Every link is $3$-move equivalent to a trivial link.
\end{conjecture}

Yasutaka Nakanishi proposed this conjecture in 1981. 
Jos\'e Montesinos analyzed 3-moves before, in connection with 
3-fold dihedral branch coverings, and asked a related but 
different question\footnote{``Is there a set of moves which do not 
change the covering manifold and
such that if two colored links have the same covering they
are related by a finite sequence of those moves?"} \cite{Mo-2}.

\begin{examples}\label{1.4}
\begin{enumerate}
\item [\textup{(i)}] 
Trefoil knots (left- and right-handed) are 3-move equivalent 
to the trivial link of two components:\\
\ \\ 
\centerline{\psfig{figure=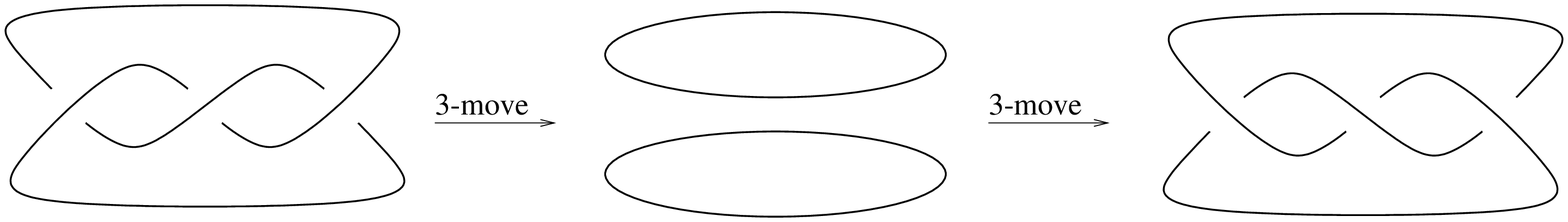,height=1.7cm}}
\centerline{\textup{Fig. 1.10}} 
\item [\textup{(ii)}] 
The figure eight knot $(4_1)$ and the knot $5_2$ are 
3-move equivalent to the trivial knot:\\
\ \\ 
\centerline{\psfig{figure=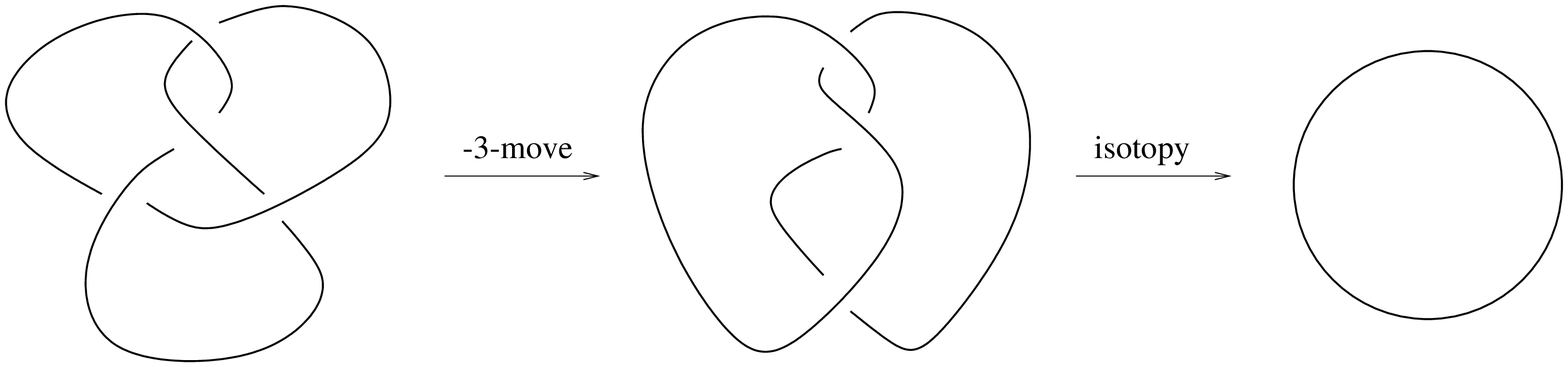,height=1.7cm}\ \ \ 
\psfig{figure=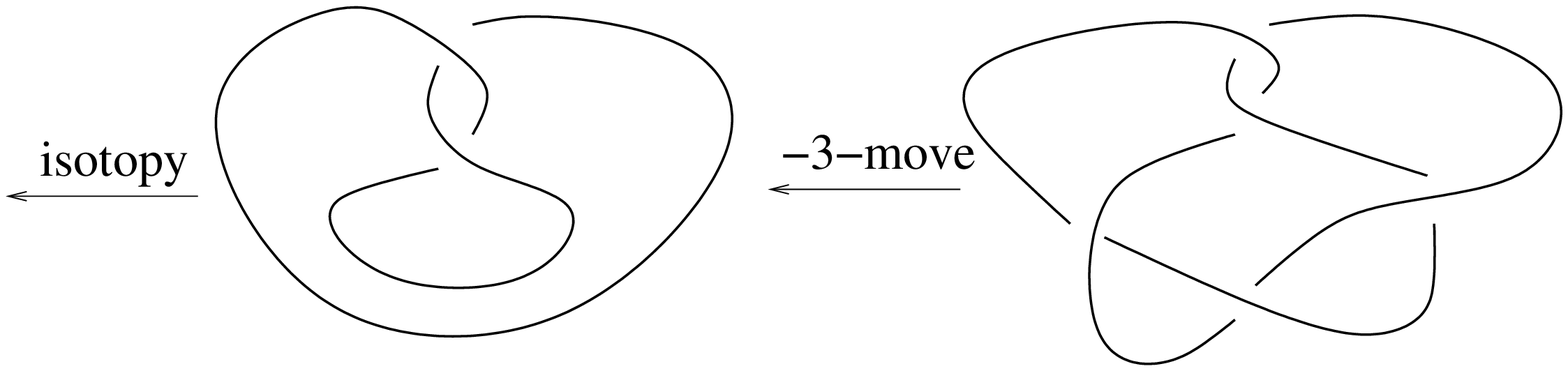,height=1.7cm}}
\centerline{\textup{Fig. 1.11}}
\item [\textup{(iii)}] The knot $5_1$ and the Hopf link are
3-move equivalent to the trivial knot\footnote{One can show that 
the knot $5_1$ cannot be reduced to the trivial knot by one ${\pm 3}$-move. 
To see this one can use the Goeritz matrix approach to the classical 
signature ($|\sigma(5_1)|=4$), see \cite{Go,G-L}.}:\\
\ \\
\centerline{\psfig{figure=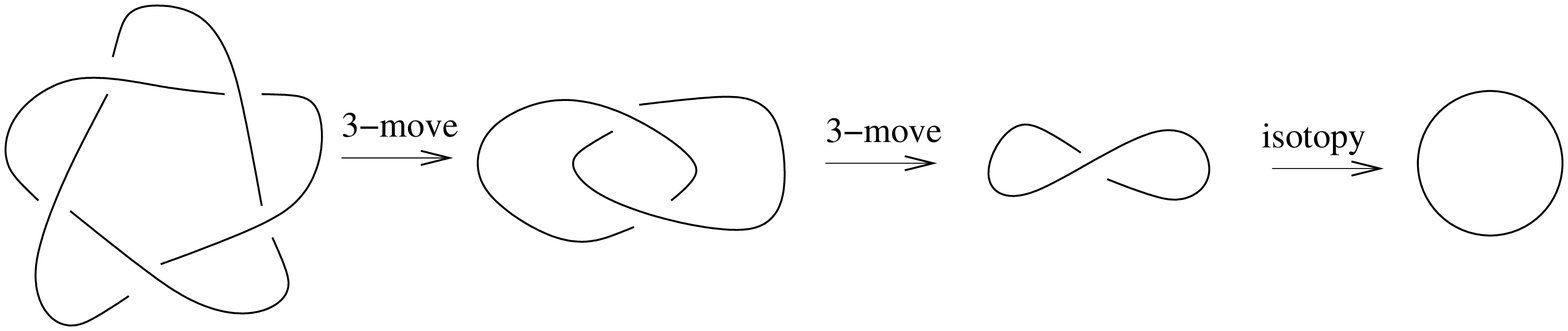,height=2.5cm}}
\centerline{\textup{Fig. 1.12}}

\end{enumerate}
\end{examples}

We will show later, in this section, that different trivial links
are not 3-move equivalent. However, in order to achieve this
conclusion we need an invariant of links preserved by 3-moves and
differentiating trivial links. 
Fox 3-coloring is such an invariant. 
We will introduce it later (today in its simplest
form and, in the second lecture, in a more general context of Fox
$n$-colorings and Alexander-Burau-Fox colorings).
 
Now let us present some other related conjectures.
\begin{conjecture}\label{Conjecture 1.5}\ \\
 Any $2$-tangle is $3$-move equivalent to one of the four 2-tangles
shown in Fig.1.13 with possible additional trivial components.
\end{conjecture}
\ \\
\centerline{\psfig{figure=S+S-S0Sinf.eps}}
\begin{center}
Fig. 1.13
\end{center}

The Montesinos-Nakanishi conjecture follows from Conjecture
1.5. More generally if Conjecture 1.5 holds for some class of 2-tangles,
then Conjecture 1.3 holds for any link obtained by closing elements
of this class, without introducing any new
crossing. The simplest interesting class of tangles for which
Conjecture 1.5 holds  are algebraic tangles in the sense of Conway
(I call them 2-algebraic tangles and in the next section 
 present a generalization). Conjecture 1.5 can be proved by induction
for 2-algebraic tangles. I will leave the proof to you 
as a pleasant exercise (compare Proposition 1.9). 
The definition you need is as follows

\begin{definition}[\cite{Co,B-S}]\label{1.6}\ 
The family of \textup{2-algebraic tangles} is the smallest family of 
2-tangles satisfying\\
\textup{(i)} Any 2-tangle with $0$ or $1$ crossing is $2$-algebraic.\\
\textup{(ii)} If $A$ and $B$ are 2-algebraic tangles
    then the $2$-tangle $r^i(A)*r^j(B)$ is also
   2-algebraic, where $r$ denotes 
the counterclockwise rotation of a tangle by $90^o$  along 
the $z$-axis, and $*$ denotes the horizontal composition of 
tangles (see the figure below). \\
\ \\
\centerline{\psfig{figure=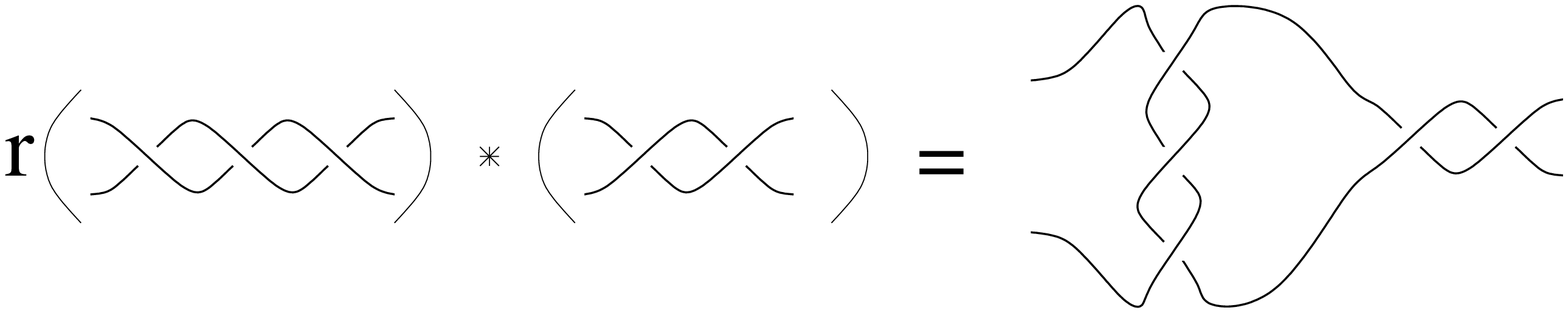,height=2.6cm}}

A link is called \textup{2-algebraic} if it can be 
obtained from a 2-algebraic tangle
by closing its ends without introducing crossings\footnote{By 
joining the top ends and then bottom ends of a tangle $T$ one obtains the 
link $N(T)$, the {\it numerator} of $T$, Fig.1.22, 1.23.  Joining 
the left-hand ends and then right-hand ends produces the 
{\it denominator} closure  $D(T)$.}.
\end{definition}

The Montesinos-Nakanishi 3-move conjecture has been proved
by my students Qi Chen and Tatsuya Tsukamoto for many special
families of links  \cite{Che,Tsu,P-Ts}.
In particular, Chen proved that the conjecture holds for all
$5$-braid links except possibly one family, containing
the square of the center of the 5-braid group, 
$\Delta_5^4=
(\sigma_1\sigma_2\sigma_3\sigma_4)^{10}$. He also found a reduction 
by $\pm 3$-moves of $\Delta_5^4$ to the 5-braid link, 
$(\sigma_1^{-1}\sigma_2\sigma_3\sigma_4^{-1}\sigma_3)^4$, with 
$20$ crossings\footnote{In the group 
$B_5/(\sigma_i^3)$ 
the calculation is as follows: $(\sigma_1\sigma_2\sigma_3\sigma_4)^{10}= 
(\sigma_1\sigma_2\sigma_3\sigma_4^2\sigma_3\sigma_2\sigma_1)^2
(\sigma_2\sigma_3\sigma_4^2\sigma_3\sigma_2)^2
(\sigma_3\sigma_4^2\sigma_3)^2 \stackrel{3}{=}
(\sigma_1\sigma_2\sigma_3\sigma_4^2\sigma_3\sigma_2\sigma_1)^2
(\sigma_2\sigma_3\sigma_4^2\sigma_3\sigma_2)^2 \stackrel{3}{=}
(\sigma_1\sigma_2\sigma_3\sigma_4^2\sigma_3\sigma_2\sigma_1)^2
(\sigma_2^{-1}\sigma_3\sigma_4^{-1}\sigma_3)^2 \stackrel{3}{=}
(\sigma_1\sigma_2\sigma_3\sigma_4^{-1}\sigma_3\sigma_2\sigma_1
\sigma_2^{-1}\sigma_3\sigma_4^{-1}\sigma_3)^2 = 
\sigma_1\sigma_2\sigma_3\sigma_4^{-1}\sigma_3
\sigma_1^{-1}\sigma_2\sigma_3\sigma_4^{-1}\sigma_3\sigma_1)^2  
\stackrel{3}{=} 
(\sigma_1^{-1}\sigma_2\sigma_3\sigma_4^{-1}\sigma_3)^4$.}, Fig.1.14.
It is now the smallest known possible counterexample 
to the Montesinos-Nakanishi 3-move 
conjecture\footnote{We proved in \cite{D-P-1} that 
Chen's link is in fact the counterexample to the Montesinos-Nakanishi 3-move 
conjecture; see Section 4. We think that it is the smallest such 
counterexample. We also demonstrated that the 2-parallel of the
Borromean rings is not 3-move equivalent to a trivial link. It is still 
possible that Chen's link with an additional trivial component is 
3-move equivalent to the 2-parallel of the Borromean rings.}.  

\ \\
\centerline{\psfig{figure=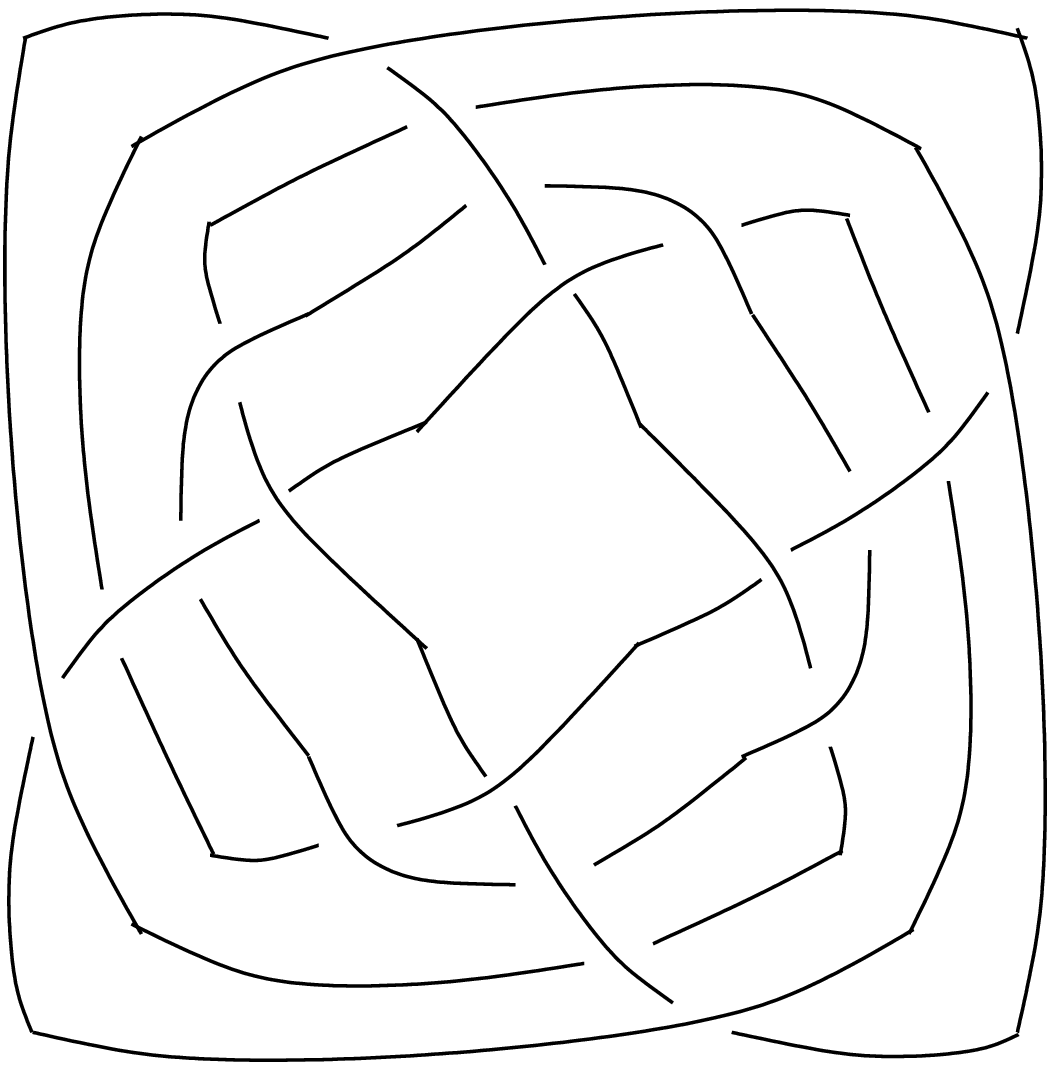,height=6.1cm}}
\begin{center}
Fig. 1.14
\end{center}
Previously Nakanishi suggested in 1994 (see \cite{Kir}), 
the 2-parallel of the Borromean rings (a 6-braid with 24 crossings) 
as a possible counterexample (Fig.1.15).

\ \\
\centerline{\psfig{figure=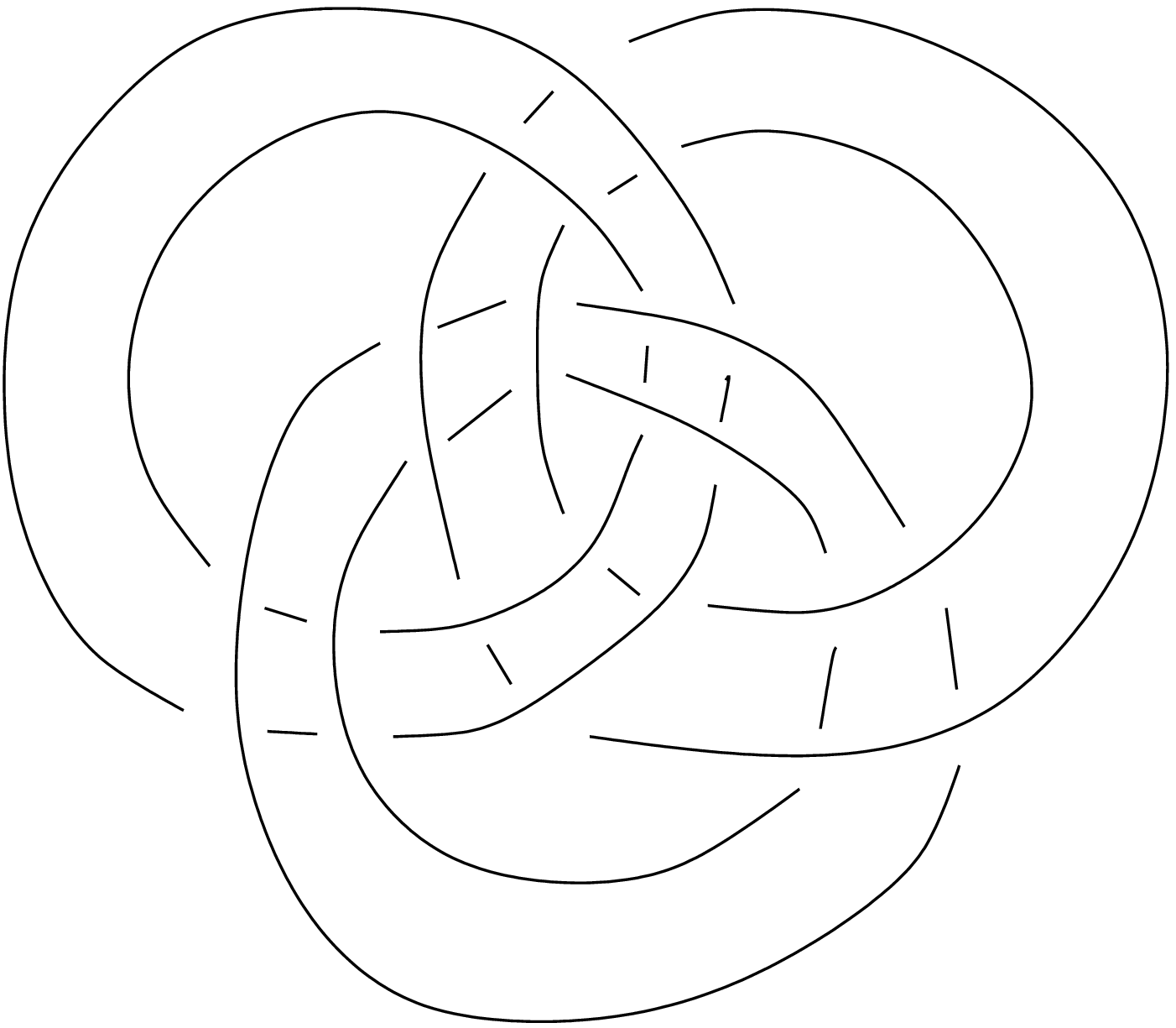,height=5.1cm}}
\begin{center}
Fig. 1.15
\end{center}

We will return to the discussion of theories motivated by 3-moves 
tomorrow. 
Now we will state some conjectures that employs other elementary moves.

\begin{conjecture}[Nakanishi, 1979]\label{1.7}\ \\
Every knot is $4$-move equivalent to the trivial knot.
\end{conjecture} 
\begin{examples}\label{1.8}
Reduction of the trefoil and the figure eight knots is 
illustrated in Fig.1.16.
\end{examples}
\centerline{\psfig{figure=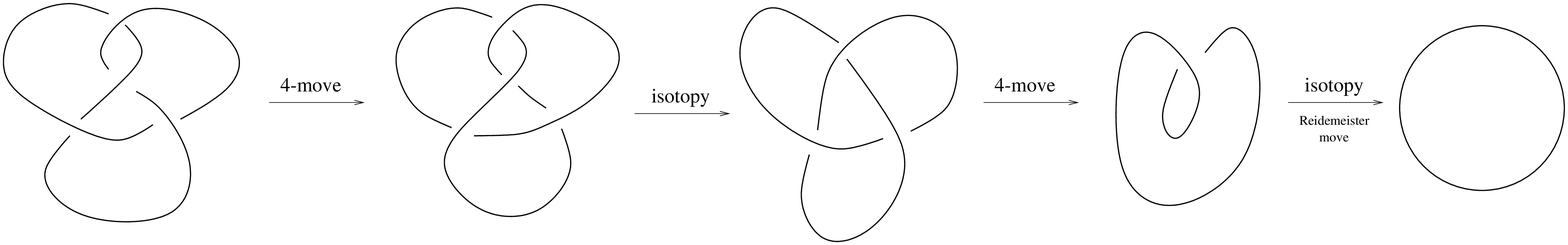,height=2.4cm}}
\centerline{Fig. 1.16}

\begin{proposition}[\cite{P-12}]\label{1.9}
\begin{enumerate}
\item[\textup{(i)}]
Every 2-algebraic tangle without a closed component can
be reduced by $\pm 4$-moves to one of the six basic 2-tangles
shown in Fig.1.17.
\item[\textup{(ii)}]
Every 2-algebraic knot can be reduced by $\pm 4$-moves to the trivial knot.
\end{enumerate}
\end{proposition}
\ \\
\centerline{\psfig{figure=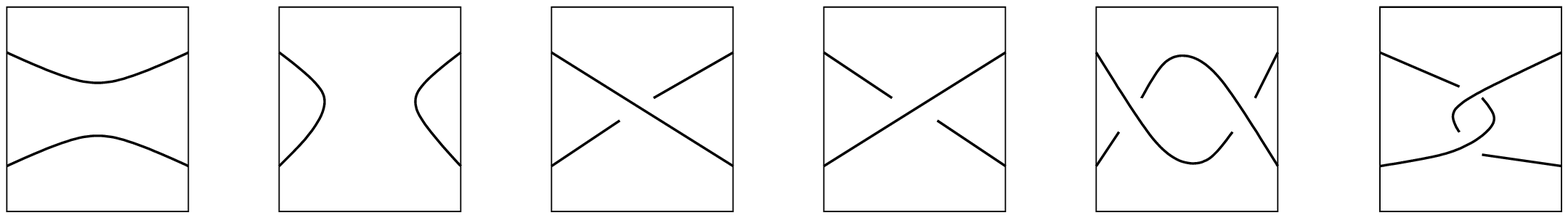,height=1.6cm}}
\centerline{\ \ \  $e_1$\ \ \ \ \ \ \ \ \ \ \ \  
 $e_2$ \ \ \ \ \ \ \ \ \ \ \ \ 
 $e_3$\ \ \ \ \ \ \ \ \ \ \ \ \  
$e_4$\ \ \ \ \ \ \ \ \ \ \ \ \ \ 
$e_5$\ \ \  \ \ \ \ \ \ \ \ \ \ $e_6$\ \ \  }
\centerline{Fig. 1.17}
\ \\

\begin{proof}
To prove (i) it suffices to show that every composition (with possible rotation)
of tangles presented in Fig.1.17 can be reduced by $\pm 4$-moves back
to one of the tangles in Fig.1.17 or it has a closed component.
These can be easily verified by inspection. Fig.1.18 is the 
multiplication table for basic tangles. We have chosen our basic tangles 
to be invariant under the rotation $r$, so it suffices to be able to reduce 
every tangle of the table to a basic tangle.
One example of such a reduction is shown in Fig.1.19. \
Part (ii) follows from (i).
\end{proof}
\ \\
\centerline{\psfig{figure=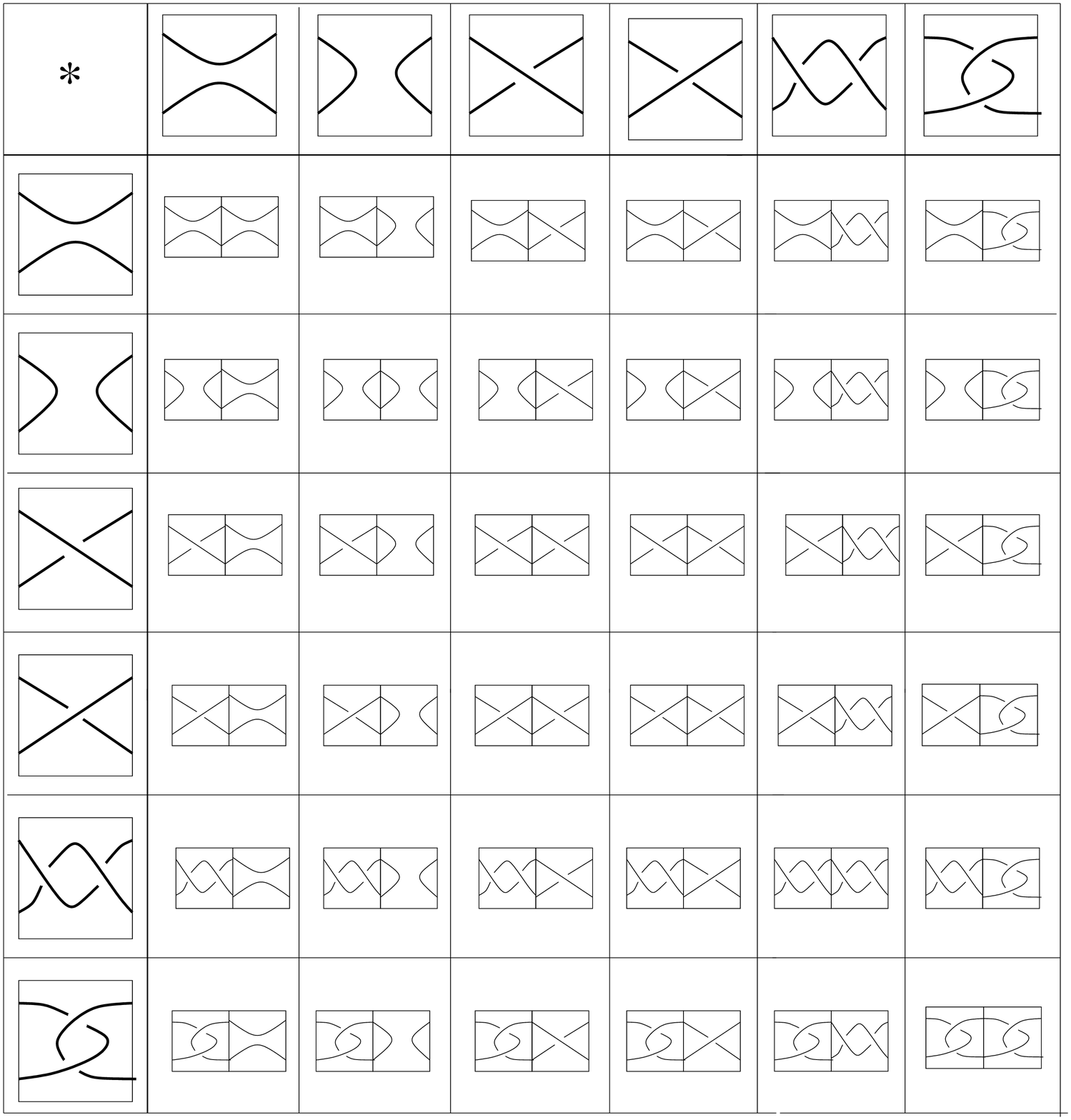,height=10.9cm}}
\centerline{Fig. 1.18}

\ \\
\centerline{\psfig{figure=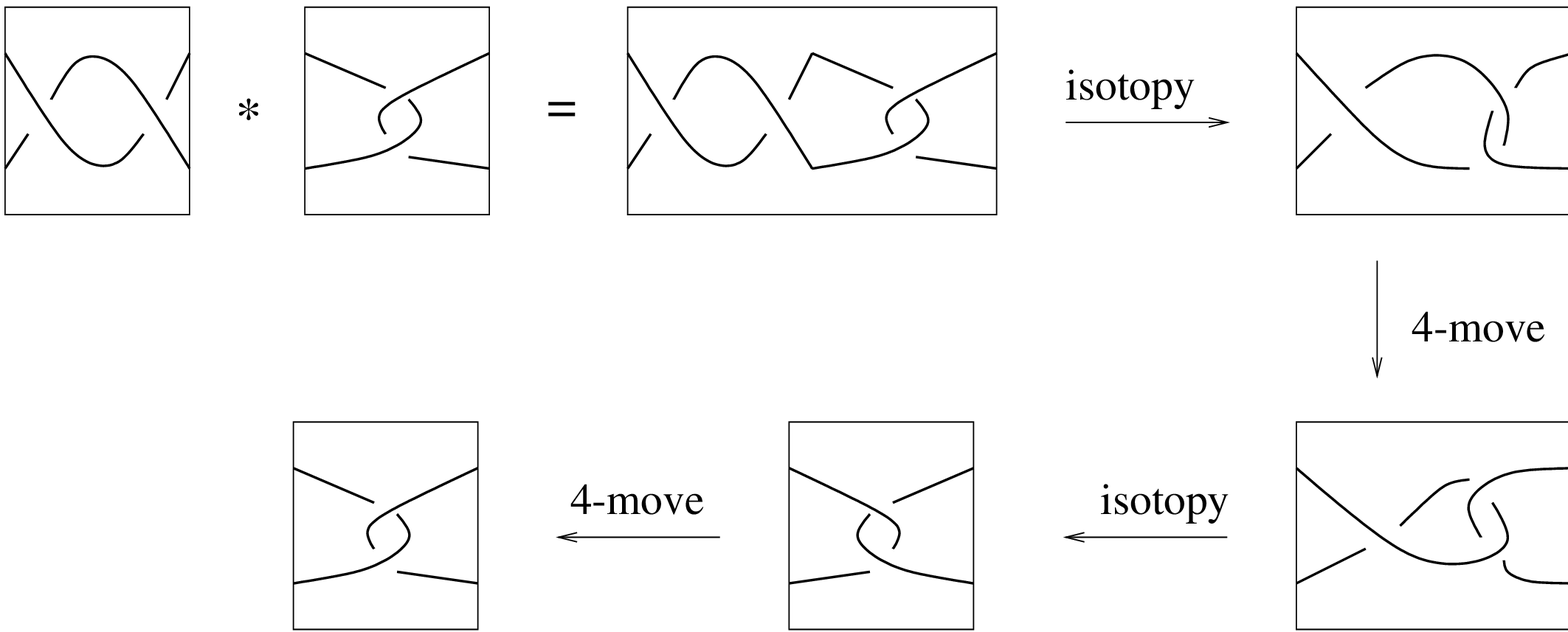,height=3.6cm}}
\centerline{Fig. 1.19}
\ \\

In 1994, Nakanishi began to suspect that a $2$-cable
of the trefoil knot cannot be simplified by 4-moves \cite{Kir}.
However, Nikolaos Askitas was able to simplify this knot \cite{Ask}.
Askitas, in turn,
suspects that the $(2,1)$-cable of the figure eight knot 
(with $17$ crossings) to be the simplest counterexample to the 
Nakanishi $4$-move conjecture.

Not every link can be reduced to a trivial link by 4-moves.
In particular, the linking matrix modulo 2 is preserved by $4$-moves.
Furthermore, Nakanishi and Suzuki demonstrated that the
Borromean rings cannot be reduced to the trivial link
of three components \cite{Na-Su}.

In 1985, after the seminar talk given by Nakanishi in Osaka, 
there was discussion about possible 
generalization of the Nakanishi 4-move conjecture for links. 
Akio Kawauchi formulated the following question for links
\begin{problem}[\cite{Kir}]\label{1.10}\
\begin{enumerate}
\item[\textup{(i)}]
Is it true
that if two links are link-homotopic\footnote{Two links $L_1$ and $L_2$ are 
{\it link-homotopic} if one can obtain $L_2$ from $L_1$ by a finite 
number of crossing changes involving only self-crossings of the components.}
then they are 4-move equivalent?
\item[\textup{(ii)}] In particular, is it true
that every two component link is 4-move equivalent to the trivial link
of two components or to the Hopf link?
\end{enumerate}
\end{problem}
We can extend the  argument used in Proposition 1.9 to show:
\begin{theorem}\label{1.11} 
Any two component 2-algebraic link is 4-move equivalent
to the trivial link of two components or to the Hopf link.
\end{theorem}
\begin{proof} 
Let $L$ be a 2-algebraic link of two components. 
Therefore, $L$ is built inductively as in Definition 1.6.  Consider the first
tangle, $T$, in the construction, which has a closed component 
(if it happens). The complement $T'$ of $T$ in the link $L$ is also
a 2-algebraic tangle but without a closed component. Therefore it 
can be reduced to one of the $6$ basic tangles shown in
 Fig.1.17, say $e_i$.
Consider the product $T*e_i$. The only nontrivial tangle $T$ to be
 considered is $e_6*e_6$ (the last tangle in Fig.1.18).
The compositions $(e_6*e_6)*e_i$ are illustrated in Fig.1.20. 
The closure of each of these product tangles (the numerator 
or the denominator)
 has two components because it is $4$-move equivalent to $L$. 
We can easily check that it reduces to the trivial 
link of two components.
\end{proof}
\ \\

\centerline{\psfig{figure=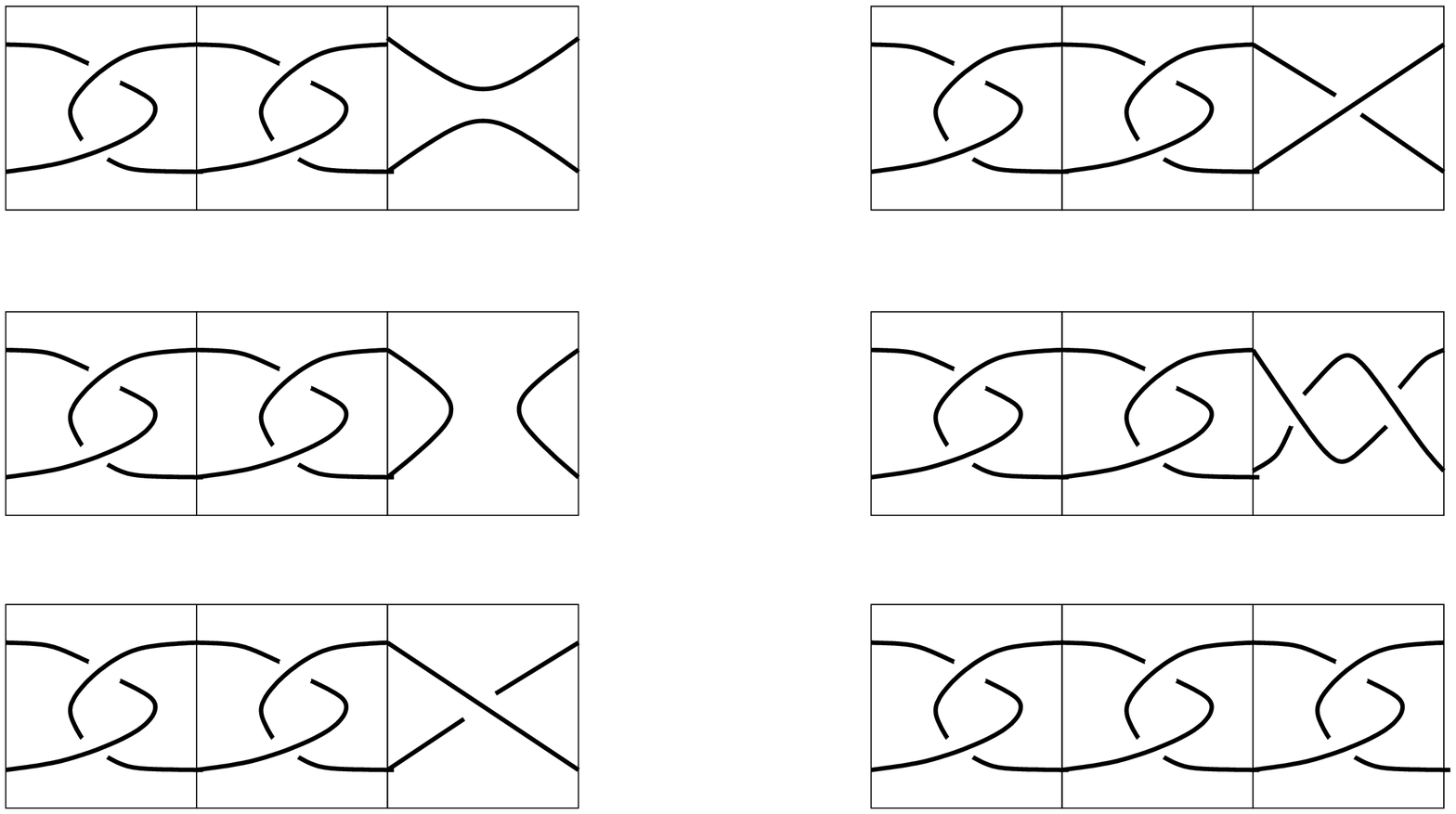,height=4.8cm}}
\centerline{Fig. 1.20}

\begin{problem}\label{1.12}
\begin{enumerate}
\item[\textup{(i)}]
Find a (reasonably small) family of 2-tangles with one closed 
component so that every 2-tangle with one closed component 
is $4$-move equivalent to one of its elements.
\item[\textup{(ii)}] 
Solve the above problem for 2-algebraic tangles with
one closed component.
\end{enumerate}
\end{problem}

Nakanishi \cite{Nak-2} pointed out that  the ``half"
2-cabling of the Whitehead link, ${\cal W}$, Fig.1.21,
was the simplest link which he could not reduce to a trivial link
 by $\pm 4$-moves but which was link-homotopic to a trivial 
link\footnote{In fact, in  June of 2002
we showed that this example cannot be reduced by $\pm 4$-moves 
 \cite{D-P-2}.}. 
\ \\
\ \\
\centerline{\psfig{figure=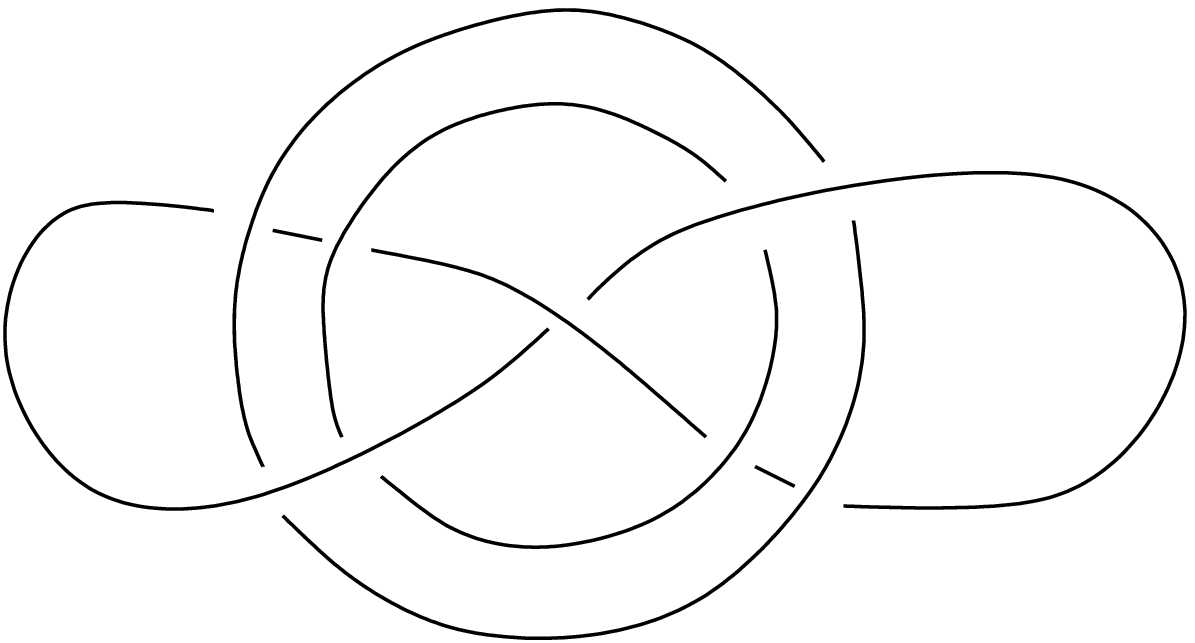,height=4.8cm}}
\centerline{Fig. 1.21}\ \\
\ \\
It is shown in Fig.1.22 that the link ${\cal W}$ is 2-algebraic. 
Similarly, the Borromean rings, $BR$, are 2-algebraic, 
Fig.1.23 (compare Fig.1.30).   \\
\ \\
\centerline{\psfig{figure=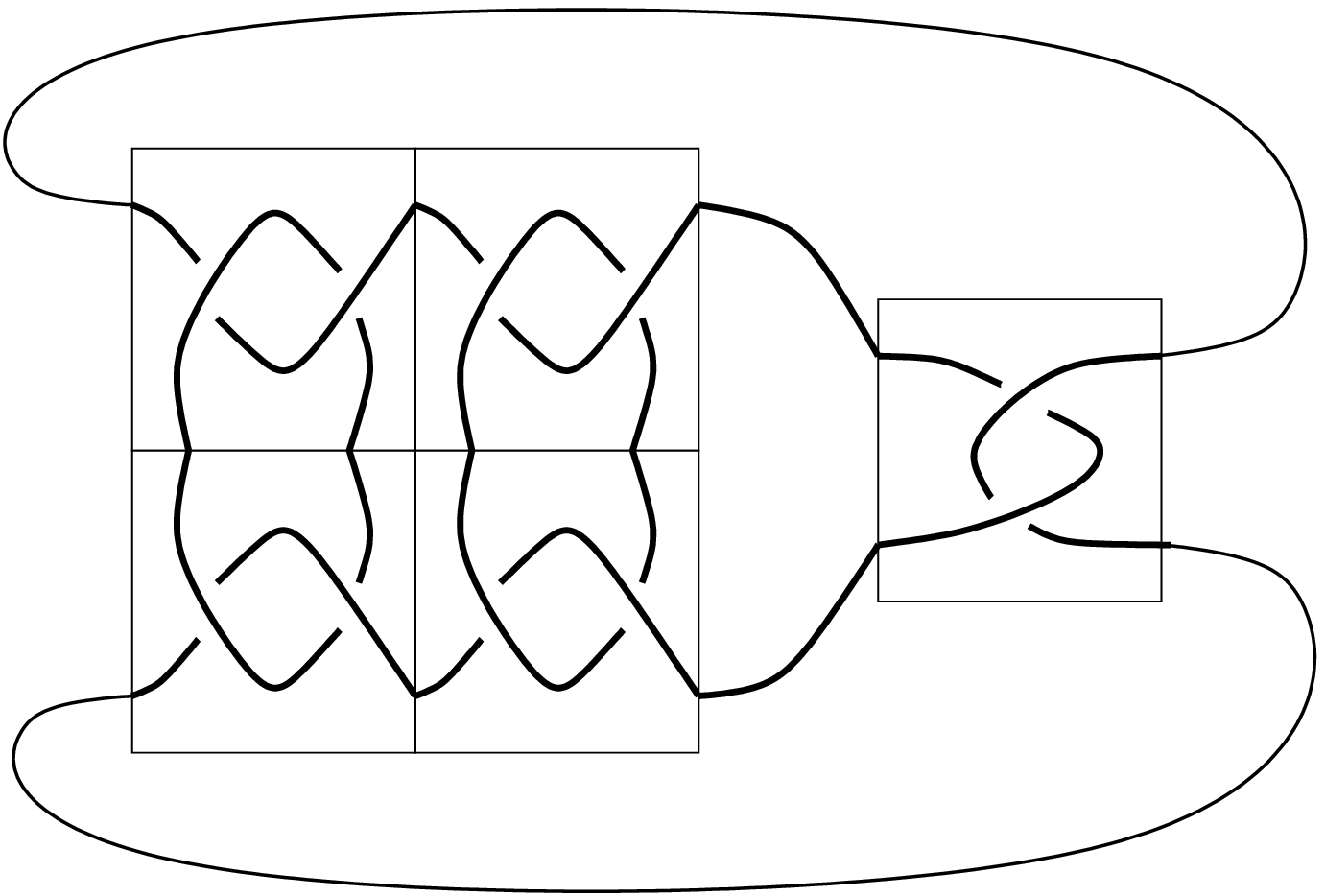,height=4.8cm}}
\centerline{Fig. 1.22; ${\cal W}=
N(r(r(e_3*e_3)*r(e_4*e_4))*r(r(e_3*e_3)*r(e_4*e_4))*r(e_3*e_3))$
}

\ \\
\centerline{\psfig{figure=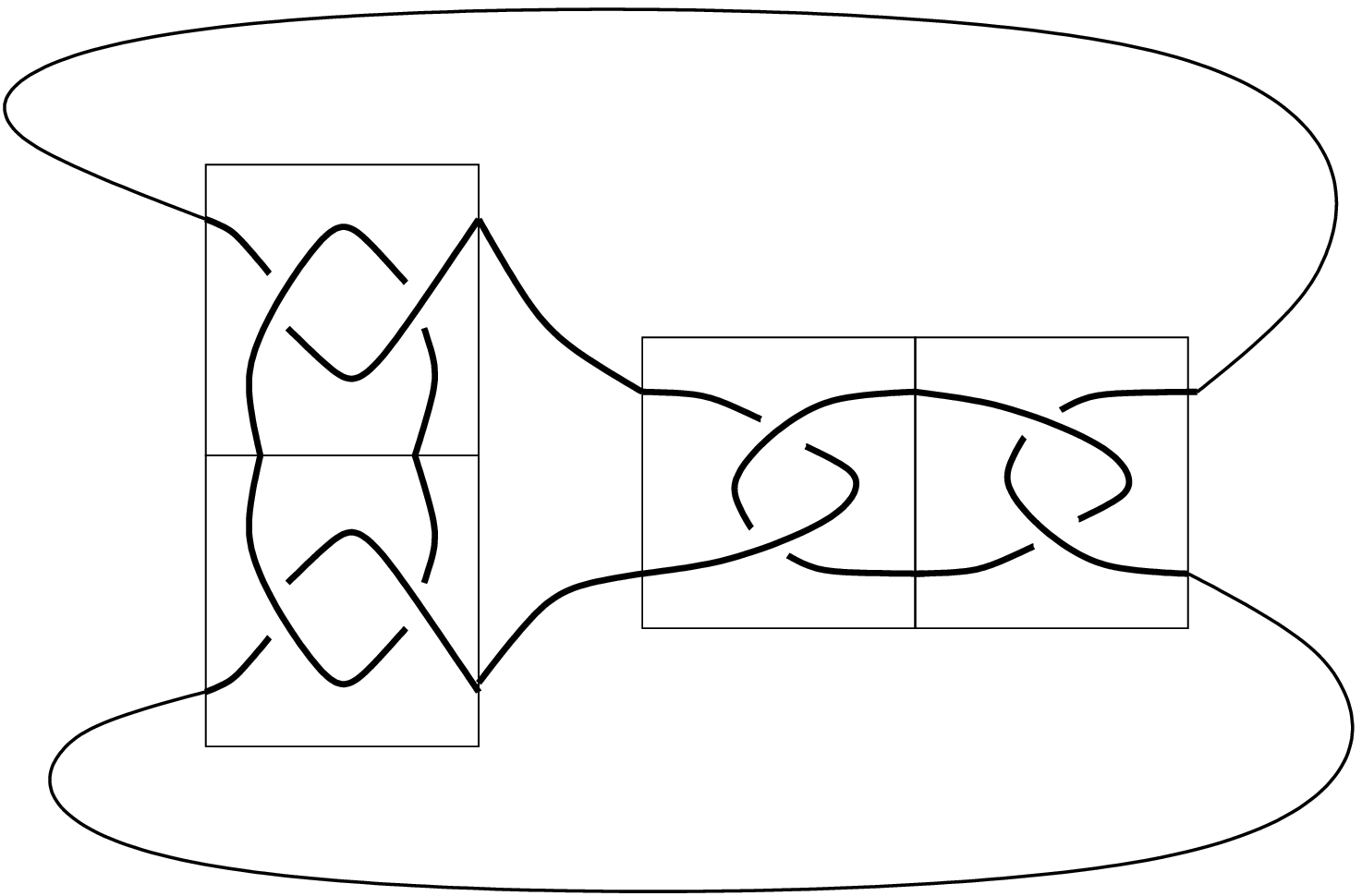,height=4.8cm}}
\centerline{Fig. 1.23;
$BR= N(r(r(e_3*e_3)*r(e_4*e_4))*r(e_3*e_3)*r(e_4*e_4))$}

\begin{problem}\label{1.13}
Is the link ${\cal W}$ the only 2-algebraic link of three components (up 
to 4-move equivalence) which 
is homotopically trivial but which is not 4-move equivalent to the 
trivial link of three components?
\end{problem}

We can also prove that the answer to Kawauchi's question is affirmative
for closed 3-braids.
\begin{theorem}\label{1.14}
\begin{enumerate}
\item[\textup{(i)}] Every knot which is a closed 3-braid is 
4-move equivalent to the trivial knot.
\item[\textup{(ii)}] Every link of two components which is
a closed 3-braid is 4-move equivalent to the trivial link of two
 components or to the Hopf link.
\item[\textup{(iii)}] 
Every link of three components which is a closed
3-braid is $4$-move equivalent either to the trivial link of three
 components, or to the Hopf
link with the additional trivial component, or to  the connected sum of two
Hopf links, or to the $(3,3)$-torus link, $\bar 6^3_{1}$, represented 
by $(\sigma_1\sigma_2)^3$  (all linking numbers are equal to 1),
or to the Borromean rings (represented by $(\sigma_1\sigma_2^{-1})^3$).
\end{enumerate}
\end{theorem}
\begin{proof}
Our proof is based on
the Coxeter theorem that the quotient group $B_3/(\sigma_i^4)$ is 
finite with $96$ elements, \cite{Cox}.
Furthermore, $B_3/(\sigma_i^4)$ has 16 conjugacy classes\footnote{$Id,\sigma_1,
\sigma_1^{-1},\sigma_1^{2}, \sigma_1\sigma_2,\sigma_1^{-1}\sigma_2,
\sigma_1^{-1}\sigma_2^{-1},
\sigma_1^{2}\sigma_2, \sigma_1^{2}\sigma_2^{-1}, \sigma_1^{2}\sigma_2^{2},
\sigma_1\sigma_2^{-1}\sigma_1\sigma_2^{-1}, 
\sigma_1\sigma_2^{2}\sigma_1\sigma_2^{-1}, 
\sigma_1\sigma_2^{-1}\sigma_1^{2}\sigma_2^{-1},\\
\sigma_1\sigma_2^{2}\sigma_1\sigma_2^{2},
\sigma_1^{-1}\sigma_2^{2}\sigma_1^{-1}\sigma_2^{2},
(\sigma_1\sigma_2^{-1})^3$ \ (checked by M.~D{\c a}bkowski
using the {\it GAP} program).}:
9 of them can be  easily identified as representing trivial links (up to 
$4$-move equivalence), and
2 of them represent the Hopf link 
($\sigma_1^{2}\sigma_2$ and  $\sigma_1^{2}\sigma_2^{-1}$), 
and $\sigma_1^{2}$ represents the Hopf link with an additional 
trivial component. 
We also have the connected sums of
Hopf links ($\sigma_1^{2}\sigma_2^{2}$). 
Finally, we are left with two representatives of the
link $\bar 6^3_{1}$ ($\sigma_1\sigma_2^{2}\sigma_1\sigma_2^{2}$ and 
$\sigma_1^{-1}\sigma_2^{2}\sigma_1^{-1}\sigma_2^{2}$)
 and the Borromean rings.
\end{proof}

Proposition 1.9 and Theorems 1.11, and 1.14 can be used to analyze 
$4$-move equivalence classes of links with small number of crossings.
\begin{theorem}\label{1.15}
\begin{enumerate}
\item[\textup{(i)}] 
Every knot of no more than $9$ crossings is $4$ move equivalent 
to the trivial knot.
\item[\textup{(ii)}] Every two component link of no more than $9$ crossings 
is 4-move equivalent to the trivial link of two components or to 
the Hopf link.
\end{enumerate}
\end{theorem}
\begin{proof} 
Part (ii) follows immediately as the only 2-component links 
with up to 9 crossings which are not 2-algebraic are
$9^2_{40}, 9^2_{41}, 9^2_{42}$ and $9^2_{61}$ and all these links
are closed 3-braids. 
There are at most $6$ knots with up to 9 crossings which are 
neither 2-algebraic nor 3-braid knots. They are: 
$9_{34},9_{39}, 9_{40}, 9_{41}, 9_{47}$ and $9_{49}$. 
We reduced three of them, $9_{39},9_{41}$ and $9_{49}$ at my Fall 2003 
Dean's Seminar. The knot $9_{40}$ was reduced in December of 
2003 by Slavik Jablan and Radmila Sazdanovic.
Soon after, my student Maciej Niebrzydowski simplified the remaining 
pair $9_{34}$ and $9_{47}$, Fig.1.24.
\end{proof}

\ \\
\centerline{\psfig{figure=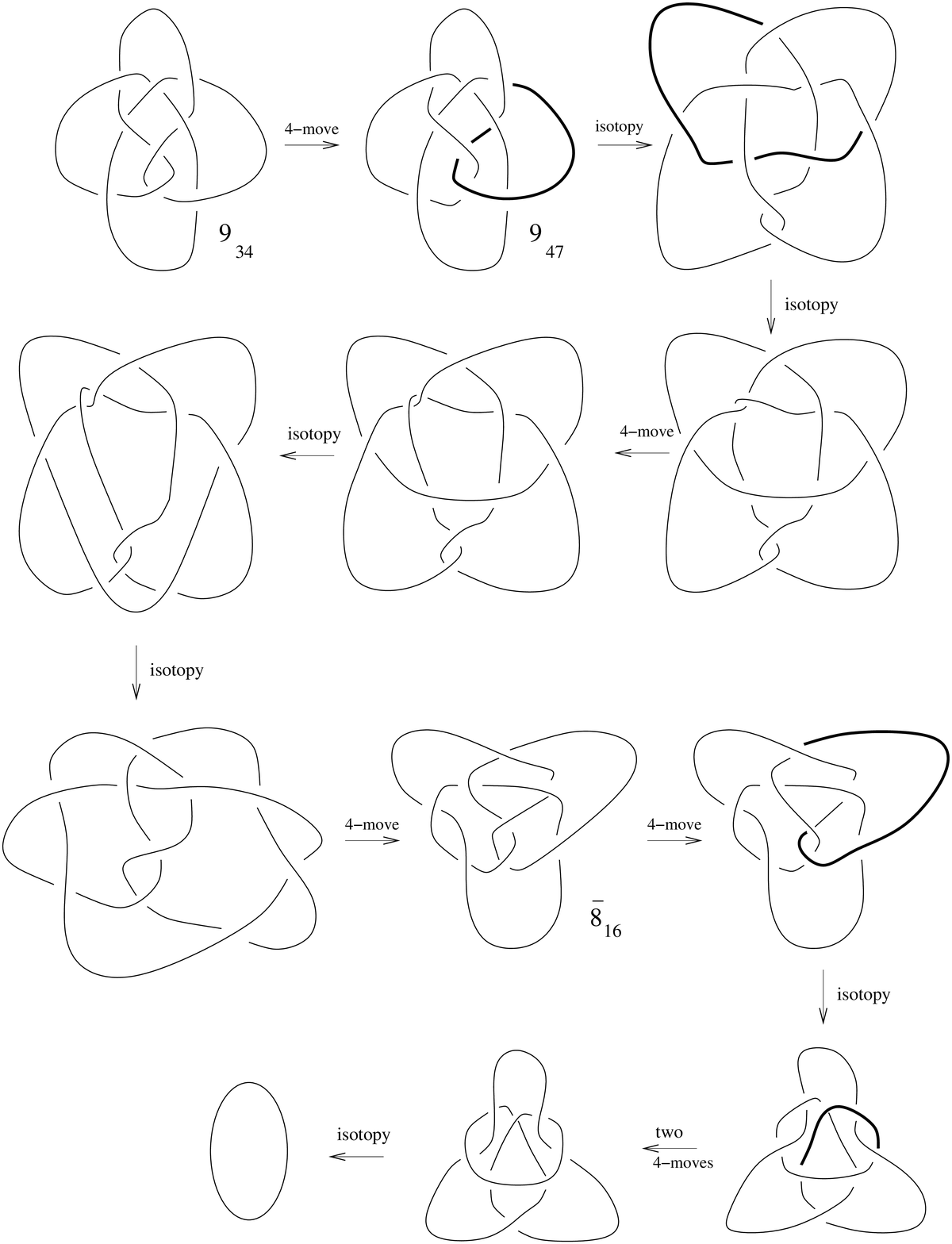,height=17.8cm}}
\ \\
\centerline{Fig. 1.24}\ \\

A weaker version of the Kawauchi question has been answered by Nakanishi 
in 1989, \cite{Nak-1}. 
If $\gamma \in B_n$ then the {\it $\gamma$-move} is the $n$-tangle move 
in which the trivial $n$-braid is replaced by the braid $\gamma$. 

\begin{theorem}[Nakanishi]\label{1.16}
If two links $L_1$ and $L_2$ have the same linking matrix modulo 2, 
then $L_2$ can be obtained from $L_1$ by a finite number of $\pm 4$-moves 
and $\Delta_3^4$-moves.
\end{theorem}
\begin{proof} The square of the center, $\Delta_3^4=(\sigma_1\sigma_2)^6$,
 of the 3-braid group $B_3$ and the Borromean braid, 
$(\sigma_1\sigma_2^{-1})^3$, are equal\footnote{We have in $B_3/(\sigma_i^4)$: 
$(\sigma_1\sigma_2)^6=
(\sigma_1^2\sigma_2\sigma_1^2\sigma_2)(\sigma_1\sigma_2^2\sigma_1\sigma_2^2)=
\sigma_1^2\sigma_2\sigma_1^2(\sigma_1\sigma_2^2\sigma_1\sigma_2^2)\sigma_2 
 \stackrel{4}{=} 
\sigma_1^2\sigma_2\sigma_1^{-1}\sigma_2^2\sigma_1\sigma_2^{-1} =
 \sigma_1\sigma_2^{-1}\sigma_1\sigma_2\sigma_2^2\sigma_1\sigma_2^{-1}
\stackrel{4}{=} (\sigma_1\sigma_2^{-1})^3$. This calculation can be 
interpreted as an illustration of Fig.28 in \cite{A-K}.}
 in $B_3/(\sigma_i^4)$. From this it also follows that $\Delta_3^4$ and 
$\Delta_3^{-4}$ are equal in $B_3/(\sigma_i^4)$. Furthermore, 
the $(\sigma_1\sigma_2^{-1})^3$-move is equivalent to $\Delta$-move 
of Nakanishi in which $\sigma_1\sigma_2^{-1}\sigma_1$ is replaced by 
$\sigma_2\sigma_1^{-1}\sigma_2$ (we can think of this move as a ``false" 
braid relation or a ``false" third Reidemeister move). Nakanishi proved 
that two oriented links are $\Delta$-move equivalent if and only if 
their linking matrices are equivalent \cite{Nak-1}. Theorem 1.16 follows.
\end{proof}

Selman Akbulut used Nakanishi's theorem to prove John Nash's 
conjecture for 3-dimensional manifolds 
\cite{A-K}\footnote{The conjecture that ``any two closed smooth connected 
manifolds of the same dimension can be made diffeomorphic after blowing them 
up along submanifolds" is an interpretation of the Nash question 
``Is there an algebraic structure on a any given smooth manifold 
which is birational to $RP^n$?"  
\cite{Nash,A-K}. The conjecture is only loosely related to the question 
mentioned in the book ``A beautiful mind" were in Chapter ``The 
`Blowing Up' Problem", it is written: ``Nash seemed, as the Fall [1963] 
unfolded, to be in far better shape
than he had been during his previous interlude at the Institute [IAS].
As he said in his Madrid lecture, he ``had had an idea which is referred
to as Nash Blowing UP which I discussed with an
eminent mathematician named Hironaka." [Letter from J.Nash to V.Nash,
1.9.66] (Hironaka eventually wrote a conjecture up.)" \cite{Nas}.}.\\
\ \\

It is not true that every link is $5$-move equivalent to a trivial link.
One can show, using the Jones polynomial, that the figure eight knot is
not 5-move equivalent to any trivial link\footnote{A $5$-move preserves 
the absolute value of the Jones polynomial 
at $t=e^{\pi i/5}$ \cite{P-1}. 
However, the Jones polynomial 
$V_{4_1}(e^{\pi i/5})=0$ but for any trivial link, $T_n$, we have 
$V_{T_n}(e^{\pi i/5})=(-e^{\pi i/10}-e^{-\pi i/10})^{n-1}\neq 0$.}. 
One can, however, introduce
a more delicate move, called $(2,2)$-move 
(\parbox{1.8cm}{\psfig{figure=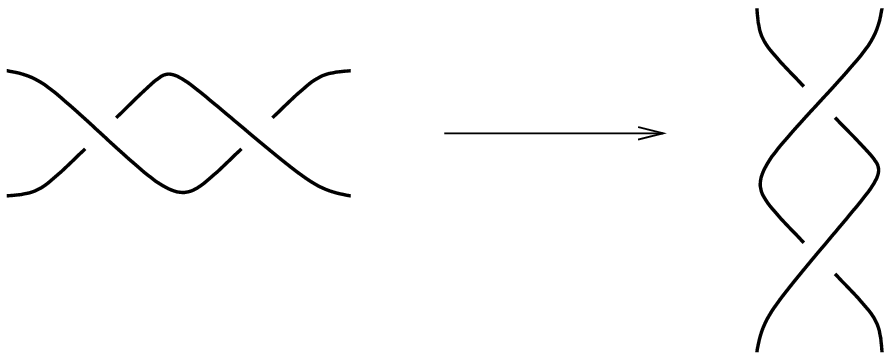,height=0.7cm}}),
 such that the $5$-move is a combination of a $(2,2)$-move 
and its mirror image $(-2,-2)$-move
(\parbox{2.1cm}{\psfig{figure=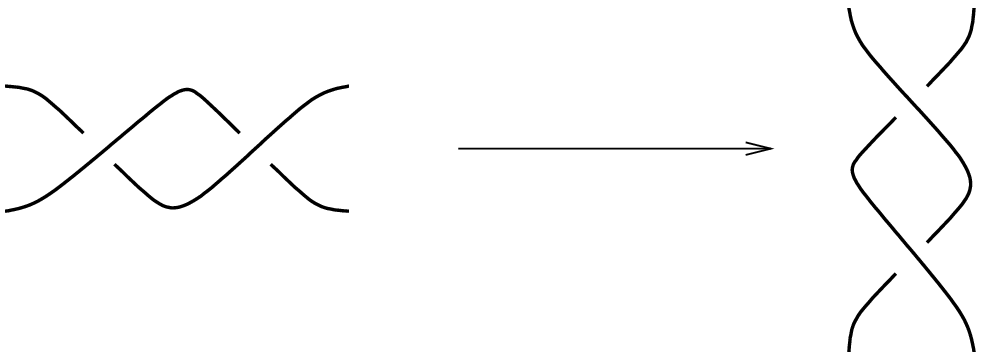,height=0.7cm}}), 
as it is illustrated in Fig.1.25 \cite{H-U,P-3}.\\
\ \\

\centerline{\psfig{figure=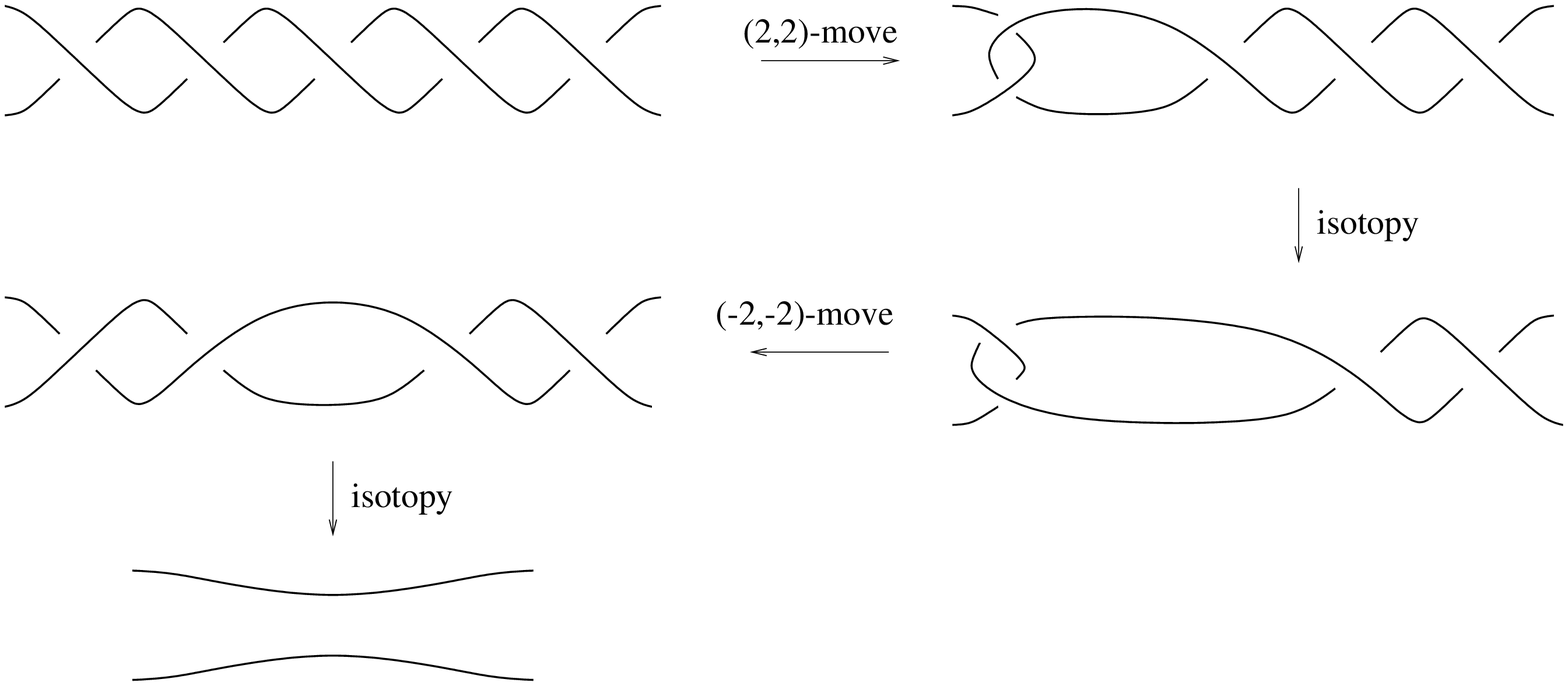,height=5.1cm}}
\centerline{Fig. 1.25}

\begin{conjecture}[Harikae, Nakanishi, Uchida, 1992]\label{1.17}\ \\
Every link is $(2,2)$-move equivalent to a trivial link.
\end{conjecture}

As in the case of 3-moves, an elementary induction shows that
the conjecture holds for 2-algebraic links. It is also known 
that the conjecture holds for all links up to $8$ crossings. 
The key element of the argument in the proof 
is the observation (going back to Conway \cite{Co}) that any link 
with up to 8 crossings (different from $8_{18}$; see  footnote 19) 
is 2-algebraic.
The reduction of the $8_{18}$ knot to
the trivial link of two components by my students, 
Jarek Buczy\'nski and Mike Veve, is 
illustrated in Fig.1.26.
\ \\
\centerline{\psfig{figure=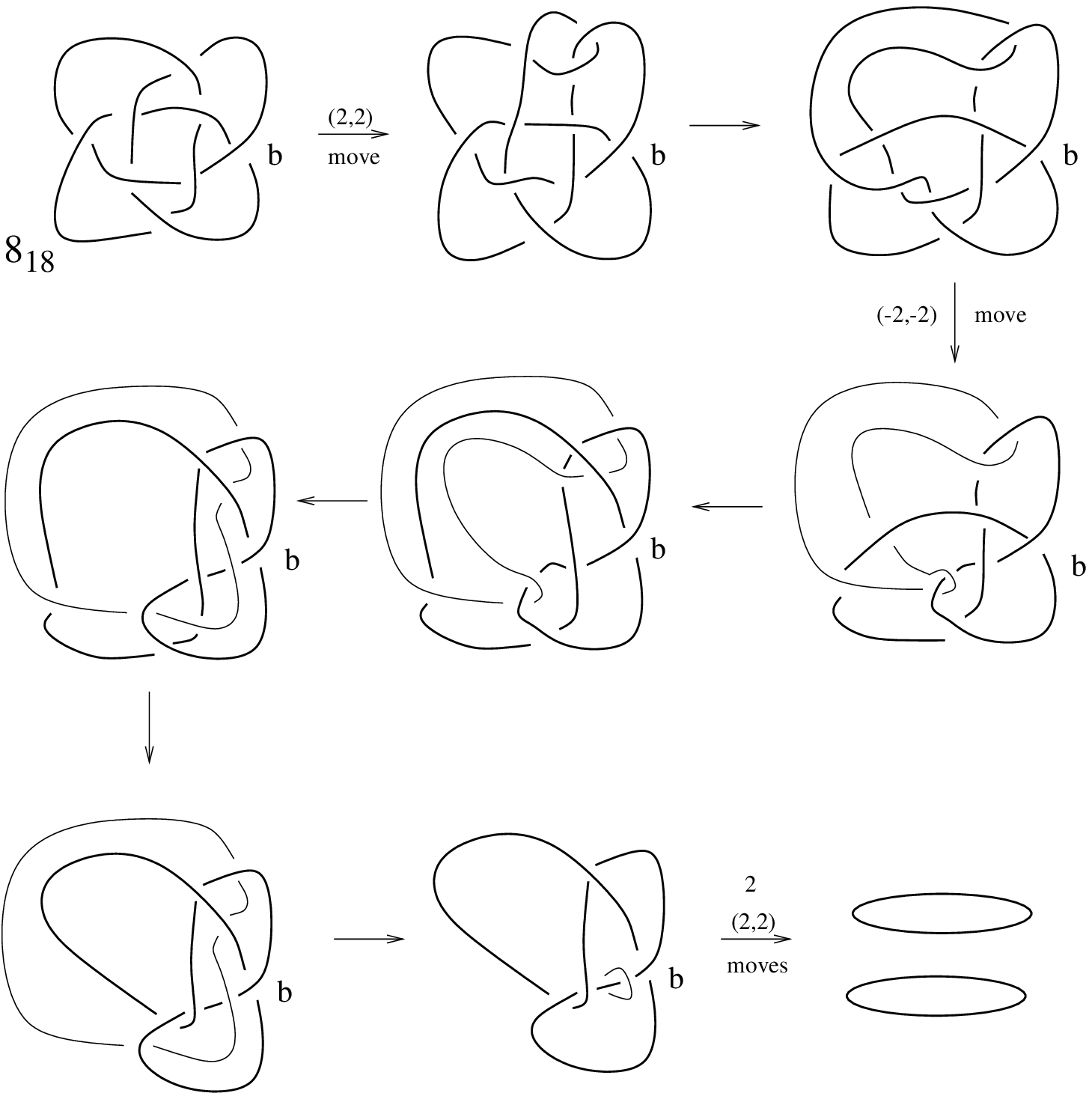,height=13.2cm}}

\centerline{Fig. 1.26;\ Reduction of the $8_{18}$ knot}
\ \\

The smallest knots that are not reduced yet are $9_{40}$ and $9_{49}$, 
Fig.1.27.\\
Possibly you can reduce them!\footnote{We showed with M.~D{\c a}bkowski 
that the knots $9_{40}$ and $9_{49}$ are not $(2,2)$-move 
equivalent to trivial links \cite{D-P-2}. Possibly you can 
prove that they are in the same $(2,2)$-move equivalence class! 
If I had to guess, I would say that it is a likely possibility.}
\ \\

\centerline{\psfig{figure=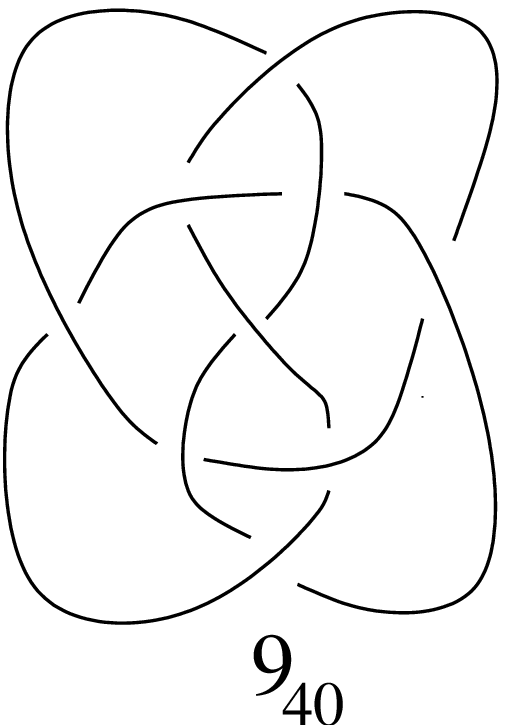,height=4.1cm}\ \ \ \ \ \  \ \ \  
\psfig{figure=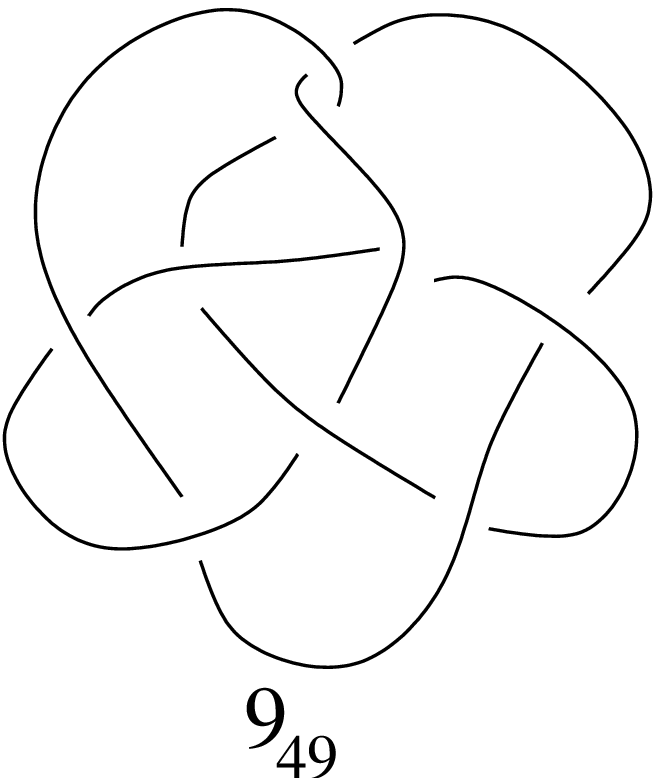,height=4.1cm}}
\centerline{Fig. 1.27}
\ \\

I am much less convinced that the answer to the next open question 
is positive, so I will not call it a ``conjecture".
First let us define a $(p,q)$-move to be a local modification 
of a link as shown in Fig.1.28.
We say that two links, $L_1$ and $L_2$, are $(p,q)$-equivalent 
if one can obtain one from the other 
by a finite number of $(p,q)$-,$(q,p)$-,$(-p,-q)$- 
and $(-q,-p)$-moves.\\
\ \\
\centerline{\psfig{figure=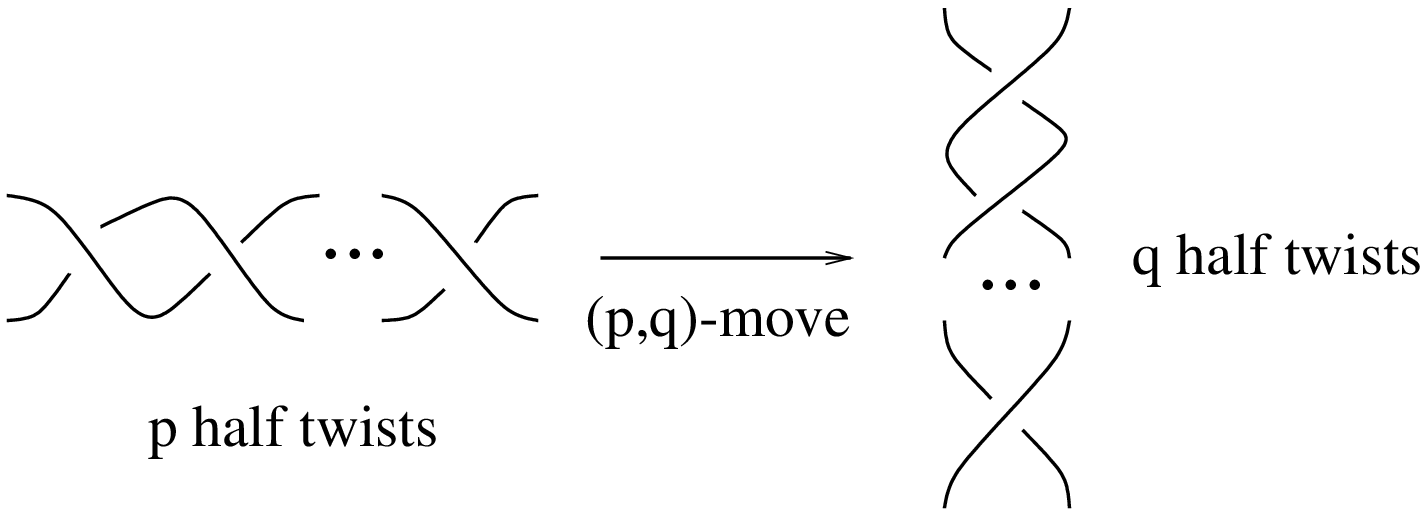,height=3.6cm}}
\centerline{Fig. 1.28}

\begin{problem}[\cite{Kir}; Problem 1.59(7), 1995]\label{1.18}
Is it true that any link is $(2,3)$-move equivalent to a trivial link?
\end{problem}

\begin{example}\label{1.19} 
Reduction of the trefoil and the figure eight knots is illustrated
in Fig.1.29. Reduction of the Borromean rings is shown in Fig.1.30.
\end{example}
\centerline{\psfig{figure=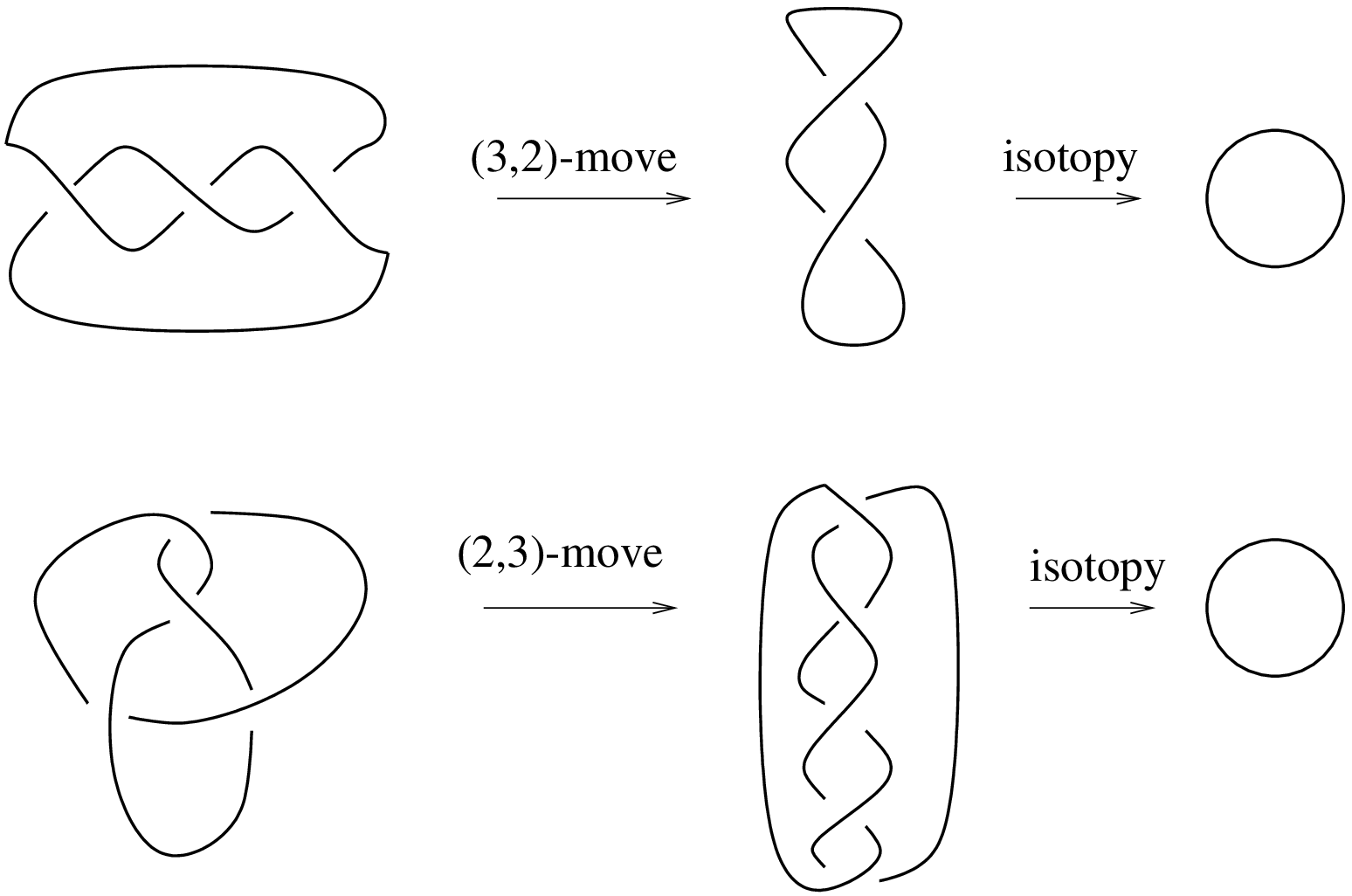,height=6.1cm}}
\centerline{Fig. 1.29}

\centerline{\psfig{figure=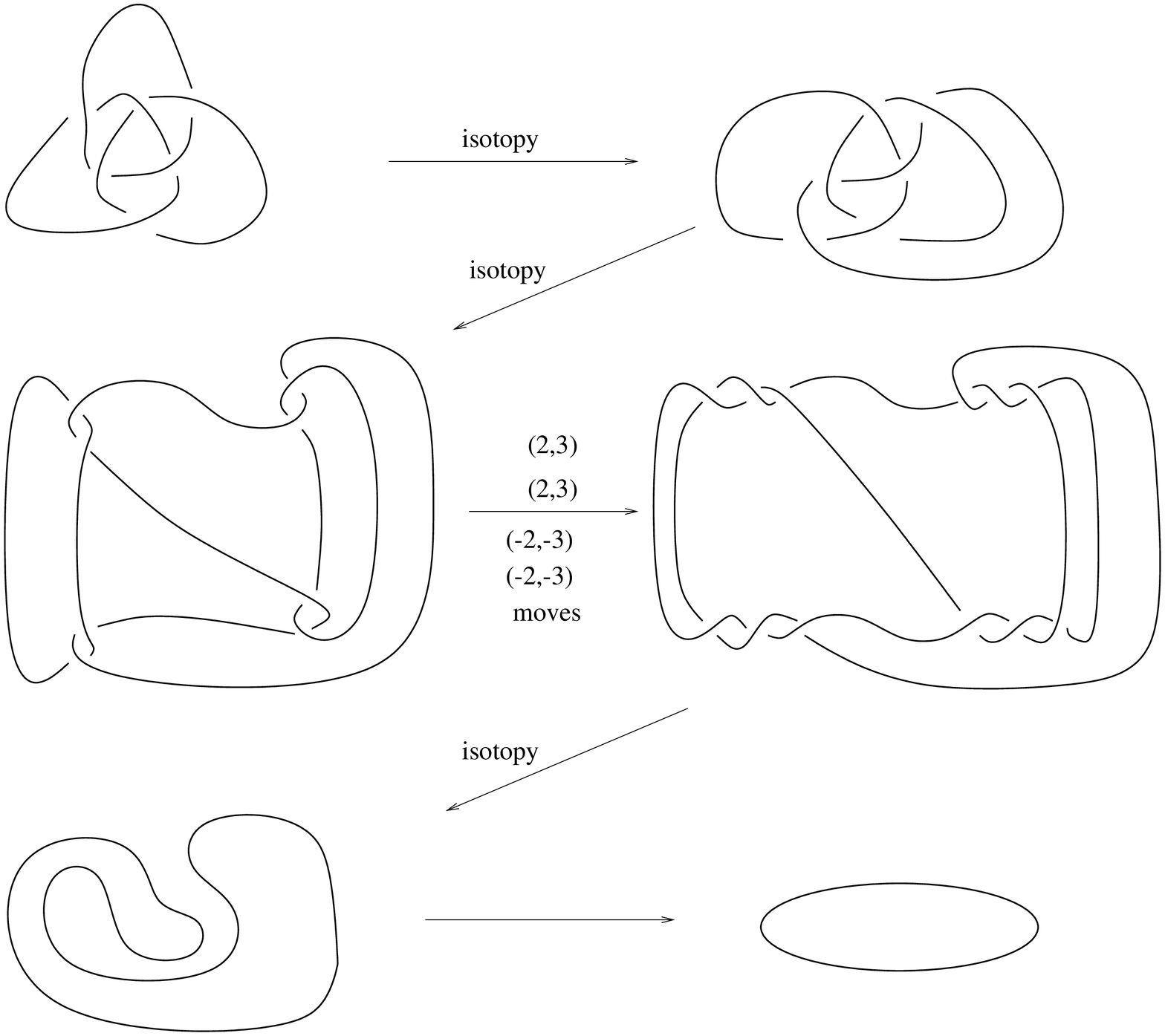,height=9.1cm}}
\centerline{Fig. 1.30}
\ \\

As in the case of Proposition 1.9, simple inductive argument 
shows that $2$-algebraic 
links are $(2,3)$-move equivalent to trivial links. 
Fig.1.31 
illustrates why the Borromean rings are 2-algebraic. By a proper
filling of black dots one can also show that all links 
with up to 8 crossings, except $8_{18}$, are 2-algebraic. 
Thus, as in the case of $(2,2)$-equivalence, 
the only link with up to $8$ crossings 
which still should be checked is the $8_{18}$ 
knot\footnote{To prove
that the knot $8_{18}$ is not 2-algebraic
one considers the 2-fold branched cover of $S^3$ branched along 
the knot, $M^{(2)}_{8_{18}}$. Montesinos proved that
algebraic knots are covered by Waldhausen graph manifolds \cite{Mo-1}.
Bonahon and Siebenmann showed (\cite{B-S}, Chapter 5) 
that $M^{(2)}_{8_{18}}$ is a hyperbolic 3-manifold so it cannot 
be a graph manifold.
This manifold is interesting from the point of view of hyperbolic
geometry because it is a closed manifold with its volume equal
to the volume of the complement of figure eight knot \cite{M-V-1}.
The knot $9_{49}$ of Fig.1.27 is not 2-algebraic either because its
2-fold branched cover is a hyperbolic 3-manifold. In
fact, it is the manifold
I suspected from 1983 to have the smallest volume among
oriented hyperbolic 3-manifolds \cite{I-MPT,Kir,M-V-2}.
In February of 2002 we (my student M.D{\c a}bkowski and myself) found 
unexpected connection between Knot Theory and the theory of Burnside groups.
This has allowed us to present simple combinatorial 
proof that the knots $9_{40}$ and $9_{49}$
are not 2-algebraic. However, our method does not work for the knot 
$8_{18}$ \cite{D-P-1,D-P-2,D-P-3}.}.
 Nobody really worked on this problem seriously, 
so maybe somebody in the audience will try this puzzle.\\
\ \\
\centerline{\psfig{figure=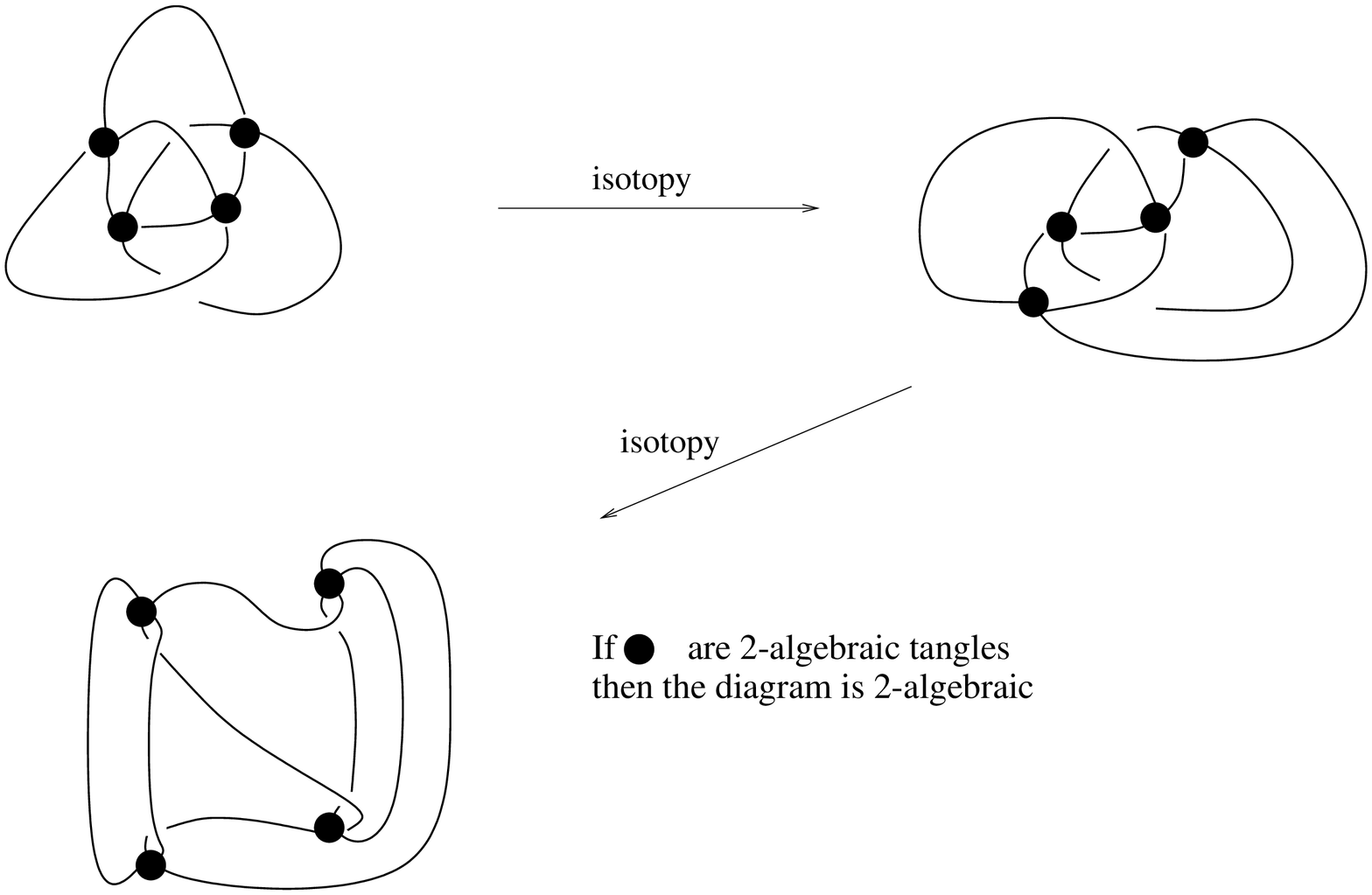,height=9.1cm}}
\centerline{Fig. 1.31}
\ \\
\ \\
{\bf Fox colorings}\\
The 3-coloring invariant which we are going to use to show that different 
trivial links are not 3-move equivalent, 
was introduced by R.~H. Fox\footnote{Ralph Hartzler Fox was born 
March 24, 1913.
A native of Morrisville, Pa., he attended Swarthmore
College for two years while studying piano
at the Leefson Conservatory of Music in Philadelphia.
He was mostly home schooled and later he was a witness in a
court case in Virginia, certifying soundness of home schooling.
He received his master's degree from the Johns Hopkins University
and his Ph.D. from the Princeton University in 1939 under
the supervision of Solomon Lefschetz.
Fox was married, when he was still a student, to Cynthia Atkinson.
They had one son, Robin.
After receiving his Princeton
doctorate, he spent the following year at
Institute for Advanced Study in Princeton.
He taught at the
University of Illinois and Syracuse University before
returning to join the Princeton University faculty in 1945
and staying there until his death. 
He was giving a series of lectures at the Instituto de Matem\'aticas
de la Universidad Nacional Aut\'onoma de M\'exico in the
summer of 1951. He was lecturing to American Mathematical Society (1949),
to the Summer Seminar of the Canadian Mathematical Society (1953),
and at the Universities of Delft and Stockholm, while on a Fulbright
grant (1952).
He died December 23, 1973 in the University
of Pennsylvania Graduate Hospital, where 
he had undergone open-heart surgery \cite{P-12}.}
 around 1956 when he was  
explaining Knot Theory to undergraduate students
at Haverford College (``in an attempt to make the subject
accessible to everyone" \cite{C-F}). It is a pleasant method 
of coding representations of the fundamental
group of a link complement into the group of symmetries 
of an equilateral triangle, however this interpretation is not needed 
for the definition and most of applications of 
3-colorings (compare \cite{Cr,C-F,Fo-1,Fo-2}).

\begin{definition}\label{1.20} (Fox 3-coloring of a link diagram).\\
Consider a coloring of a link diagram using colors r (red), 
y (yellow), and b (blue) 
in such a way that an arc of the diagram (from 
an undercrossing to an undercrossing) is colored by one color and at each 
 crossing one uses either only one or all three colors. Such a coloring is
called a \textup{Fox 3-coloring}. If the whole diagram is colored by  just one 
color we say that we have a \textup{trivial coloring}. 
The number of different Fox 3-colorings of $D$ is denoted by $tri(D)$. 
\end{definition}
\begin{example}
\begin{enumerate}
\item[\textup{(i)}] $tri(
\parbox{.6cm}{\psfig{figure=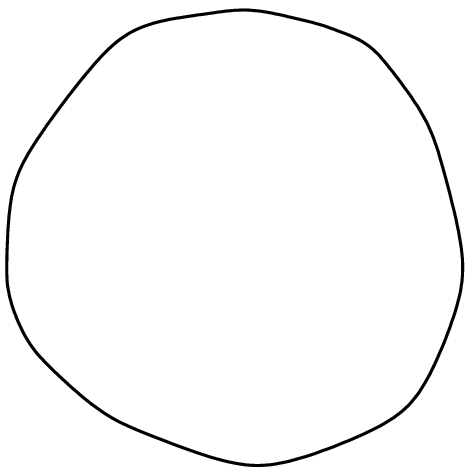,height=0.5cm}}
) = 3$ as the trivial link diagram has only trivial colorings.
\item[\textup{(ii)}] $tri(
\parbox{1.2cm}{\psfig{figure=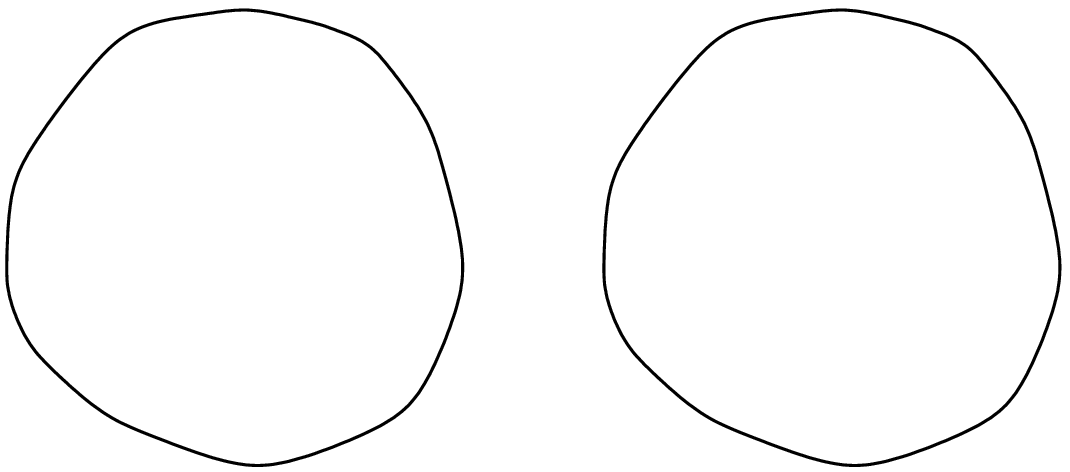,height=0.5cm}}
)=9$, and more generally, for the trivial link diagram of
$n$ components, $T_n$, one has $tri(T_n) = 3^n$.
\item[\textup{(iii)}] 
For the standard diagram of the right-handed trefoil knot we have 
three trivial colorings
and six nontrivial colorings. One of them is presented in Fig.1.32   
(all the others differ
from this one by permutations of colors). Thus, $tri(
\parbox{.9cm}{\psfig{figure=aatrefoil.eps,height=.7cm}})
 = 3+6=9$.
\end{enumerate}
\end{example}
\centerline{\psfig{figure=trefoil3-col.eps,height=3.5cm}}
\centerline{Fig. 1.32; \ Different colors are marked by lines of 
different thickness.}
\ \\

Fox 3-colorings were defined for link diagrams. They are, however,  
invariants of links. One only needs to show that $tri(D)$ 
is unchanged by Reidemeister moves.

The invariance under $R_1$ and $R_2$ is illustrated in Fig.1.33 and the
invariance under $R_3$ is illustrated in Fig.1.34.
\ \\
\ \\
\ \\
\centerline{\psfig{figure=R1R23-col.eps,height=1.9cm}}
\begin{center}
Fig. 1.33
\end{center}

\centerline{\psfig{figure=R33-col.eps,height=6.3cm}}
\begin{center}
Fig. 1.34
\end{center}

The next property of Fox 3-colorings is the key in proving that
different trivial links are not 3-move equivalent.

\begin{lemma}[\cite{P-1}]\label{1.22}  
3-moves do not change $tri(D)$.
\end{lemma}
The proof of the lemma is illustrated in Figure 1.35.

\ \\
\centerline{\psfig{figure=3-movecol.eps}}
\begin{center}
Fig. 1.35
\end{center}
The lemma also explains the fact that the trefoil knot has nontrivial
Fox 3-colorings: the trefoil knot is 3-move equivalent to
the trivial link of two components (Example 1.4(i)). 

Tomorrow, I will place the theory of Fox colorings in a more
general (sophisticated) context, and apply it to the analysis 
of 3-moves (and $(2,2)$- and $(2,3)$-moves) on $n$-tangles.
Interpretation of tangle colorings as Lagrangians in symplectic
spaces is our main (and new) tool. In the third section,
I will discuss another motivation for studying 3-moves: 
 understanding skein modules based on their deformation.

\section{Lagrangian approximation of Fox $p$-colorings of tangles}\label{2}

We just had the opportunity to listen to a beautiful and elementary 
talk by Lou Kauffman. I hope to follow this example by making my talk
elementary and deep at the same time. I will use several
results introduced by Lou, like classification of rational tangles, 
and also I am going to build on my yesterday's talk. 
I will culminate today's talk with introduction of the symplectic 
structure on the boundary of a tangle in such a way that tangles 
will yield Lagrangians in the associated symplectic space. 
I could not dream of this connection a year ago; however, now, 
 10 months after, I see the symplectic structure 
as a natural development.

Let us start our discussion slowly using my personal perspective 
and motivation.
In the Spring of 1986, I was analyzing behavior of Jones type invariants
of links when modified by  $k$-moves (or $t_k$-, $\bar t_{2k}$-moves 
in the oriented case). My interest had its roots in the
fundamental paper by Conway \cite{Co}. In July of 1986, I gave 
a talk at the ``Braids" conference in Santa Cruz. After 
my talk, I was told by Kunio Murasugi and Hitoshi Murakami 
about the Nakanishi's 3-move conjecture.
It was suggested to me by R. Campbell (Rob Kirby's student in 1986) 
to consider the effect of 3-moves on Fox colorings. 
Several years later, when writing \cite{P-3} in 1993, 
I realized that Fox colorings can be successfully used to analyze 
moves on tangles by considering not only the space of colorings 
but also the induced colorings of boundary points. 
More of this later, but let us now define Fox $k$-colorings first.
\begin{definition}\label{2.1}
\begin{enumerate}
\item [\textup{(i)}]
We say that a link (or a tangle) diagram is k-colored if every
arc is colored by one of the numbers $0,1,...,k-1$ (elements of the 
group ${\mathbb Z}_k$) in such a way that
at each crossing the sum of the colors of the undercrossings equals
twice the color of the overcrossing modulo $k$; see Fig.2.1. 
\item [\textup{(ii)}]
The set of $k$-colorings forms an abelian group, denoted by $Col_k(D)$ 
(we can also think of $Col_k(D)$ as a module over $\mathbb{Z}_k$).
The cardinality of the group will be denoted by $col_k(D)$.
For an $n$-tangle $T$ each Fox $k$-coloring of $T$ yields a
coloring of boundary points of $T$ and we have the homomorphism
$\psi :Col_k(T) \rightarrow \mathbb{Z}_k^{2n}$
\end{enumerate}
\end{definition}.
\centerline{\psfig{figure=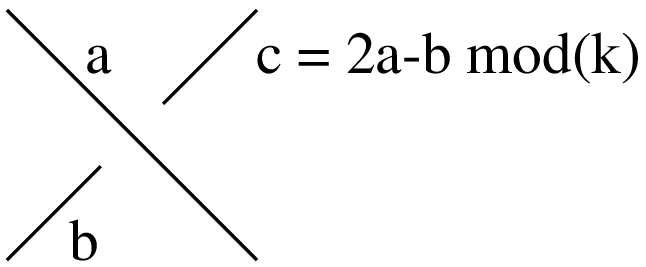,height=2.6cm}}
\centerline{Fig. 2.1}
\ \\

It is a pleasant exercise to show that $Col_k(D)$ is unchanged
by Reidemeister moves, so I am going to leave it for you.
The invariance  under $k$-moves is explained in Fig.2.2. 
\\
\ \\
\centerline{\psfig{figure=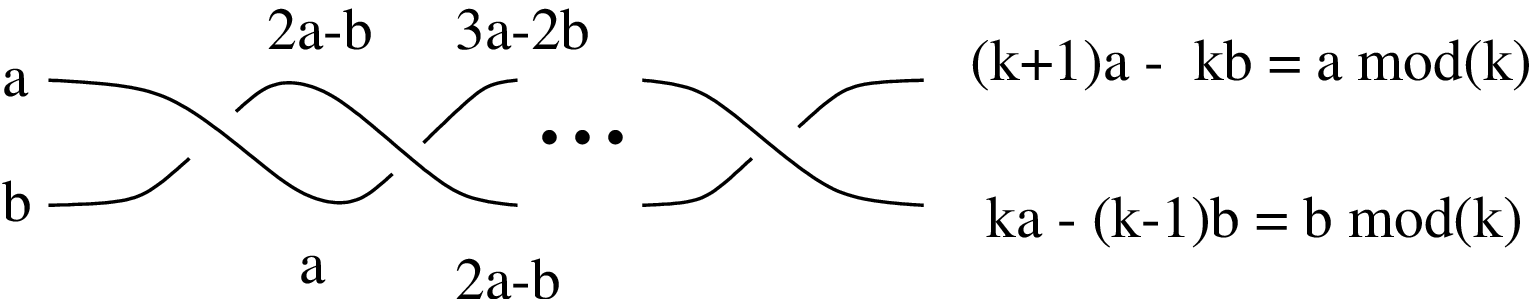,height=2.6cm}}
\centerline{Fig. 2.2}
\ \\

Having observed that $k$-moves
preserve the space of Fox $k$-colorings, 
let us take a closer look at the unlinking conjectures described 
before.
We discussed the 3-move conjecture, 
the 4-move conjecture for knots, and the Kawauchi's  
question for links. 
As I mentioned yesterday, not every link can be simplified using
$5$-moves, but the $5$-move is a combination of $\pm(2,2)$-moves 
and these moves might be sufficient to reduce every link to trivial links. 
Similarly not every link can be reduced via
$7$-moves, but again each $7$-move is a combination of 
$(2,3)$-moves\footnote{To be precise, a $7$-move is a combination
of $(-3,-2)$- and $(2,3)$-moves; compare Fig.1.25.}
which still might be sufficient for reduction. We stopped at this
point yesterday, but what could be used instead of general $k$-moves?
Let us consider the case of $p$-moves, where $p$ is a prime number.
I suggest (and state publicly for the first time)
that possibly one should consider {\it rational moves} instead, 
that is, moves in which
a rational $\frac{p}{q}$-tangle of Conway is substituted in
place of the identity tangle\footnote{The move was first
considered by J.~M. Montesinos \cite{Mo-2}; 
compare also Y.~Uchida \cite{Uch}.}. 
The most important observation for us is that
$Col_{p}(D)$ is preserved by $\frac{p}{q}$-moves. Fig.2.3
illustrates, for example, the fact that $Col_{13}(D)$ is unchanged by
a $\frac{13}{5}$-move.
\\
\ \\
\centerline{\psfig{figure=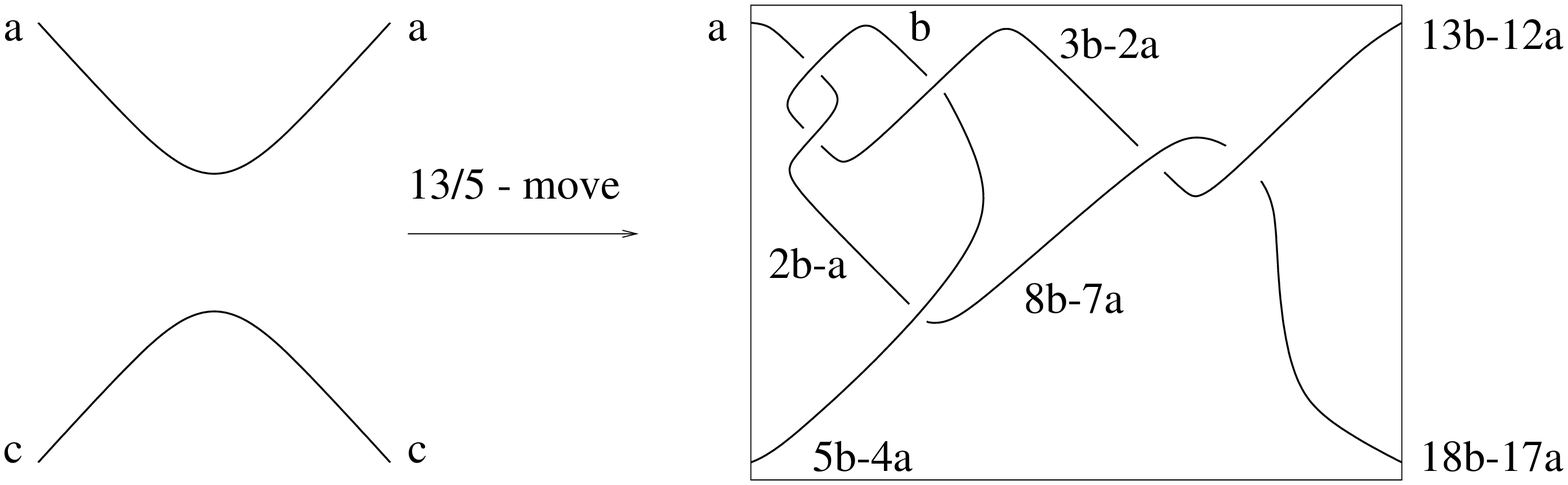,height=3.6cm}}
\centerline{Fig. 2.3}
\ \\

We also should note that $(m,q)$-moves are equivalent
to $\frac{mq+1}{q}$-moves (Fig.2.4) 
so the space of Fox $(mq+1)$-colorings
is preserved by them too.
\\
\ \\
\centerline{\psfig{figure=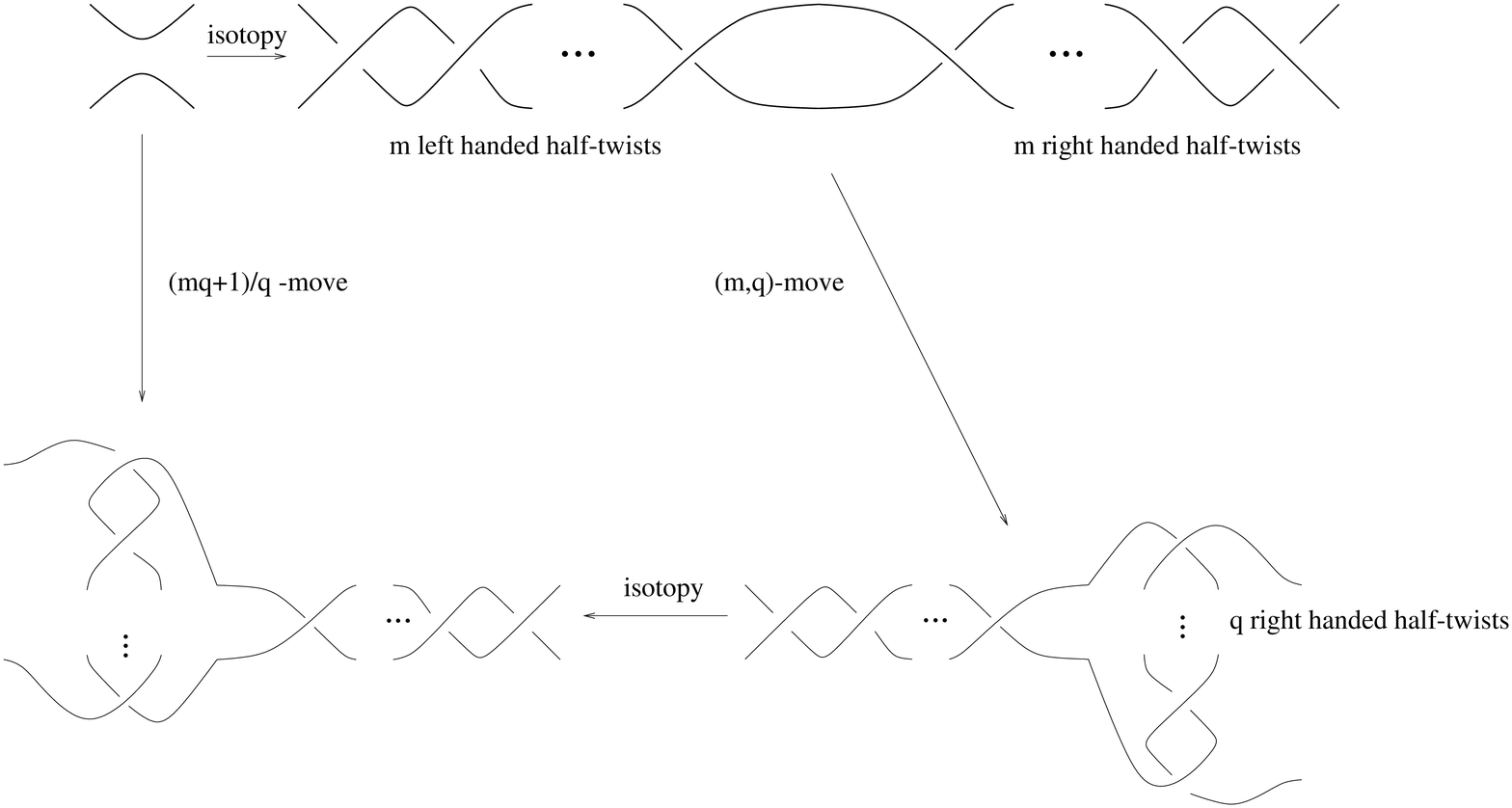,height=4.9cm}}
\centerline{Fig. 2.4}
\ \\

We have just heard about Conway's  classification
of rational tangles in Lou's talk\footnote{L.Kauffman's talk 
in Cuautitlan , March, 2001; compare \cite{K-L}.}, so  I will
just briefly sketch necessary definitions and introduce basic notation.
The 2-tangles shown in Fig.2.5 are called rational tangles
 -- in Conway's notation, $T(a_1,a_2,...,a_n)$. 
A rational tangle is 
$\frac{p}{q}$-tangle if $\frac{p}{q} = 
a_n + \frac{1}{a_{n-1}+...+\frac{1}{a_1}}$.\footnote{$\frac{p}{q}$
is called the slope of the tangle and can be easily
identified with the slope of the meridian disk of the solid torus
being the branched double cover of the rational tangle.} Conway
proved that two rational tangles are ambient isotopic
(with boundary points fixed) if and only if their slopes
are equal (compare \cite{Kaw}).
\\
\ \\
\centerline{\psfig{figure=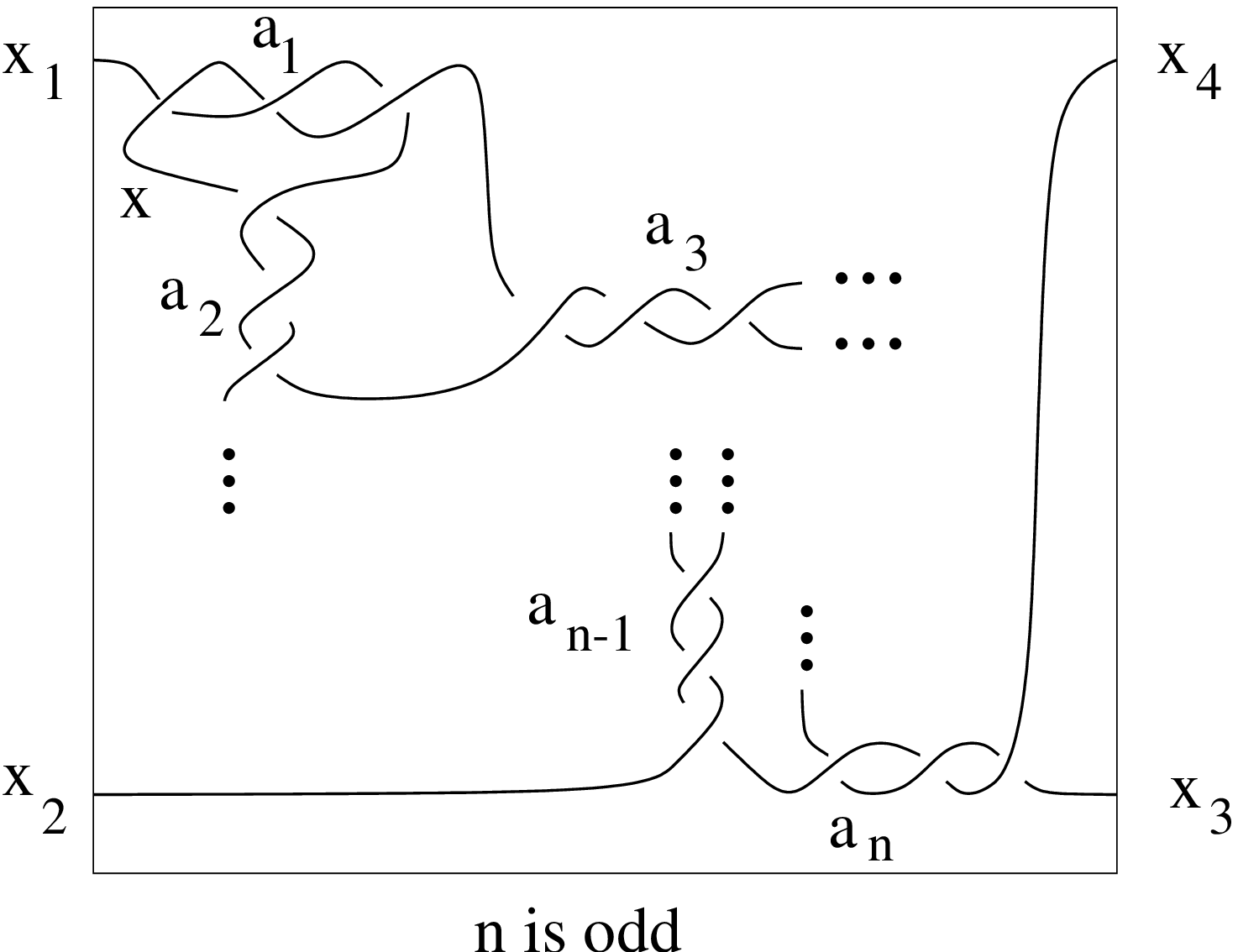,height=4.9cm}\ \ \ 
\psfig{figure=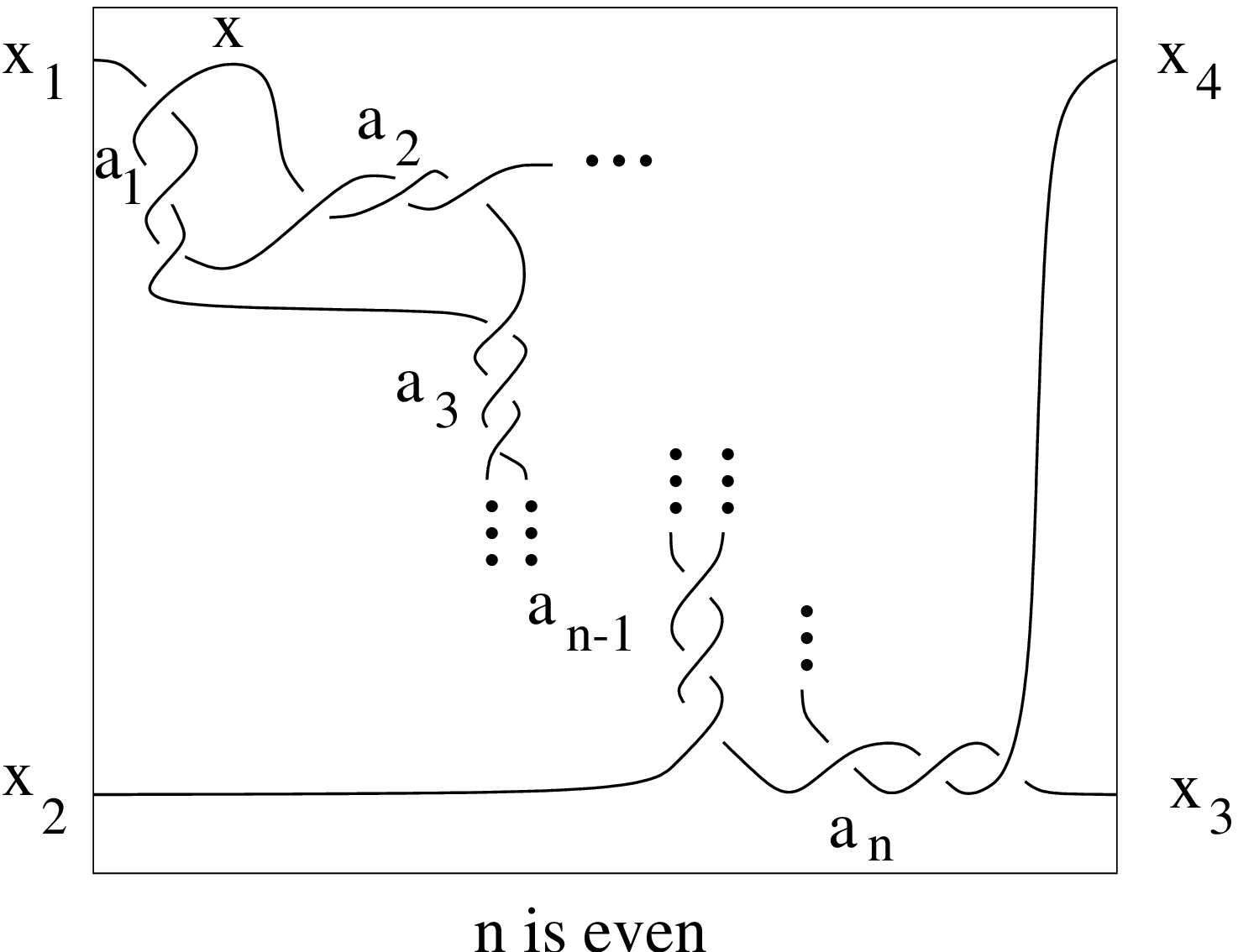,height=4.9cm}}
\centerline{Fig. 2.5}
\ \\

For a given Fox coloring of the rational $\frac{p}{q}$-tangle with
boundary colors $x_1,x_2,x_3,x_4$ (Fig.2.5), one has 
$x_4-x_1 = p(x-x_1)$, $x_2-x_1=q(x-x_1)$ and $x_3 = 
x_2 + x_4 -x_1$. If a coloring is nontrivial ($x_1\neq x$) 
then $\frac{x_4-x_1}{x_2-x_1} = \frac{p}{q}$ as it 
has been explained by Lou.

\begin{conjecture}\label{2.2}\ \\
Let $p$ be a fixed prime number, then\footnote{I decided to keep the 
word ``Conjecture" as it was used in my talk. However, in Spring 
of 2002, we disproved it for any $p$, \cite{D-P-1,D-P-2,D-P-3}. 
The talks in Mexico and Canada were essential for clarifying ideas 
and finally in constructing counterexamples.}
\begin{enumerate} 
\item[\textup{(i)}] Every link can
be reduced to a trivial link
by rational $\frac{p}{q}$-moves ($q$ any integer).
\item[\textup{(ii)}] There is a function
$f(n,p)$ such that any $n$-tangle can be reduced to one of
``basic" $f(n,p)$ $n$-tangles (allowing additional trivial components)
by rational $\frac{p}{q}$-moves.
\end{enumerate}
\end{conjecture}
First we observe that it suffices to use $\frac{p}{q}$-moves with
$|q| \leq \frac{p}{2}$, as they generate all the other $\frac{p}{q}$-moves 
follow. Namely, we have $\frac{p}{p-q} = 1 + \frac{1}{-1 + \frac{p}{q}}$ 
and $\frac{p}{-(p+q)} = -1 + \frac{1}{1 + \frac{p}{q}}$. 
Thus $\frac{p}{q}$-moves generate $\frac{p}{-q\pm p}$-moves (e.g.,  
$\frac{p}{p-q}$ tangle is reduced by an inverse of a 
$\frac{p}{q}$-move to the $0$-tangle,  $1 + \frac{1}{-1 + 0} = 0$). 
Furthermore, we know
that for odd $p$ the $\frac{p}{1}$-move is a combination 
of $\frac{p}{2}$ and $\frac{p}{-2}$-moves (compare Fig.1.25). 
Thus, in fact, 3-move, $(2,2)$-move and
$(2,3)$-move conjectures are special cases of Conjecture 2.2(i).
For $p=11$ we have $\frac{11}{2} = 5+ \frac{1}{2}$,
$\frac{11}{3} = 4 - \frac{1}{3}$, $\frac{11}{4} = 3 - \frac{1}{4}$,
$\frac{11}{5} = 2+ \frac{1}{5}$. Thus:
\begin{conjecture}\label{2.3}\ \\
Every link can be reduced to a trivial link (with the same 
space of 11-colorings)
by $(2,5)$- and $(4,-3)$-moves, their inverses and their mirror 
images\footnote{As mentioned in the footnote 25, Conjecture 2.3 does 
not hold. The closure of the $3$-braid $\Delta^4_3$ provides the 
simplest counterexample \cite{D-P-2,D-P-3}. However, it holds for 
2-algebraic links; see Proposition 2.7.}.
\end{conjecture}

What about the number $f(n,p)$?
We know that because $\frac{p}{q}$-moves preserve $p$-colorings, 
therefore $f(n,p)$ is bounded
from below by the number of subspaces of $p$-colorings 
of the $2n$ boundary points induced by Fox $p$-colorings 
of $n$-tangles (that is by the number
of subspaces $\psi (Col_p(T))$ in $\mathbb{Z}_p^{2n}$).
I noted in \cite{P-3} that for 2-tangles this number is equal to $p+1$
(even in this special case my argument was complicated).
For $p=3$ and $n=4$
the number of subspaces followed from the work of my student
Tatsuya Tsukamoto and is equal to $40$ \cite{P-Ts}.
The combined effort of Mietek D{\c a}bkowski and Tsukamoto
gave the number $1120$ for subspaces $\psi (Col_3(T))$ and 4-tangles.
That was my knowledge at the early Spring of 2000.
On May 2nd and 3rd I attended talks on Tits buildings
(at the Banach Center in Warsaw) by Janek Dymara and Tadek Januszkiewicz. 
I realized that the topic may have some 
connection to my work.
I asked Januszkiewicz whether he sees relation 
and I gave him numbers $4$, $40$, $1120$ for $p=3$.
He immediately answered that most likely I was counting the number
of Lagrangians in $\mathbb{Z}_3^{2n-2}$ symplectic space, and that
the number of Lagrangians in $\mathbb{Z}_p^{2n-2}$ is known to be
equal to $\prod_{i=1}^{n-1} (p^i +1)$.
Soon I constructed the appropriate symplectic form (as did
 Dymara). I will spend most of this talk on this 
construction and end with discussion of classes
of tangles for which it has been proved that 
$f(n,p) = \prod_{i=1}^{n-1} (p^i +1)$.

 Consider $2n$ points on a circle (or a square) and a
field $\mathbb{ Z}_p$ of $p$-colorings of a point. The colorings of
$2n$ points form the linear space $\mathbb{Z}_{p}^{2n}$. 
Let $e_1,\ldots , e_{2n}$ be the standard basis of $\mathbb{Z}_{p}^{2n}$,
\begin{center}
\centerline{\psfig{figure=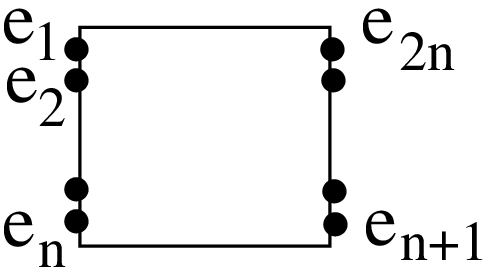,height=2.6cm}}
\end{center}
\begin{center} Fig. 2.6 \end{center}
\ \\
 $e_i=(0,\ldots , 1,\ldots ,0)$, where 1
occurs in the $i$-th position. Let $\mathbb{ Z}_{p}^{2n-1}
\subset\mathbb{ Z}_p^{2n}$ be the subspace of vectors $\sum a_i e_i$ 
satisfying $\sum (-1)^i a_i=0$ (alternating condition). 
Consider the basis $f_1,\ldots , f_{2n-1}$ of ${\bf
Z}_{p}^{2n-1}$ where $f_k=e_k+e_{k+1}$. 

Let 

$$\phi = \left( \begin{array}{cccc} 0 & 1 & \ldots &\ldots \\
 -1 & 0 & 1 &\ldots \\ \ldots & \ldots & \ldots & \ldots \\ \ldots & \ldots
 & -1 & 0
\end{array} \right)$$
\ \\

be a skew-symmetric
form $\phi$ on $\mathbb{ Z}_{p}^{2n-1}$ of nullity 1,
that is,
 \vspace{1mm}
   \renewcommand{\arraystretch}{3}
   $$\phi (f_i, f_j) =\left\{
   \begin{array}{lr}
   0 &
   {\rm if}\ |j-i|\neq 1 \\
  1 &  {\rm if\ }\ j=i+1\\ -1 &
   {\rm if}\ j=i-1. \\
   \end{array}
  \right .$$
   \par
   \vspace{2mm}

Notice that the vector $e_1+ e_2+ \ldots + e_{2n}$
($=f_1+f_3+\ldots + f_{2n-1}=f_2+f_4+\ldots + f_{2n})$
is $\phi$-orthogonal to any other vector.
If we consider $\mathbb{ Z}_{p}^{2n-2}= \mathbb{ Z}_{p}^{2n-1}/
\mathbb{ Z}_{p}$, where the subspace $\mathbb{ Z}_{p}$ is generated by
$e_1+\ldots + e_{2n}$, that is, $\mathbb{ Z}_{p}$ consists
of monochromatic (i.e., trivial) colorings,
then $\phi$ descends to the symplectic form 
$\hat\phi$ on $\mathbb{ Z}_{p}^{2n-2}$. 
Now we can analyze isotropic subspaces of
$(\mathbb{ Z}_{p}^{2n-2},\hat\phi)$, that is subspaces on
which  $\hat\phi$ is $0$
($W\subset \mathbb{Z}_{p}^{2n-2}, $ where $ 
\hat\phi (w_1,w_2)=0$ for all $w_1,w_2\in W$).
The maximal isotropic subspaces of 
$\mathbb{ Z}_{p}^{2n-2}$ are $(n-1)$-dimensional and they are
called Lagrangian subspaces (or maximal totally degenerated subspaces). 
There are $\prod_{i=1}^{n-1}(p^i+1)$ of them.

Our local condition on Fox colorings (Fig.2.1) guarantees 
that for any tangle $T$,
$\psi (Col_p T)\subset \mathbb{Z}_{p}^{2n-1}$. 
Furthermore, the space of trivial colorings, 
$\mathbb{ Z}_{p}$ is always in $Col_p T$.
Thus $\psi$ descends to 
$\hat\psi :Col_p T/\mathbb{ Z}_{p} \rightarrow
\mathbb{ Z}_{p}^{2n-2}= \mathbb{ Z}_{p}^{2n-1}/\mathbb{ Z}_{p}$. 
Now we answer
the fundamental question:
Which subspaces of
$\mathbb{ Z}_{p}^{2n-2}$ are yielded by $n$-tangles?

\begin{theorem}\label{2.4}\ \\
$\hat\psi (Col_p T/\mathbb{ Z}_{p})$ is a Lagrangian
subspace of $\mathbb{ Z}_{p}^{2n-2}$ with the symplectic form $\hat\phi$.
\end{theorem}
The natural question is whether every Lagrangian subspace
can be realized as a space of induced colorings on the boundary for some
 tangle. The answer is negative for
$p=2$ and positive for $p>2$.
\begin{theorem}[\cite{D-J-P}]\label{2.5}\ \\
\begin{enumerate}
\item[\textup{(i)}] For an odd prime number $p$, every Lagrangian 
in $(\mathbb{ Z}_{p}^{2n-2},
\hat\phi)$ is realized as $\hat\psi (Col_p T/\mathbb{ Z}_{p})$ for some 
$n$-tangle $T$. Furthermore, $T$ can be chosen to be a rational 
$n$-tangle\footnote{An $n$-tangle is a rational (or $n$-bridge) tangle 
if it is homeomorphic to a tangle without crossing and trivial 
components (we allow homeomorphism moving the boundary of the 3-ball). 
Alternatively, we can use an inductive definition modifying 
Definition 2.9 in such a way that we start from a tangle without a 
crossing and a trivial component and we assume in condition (i)(1) that 
$B$ has exactly one crossing (which is not nugatory, that is, it cannot be 
eliminated by a first Reidemeister move).}.
\item[\textup{(ii)}]
For $p=2$, $n>2$, not every Lagrangian is realized as 
$\hat\psi (Col_2 T/\mathbb{ Z}_{2})$. We have $f(n,2)= \prod_{i=1}^{n-1}(2i+1)$ 
(a 2-coloring is unchanged by a crossing change) but the number of Lagrangians 
is equal to $\prod_{i=1}^{n-1}(2^i+1)$.
\end{enumerate}
\end{theorem}

As a corollary we obtain a fact which was considered 
to be difficult before, even for 2-tangles (compare \cite{P-3,J-P}.
\begin{corollary}\label{2.6}\ \\
For any $p$-coloring $\bf x$ of a tangle boundary points satisfying 
the alternating property (i.e., ${\bf x} \in \mathbb{ Z}_{p}^{2n-1}$) 
there is an $n$-tangle and
 $p$-coloring of it that yields ${\bf x}$.
In other words: $\mathbb{ Z}_{p}^{2n-1} = \bigcup_T \psi_T(Col_pT)$.
Furthermore, the space $\psi_T(Col_pT)$ is $n$-dimensional for any $T$.
\end{corollary}
We can say that we understand the lower bound for  
the function $f(n,p)$, but when does Conjecture 2.2 holds with
 $f(n,p)= \prod_{i=1}^{n-1}(p^i+1)$?

In \cite{D-I-P} we discuss Conjecture 2.2 for 2-algebraic 
tangles. Here we sketch a proof of a simpler fact.
\begin{proposition}\label{2.7}\ \\
Let $p$ be a fixed prime number and let be $H_p$ the family of 2-tangles:
$\frac{1-p}{2}$, $\frac{3-p}{2}$,..., $0$,...,  $\frac{p-3}{2}$,
$\frac{p-1}{2}$ and $\infty$ (horizontal family), and let $V_p$ be the
 vertical family of 2-tangles, $V_p=r(H_p)$;
then
\begin{enumerate}
\item[\textup{(i)}]
Every $2$-algebraic tangle can be reduced to a 2-tangle from the family
$H_p$ (resp. $V_p$) with possible additional trivial components by
$(\frac{sp}{q})$-moves, where $s$ and $q$ are any integers such that
$sp$ and $q$ are relatively prime. 
Furthermore, for $p \leq 13$ one can assume that $s=1$. 
\item[\textup{(ii)}] 
Every $2$-algebraic link can be reduced to a trivial link by
$(\frac{sp}{q})$-moves, where $s$ and $q$ are any integers such that
$sp$ and $q$ are relatively prime.
 Furthermore, for $p \leq 13$ one can assume that $s=1$.
\end{enumerate}
\end{proposition}
{\bf Outline of the proof.}\
We use the structure of 2-algebraic tangles 
to perform an inductive proof similar to that of 
Proposition 1.9. The main problem in the proof is to show that the family 
$V_p$ can be reduced to the family $H_p$ by our moves. 
Consider the vertical tangle $r(k)$  where $k$ is relatively prime to $p$. 
There are integers $k'$ and $s$ such that $kk'+1 =sp$ or equivalently 
$k'+ \frac{1}{k} = \frac{sp}{k}$. Therefore the $\frac{sp}{k}$-move 
(equivalently $(k,k')$-move) is changing $r(k)$ to the horizontal 
tangle $k'$. In this reasoning 
we do not have a control over $s$. 
Consider now the case of $p=13$ and $s=1$. 
By considering fractions $\frac{13}{2} = 6 + \frac{1}{2}$,  
$\frac{13}{3}=4+\frac{1}{3}$,$\frac{13}{4}=3+\frac{1}{4}$, 
$\frac{13}{6}=2+\frac{1}{6}$, we are able to work with all $r(k)$ 
except $k=5$. $5+\frac{1}{5}= \frac{26}{5}$ so $s=2$ in this case. 
We can, however, realize $\frac{26}{5}$-move as a combination of
$\frac{13}{3}$-move and $\frac{13}{2}$-move as illustrated in
Fig.2.7 (we start by presenting $\frac{26}{5}$ as 
$6 + \frac{1}{-1 -\frac{1}{4}}$).\\
\ \\
\centerline{\psfig{figure=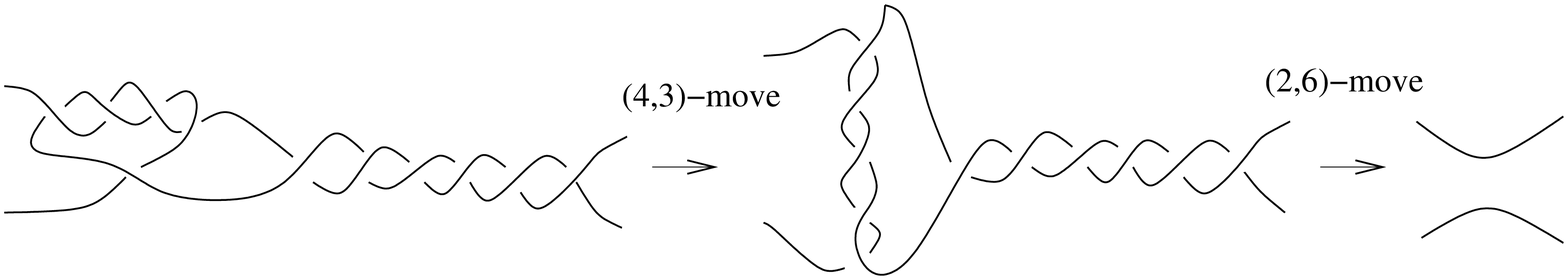,height=2.5cm}}
\centerline{Fig. 2.7}
\ \\

\begin{corollary}\label{2.8}\ \\
Every $(\frac{sp}{q})$-rational tangle, $p$ odd prime,
 can be reduced to the $0$ 2-tangle
by $(k,k')$-moves where $|k| <\frac{p}{2}$ and $kk'+1=sp$ for some $s$.
\end{corollary}

In order to be able to use induction for 
$n$-tangles with $n>2$, we generalize the notion of the algebraic tangle.
\begin{definition}\label{2.9}\ \\
\begin{enumerate}
\item[\textup{(i)}] The family of $n$-algebraic tangles 
is the smallest family 
of n-tangles which satisfies:\\
\textup{(0)} Any n-tangle with 0 or 1 crossing is n-algebraic.\\
\textup{(1)} If $A$ and $B$ are n-algebraic tangles then $r^i(A)*r^j(B)$ 
is n-algebraic, where $r$ denotes here the rotation of a tangle by 
$\frac{2\pi}{2n}$ angle, and $*$ denotes horizontal composition of tangles.

\item[\textup{(ii)}] If in the condition (1), $B$ is restricted 
to tangles with no more than $k$ crossings, we obtain 
the family of $(n,k)$-algebraic tangles.

\item[\textup{(iii)}] 
If an $m$-tangle, $T$, is obtained from an $(n,k)$-algebraic 
tangle (resp. $n$-algebraic tangle)
by partially closing its endpoints ($2n-2m$ of them) 
without introducing any new crossings, then $T$ is called
an $(n,k)$-algebraic (resp. $n$-algebraic) $m$-tangle. For $m=0$
we obtain an $(n,k)$-algebraic (resp. $n$-algebraic) link.
\end{enumerate}
\end{definition}
Conjecture 2.2, for $p=3$, has been proved for 3-algebraic 
tangles \cite{P-Ts}  ($f(3,3)=40$) and $(4,5)$-algebraic 
tangles \cite{Tsu} ($f(4,3)=1120$). In particular the 
Montesinos-Nakanishi 3-move conjecture holds for 3-algebraic and
$(4,5)$-algebraic links. 40 ``basic" 3-tangles are
shown in Fig. 2.8.\\
\ \\
\centerline{\psfig{figure=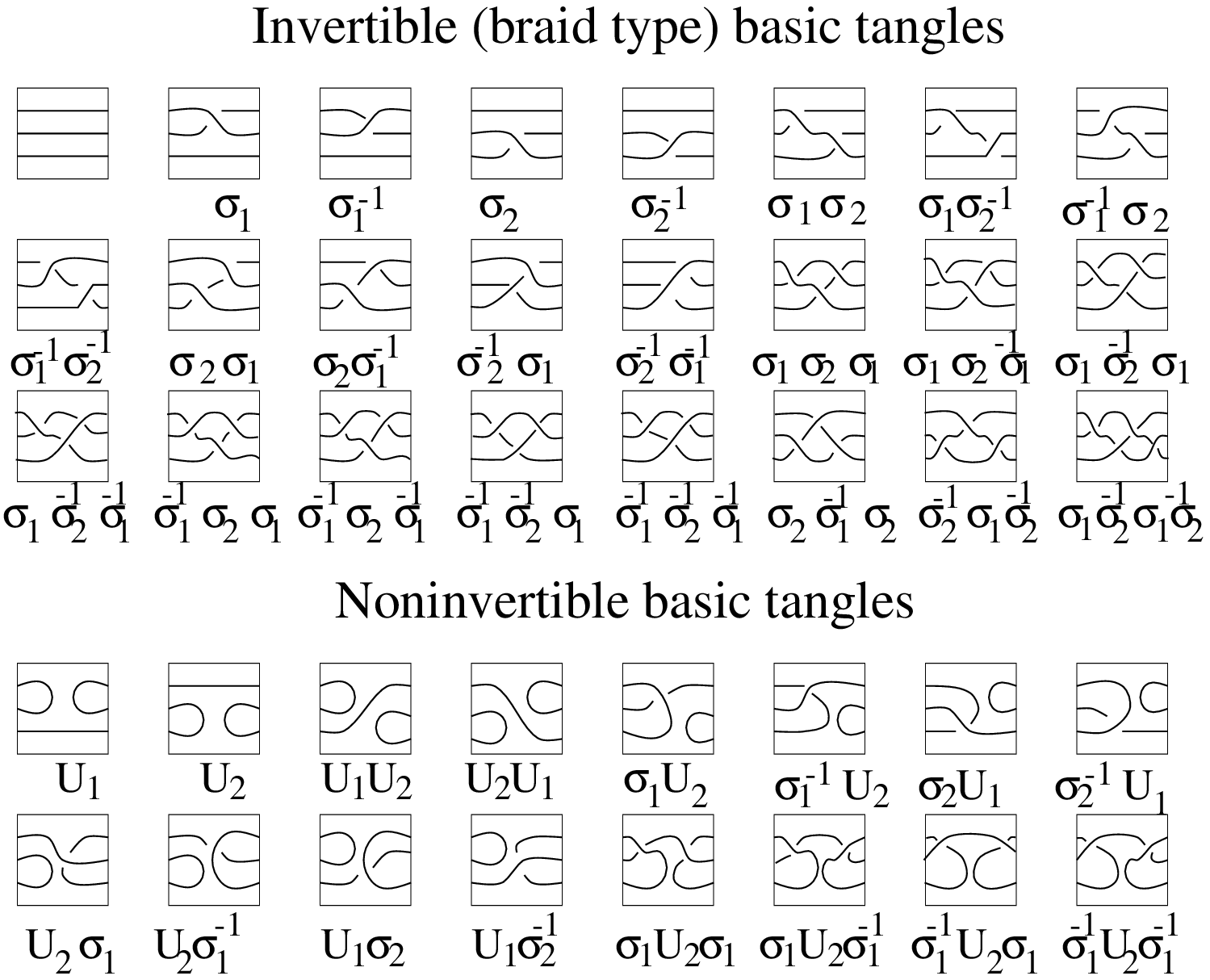,height=9.8cm}}
\centerline{Fig. 2.8}

The simplest $4$-tangles which cannot be distinguished by
3-colorings for which 3-move equivalence is not yet established
are illustrated in Fig.2.9. As for $(2,2)$-moves,
the equivalence of 2-tangles in Fig.2.10 is an open 
problem\footnote{The 4-tangles in Fig.2.9 are not 3-move equivalent. 
This follows from the fact that the Borromean rings and the Chen's 
link are not 3-move equivalent to trivial links \cite{D-P-1,D-P-3}. 
Similarly, the fact that 2-tangles of Fig.2.10 are not $(2,2)$-move 
equivalent follows from the result proven in \cite{D-P-2,D-P-3} 
that the knot $9_{49}$ and the link $9^2_{40}$ are not $(2,2)$-move 
equivalent to the trivial link of three components.}.\\
\ \\

\centerline{\psfig{figure=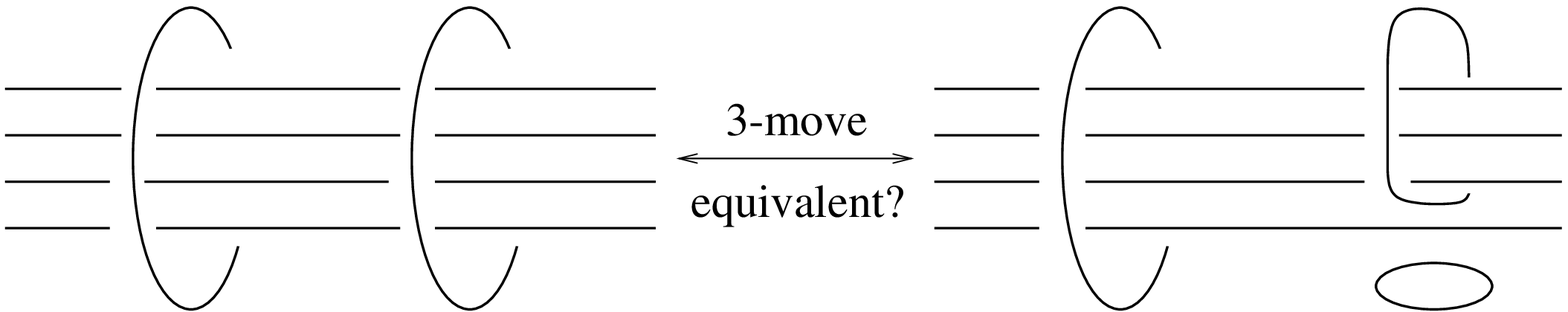,height=2.2cm}}
\centerline{Fig. 2.9}
\centerline{\psfig{figure=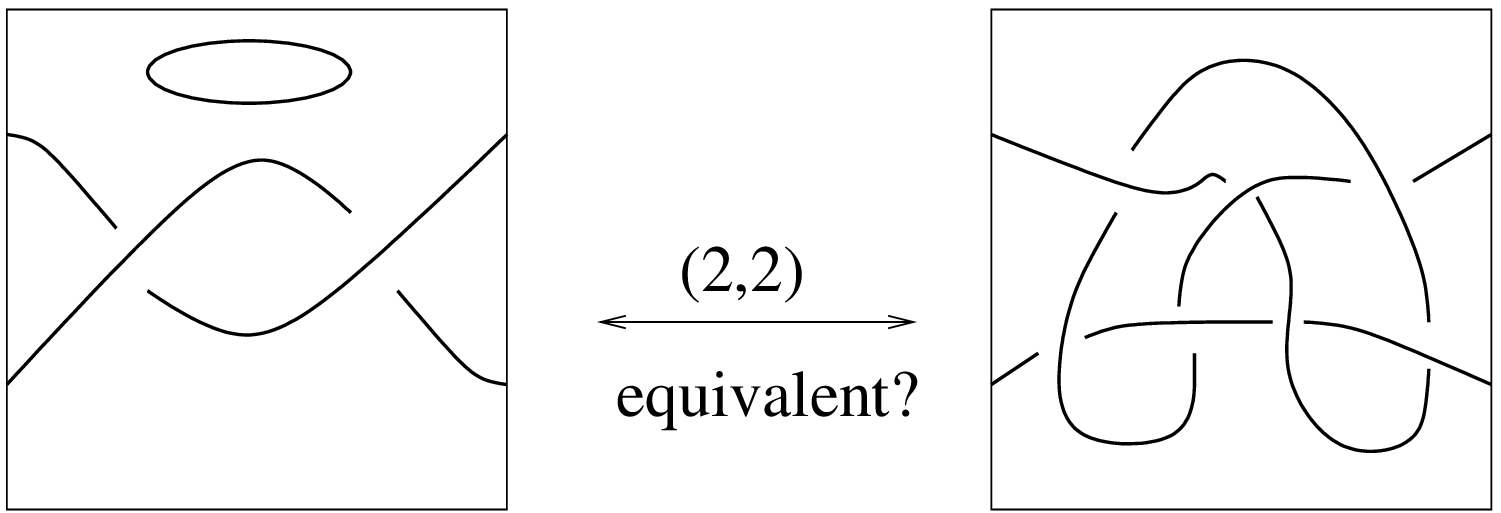,height=3.6cm}}
\ \\
\ \\
\centerline{Fig. 2.10}

A weaker version of the Montesinos-Nakanishi $3$-move conjecture
 has been proved by Bronek Wajnryb in 1985 \cite{Wa-1,Wa-2} (compare 
Theorem 1.16).
\begin{theorem}[Wajnryb]\label{2.10}
Every link can be reduced to a trivial link
 by a finite number of $\pm 3$-moves
and $\Delta_5^4$-moves.
\end{theorem}

 Let me complete this talk by mentioning two generalizations
of the Fox $k$-colorings. 

In the first generalization we consider any commutative
ring with the identity, ${\cal R}$, instead of $\mathbb{Z}_k$. 
We construct $Col_{\cal R}T$
in the same way as before with the relation at
each crossing, Fig.2.1,
having the form $c=2a-b$ in ${\cal R}$.
The skew-symmetric form $\phi$ on ${\cal R}^{2n-1}$, the symplectic
form $\hat\phi$ on ${\cal R}^{2n-2}$ and the homomorphisms $\psi$ and 
$\hat\psi$ are defined in the same manner as before.
Theorem  2.4 generalizes as follows (\cite{D-J-P}):
\begin{theorem}\label{2.11}
Let ${\cal R}$ be a Principal Ideal Domain (PID). Then,
$\hat\psi (Col_{\cal R} T/{\cal R})$ is a virtual Lagrangian
submodule of ${\cal R}^{2n-2}$ with the symplectic form $\hat\phi$.
That is, $\hat\psi (Col_{\cal R} T/{\cal R})$ is a finite index submodule
of a Lagrangian in ${\cal R}^{2n-2}$.
\end{theorem}
This result can be used to analyze embeddability of tangles in links. 
It gives an alternative proof of Theorem 2.2 in \cite{P-S-W} 
in the case of the 2-fold cyclic cover of $B^3$ branched over a tangle.

\begin{example}\label{2.12}
Consider the pretzel tangle $T= (p,-p)$, Fig.2.11, and the 
ring ${\cal R}=\mathbb{Z}$.
Then the virtual Lagrangian $\hat\psi(Col_\mathbb{Z}T/Z)$ 
 has index $p$. Namely, coloring of $T$, as illustrated in 
Fig.2.11, forces us to have $a=b$ and modulo trivial 
colorings the image $\hat\psi(Col_\mathbb{Z}T/\mathbb{Z})$ 
is generated by the vector $(0,p,p,0)= p(e_2+e_3)$. 
The symplectic space $(\mathbb{Z}^{4-2},\hat\phi)$ 
has a basis $(e_1+e_2,e_2+e_3)$. Thus, 
$\hat\psi(Col_\mathbb{Z}T/\mathbb{Z})$ is a 
virtual Lagrangian of index $p$.
\end{example}

\centerline{\psfig{figure=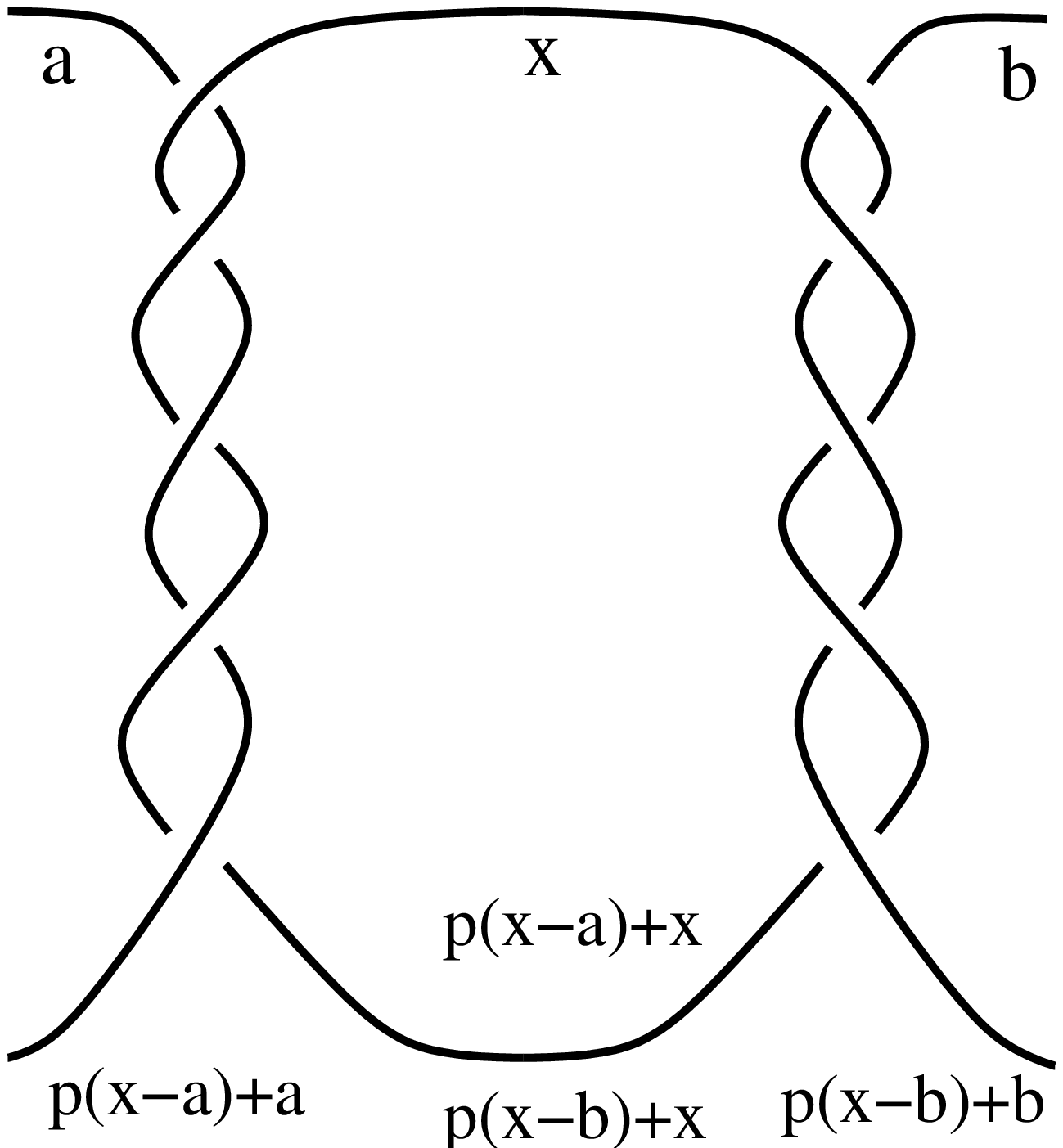,height=4.1cm}} 
\centerline{Fig. 2.11; $p=5$}
\ \\

\begin{corollary}\label{2.13}\ \\
If $\hat\psi(Col_\mathbb{Z}T/\mathbb{Z})$ is a
virtual Lagrangian of index $p>1$, then $T$ does not embed 
in the trivial knot.
\end{corollary}
 \ \\

The second generalization leads to racks and quandles \cite{Joy,F-R}, 
but we restrict our remarks to the abelian case -- Alexander-Burau-Fox
colorings\footnote{The related approach was first outlined in the letter
of J.~W. Alexander to O. Veblen, 1919 \cite{A-V}. Alexander was probably
influenced by Poul Heegaard's dissertation, 1898, which he reviewed for
the French translation \cite{Heeg}. Burau was considering a braid
representation, but locally his relation was the same as that of Fox.
According to J. Birman, Burau learned about the representation 
from Reidemeister or Artin \cite{Ep}, p.330.}.
An ABF-coloring uses colors from a ring ${\cal R}$ with an invertible 
element $t$ (e.g., ${\cal R}= \mathbb{Z}[t^{\pm 1}]$). The relation 
in Fig.2.1 is modified to the relation $c=(1-t)a + tb$ in ${\cal R}$ 
at each crossing of an oriented link diagram; see Fig. 2.12.
\ \\
\ \\
\centerline{\psfig{figure=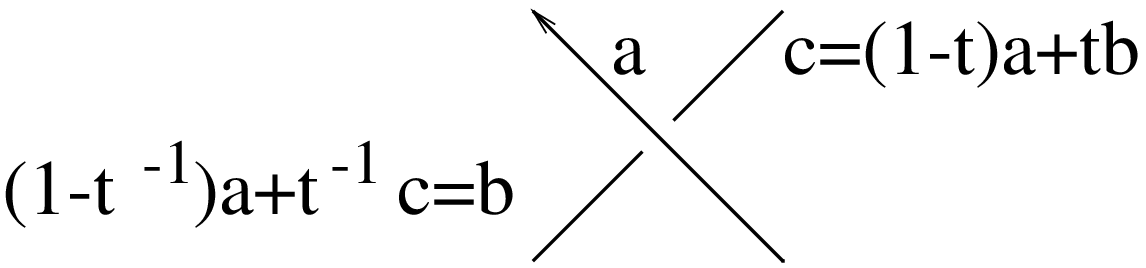,height=2.6cm}}
\centerline{Fig. 2.12}
\ \\

The space ${\cal R}^{2n-2}$ has a natural Hermitian structure \cite{Sq},
but one can also find a symplectic structure and prove a 
version of Theorem 2.11 in this setting \cite{D-J-P}.

\section{Historical Introduction to Skein Modules}\label{3}

In my last talk of the conference, I will discuss
{\it skein modules}, or as I prefer to say more generally, 
{\it algebraic topology based on knots}. It was my mind's child,  
 even if the idea was also conceived by
other people (most notably Vladimir Turaev), and was
envisioned by John H. Conway (as ``linear skein") a
decade earlier.
Skein modules have their origin in the observation made by Alexander
\cite{Al}, that his polynomials ({\it Alexander polynomials}) of three
links, $L_+,L_-$ and $L_0$ in $R^3$ are linearly related (Fig.3.2).

For me it started in Norman, Oklahoma in April of 1987,
when I was enlightened to see that the multivariable
version of the Jones-Conway (Homflypt) polynomial
analyzed by Jim Hoste and Mark Kidwell is really a module
of links in a solid torus (or more generally, in
the connected sum of solid tori).

I would like to discuss today, in more detail,  
skein modules related to the (deformations) 
of 3-moves and the Montesinos-Nakanishi 3-move 
conjecture, but first I will give the general definition
and I will make a short tour of the world of skein modules.

Skein Module is an algebraic object associated with a manifold, 
usually constructed as
a formal linear combination of embedded (or immersed) submanifolds,
modulo locally defined relations.
In a more restricted setting a skein module
is a module associated with a 3-dimensional manifold, by considering
linear combinations of links in the manifold, modulo properly chosen
(skein) relations. It is the main object of the {\it algebraic topology
based on knots}. When choosing relations one takes
into account several factors:
\begin{enumerate}
\item [\textup{(i)}]  Is the module we obtain accessible (computable)?
\item [\textup{(ii)}]
How precise are our modules in distinguishing 3-manifolds and links
in them?
\item [\textup{(iii)}] 
Does the module reflect topology/geometry of a 3-manifold
(e.g., surfaces in a manifold, geometric decomposition of a manifold)?
\item [\textup{(iv)}]
Does the module admit some additional structure (e.g., filtration,
gradation, multiplication, Hopf algebra structure)? Is it leading
to a Topological Quantum Field Theory (TQFT) by taking a finite
dimensional quotient?
\end{enumerate}
\ \\
One of the simplest skein modules is a $q$-deformation of the first
homology group of a 3-manifold $M$, denoted by ${\cal S}_2(M;q)$.
It is based on the skein relation (between oriented framed links in $M$):
$L_+=qL_0$; it also satisfies the framing relation 
$L^{(1)}-qL$, where $L^{(1)}$
denote a link obtained from  $L$ by twisting the framing of $L$ once
in the positive direction. 
 This easily defined skein module ``sees" already
nonseparating surfaces in $M$. These surfaces are responsible for
torsion part of our skein module \cite{P-10}.

There is a more general pattern: most of the analyzed skein modules reflect
various surfaces embedded in a manifold.

The best studied skein modules use skein relations
which worked successfully in the classical
Knot Theory (when defining polynomial invariants of links in $R^3$).
\begin{enumerate}
\item[\textup{(1)}] The Kauffman bracket skein module, KBSM.\\
The skein module based on the {\it Kauffman bracket skein relation},
$L_+ = AL_- + A^{-1}L_{\infty}$, and denoted by 
$S_{2,\infty}(M)$,
is the best understood among the Jones type skein modules.
It can be interpreted as a quantization of the co-ordinate ring
of the character variety of $SL(2,\mathbb{C})$ representations of the
fundamental group of the manifold $M$, \cite{Bu-2,B-F-K,P-S}.
For $M= F\times [0,1]$,
KBSM is an algebra (usually noncommutative). It is finitely generated
algebra for a compact $F$ \cite{Bu-1}, and has no zero divisors \cite{P-S}.
The center of the algebra is generated by boundary 
components of $F$ \cite{B-P,P-S}.
Incompressible tori and 2-spheres in $M$ yield torsion in KBSM;
it is the question of fundamental importance whether other surfaces
could yield torsion as well. The Kauffman bracket skein modules of the 
 exteriors of 2-bridge knots have been recently (April 2004) 
computed by Thang Le \cite{Le}. 
For a 2-bridge (rational) knot $K_{\frac{p}{m}}$ the 
skein module is the free $\mathbb{Z}[A^{\pm 1}]$ module with the basis 
$\{x^iy^j\}$, $0\leq i$, $0\leq j \leq \frac{p-1}{2}$, where $x^iy^j$ 
denotes the element of the skein module represented by the link 
composed of $i$ parallel copies of the meridian curve $x$ and $j$ 
parallel copies of the curve $y$; see Fig.3.1.
Le's theorem generalizes results in \cite{Bu-3} and \cite{B-L}. 
\ \\
\ \\
\centerline{\psfig{figure=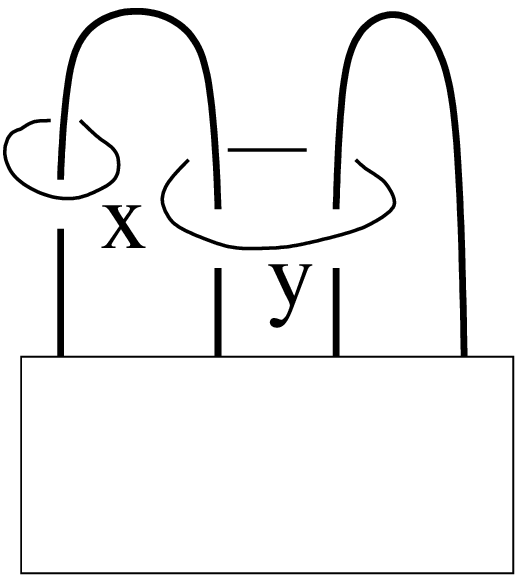,height=5.2cm}}
\centerline{Fig. 3.1}

\item[\textup{(2)}] 
Skein modules based on the Jones-Conway (Homflypt) relation.\\
$v^{-1}L_+ - vL_- = z L_0$, where $L_+,L_-,L_0$ are oriented links 
(Fig.3.2).
These skein modules are denoted by $S_3(M)$ and generalize skein modules
based on Conway relation which were hinted at by Conway. 
For $M= F\times [0,1]$,
$S_3(M)$ is a {\it Hopf algebra} (usually neither
commutative nor co-commutative),
\cite {Tu-2,P-6}. $S_3(F\times [0,1])$ is a free module and can be
interpreted as a quantization \cite{H-K,Tu-1,P-5,Tu-2}.
$S_3(M)$ is related to the algebraic set of
$SL(n,\mathbb{C})$ representations of the
fundamental group of the manifold $M$, \cite{Si}.\\
\ \\
\centerline{\psfig{figure=L+L-L0.eps}}
\centerline{Fig. 3.2}
\item[\textup{(3)}] 
Skein modules based on the {\it Kauffman polynomial} relation.\\
$L_{+1}  + L_{-1} = x (L_0 + L_{\infty})$ (see Fig.3.3) 
and the framing relation $L^{(1)}-aL$.
This module is denoted by $S_{3,\infty}$ and is known 
to be free for $M= F\times [0,1]$.

\item[\textup{(4)}] 
Homotopy skein modules.\\
 In these skein modules, $L_+ = L_-$ for
self-crossings. The best studied example is the q-homotopy skein module
with the skein relation $q^{-1}L_+ -qL_- =zL_0$ for mixed crossings.
For $M= F\times [0,1]$ it is a quantization, \cite{H-P-1,Tu-2,P-11}, and
as noted by Uwe Kaiser they can be almost completely understood using
singular tori technique introduced by Xiao-Song Lin.

\item[\textup{(5)}]
Skein modules based on Vassiliev-Gusarov filtration.\\
We extend the family of knots, $\cal K$, by allowing singular knots,
and resolve a singular crossing by $K_{cr}= K_+ -K_-$. These allow us
to define the Vassiliev-Gusarov filtration:
$... \subset C_3  \subset C_2  \subset C_1  \subset C_0= R{\cal K}$,
where $C_k$ is generated by knots with $k$ singular crossings.
The $k$th Vassiliev-Gusarov skein module is defined to be a quotient:
$W_k(M) = R{\cal K}/C_{k+1}$. The completion of the space of knots
with respect to the Vassiliev-Gusarov filtration, $\hat{R{\cal K}}$,
 is a {\it Hopf algebra} (for $M=S^3$). Functions dual to
Vassiliev-Gusarov skein modules are called {\it finite type} or
{\it Vassiliev invariants} of knots, \cite{P-7}.

\item[\textup{(6)}]
Skein modules based on relations deforming n-moves.\\
${\cal S}_n(M)= R{\cal L}/(b_0L_0 + b_1L_1 + b_2L_2 +...+b_{n-1}L_{n-1})$.
In the unoriented case, we can add to the relation the 
term $b_{\infty}L_{\infty}$ 
to get ${\cal S}_{n,\infty}(M)$, and also,
possibly, a framing relation.  The case $n=4$, on which I am working
with my students  will be described, in greater detail 
in a moment.
\end{enumerate}
\
\\
Examples (1)-(5) gave a short description of skein 
modules studied extensively until now. I will now spend  
more time on two other examples which only recently
have been considered in more detail. The first example is based on
a deformation of the 3-move and the second on the deformation of
the $(2,2)$-move. The first one has been studied by
my students Tsukamoto and Veve. I denote the skein module
described in this example by 
${\cal S}_{4,\infty}$ since it involves (in the skein relation) 
4 horizontal positions and the vertical ($\infty$) smoothing.
\begin{definition}
Let $M$ be an oriented 3-manifold and let ${\cal L}_{fr}$ be 
the set of unoriented
framed links in $M$ (including the empty link, $\emptyset$), and let
$R$ be any commutative ring with identity.  Then we define the
$(4,\infty)$ skein module as:
${\cal S}_{4,\infty}(M;R) = R{\cal L}_{fr}/I_{(4,\infty)}$,
where $I_{(4,\infty})$ is the submodule of $R{\cal L}_{fr}$
generated by the skein relation:\\
$b_0L_0  + b_1L_1 + b_2L_2 + b_3L_3 + b_{\infty}L_{\infty} = 0$
and
the framing relation:\\
$L^{(1)} = a L$ where $a,b_0,b_3$ are invertible elements in $R$ and
$b_1,b_2,b_{\infty}$ are any fixed elements of $R$ (see Fig.3.3).
\end{definition}
\centerline{\psfig{figure=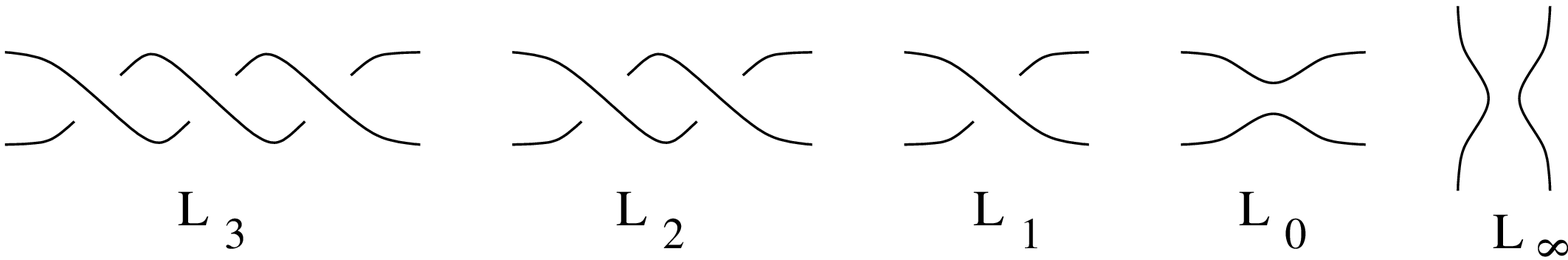,height=2.2cm}}
\centerline{Fig. 3.3}
\ \\

The generalization of the Montesinos-Nakanishi 3-move conjecture
says that ${\cal S}_{4,\infty}(S^3,R)$ is generated by
trivial links\footnote{The counterexamples to the Montesinos-Nakanishi 
3-move conjecture, \cite{D-P-1}, can be used to show that trivial links 
``generically" do not generate ${\cal S}_{4,\infty}(S^3,R)$. This happen, 
for example, if there is a proper ideal ${\cal I}\in R$ such that 
$b_1,b_2 $ and $b_{\infty}$ are in ${\cal I}$.}
 and that for $n$-tangles our skein
module is generated by $f(n,3)$ basic tangles (with
possible trivial components). This would give a
generating set for our skein module of $S^3$ or 
$D^3$ with
$2n$ boundary points (an $n$-tangle).
In \cite{P-Ts} we analyzed extensively the possibility
that trivial links, $T_n$, are linearly independent.
This may happen if $b_{\infty} = 0$ and $b_0b_1=b_2b_3$.
These lead to the following conjecture:
\begin{conjecture}\label{3.2}
\begin{enumerate}
\item [\textup{(1)}] There is a polynomial invariant 
of unoriented links in $S^3$,
$P_1(L) \in Z[x,t]$, which satisfies:
\begin{enumerate}
\item [\textup{(i)}] 
Initial conditions: $P_1(T_n) = t^n$, where $T_n$ is a trivial
link of $n$ components.
\item [\textup{(ii)}] 
Skein relation: $P_1(L_0) + xP_1(L_1) - xP_1(L_2) - P_1(L_3)=0$, 
where $L_0,L_1,L_2,L_3$ is a standard, unoriented skein quadruple
($L_{i+1}$ is obtained from $L_{i}$ by a right-handed half-twist on
two arcs involved in $L_{i}$; compare Fig.3.3).
\end{enumerate}
\item [\textup{(2)}] 
There is a polynomial invariant of unoriented framed links,
$P_2(L) \in Z[A^{\pm 1},t]$ which satisfies:
\begin{enumerate}
\item [\textup{(i)}] Initial conditions: $P_2(T_n) = t^n$,
\item [\textup{(ii)}] 
Framing relation: $P_2(L^{(1)})=-A^3P_2(L)$ where $L^{(1)}$ is
obtained from a framed link $L$ by a positive half twist on its framing.
\item [\textup{(iii)}] 
Skein relation: $P_2(L_0) + A(A^2 + A^{-2})P_2(L_1) +
(A^2 + A^{-2})P_2(L_2) + AP_2(L_3)=0$.
\end{enumerate}
\end{enumerate}
\end{conjecture}
The above conjectures assume that $b_{\infty}=0$ in our
skein relation. Let us consider, for a moment, the possibility that
$b_{\infty}$ is invertible in $R$. Using the ``denominator"
of our skein relation (Fig.3.4) we obtain the relation which
allows us to compute the effect of adding a trivial component
to a link $L$ (we write $t^n$ for the trivial link $T_n$):
$$(*)\ \ \ (a^{-3}b_3 + a^{-2}b_2 + a^{-1}b_1 + b_0 + b_{\infty}t)L=0$$
When considering the ``numerator" of the relation and its 
mirror image (Fig.3.4)
we obtain formulas for Hopf link summands, and because the unoriented
Hopf link is amphicheiral we can eliminate it from
our equations to get the following formula (**):
$$b_3(L\#H) + (ab_2 + b_1t +a^{-1}b_0 + ab_{\infty})L =0.$$
$$b_0(L\#H) + (a^{-1}b_1 +b_2t + ab_3 + a^2b_{\infty})L =0.$$
$$(**)\ \ \ ((b_0b_1 - b_2b_3)t + (a^{-1}b_0^2 -a b_3^2) +
(ab_0b_2 - a^{-1}b_1b_3) + b_{\infty}(ab_0 -a^2b_3))L = 0.$$
\ \\
\centerline{\psfig{figure=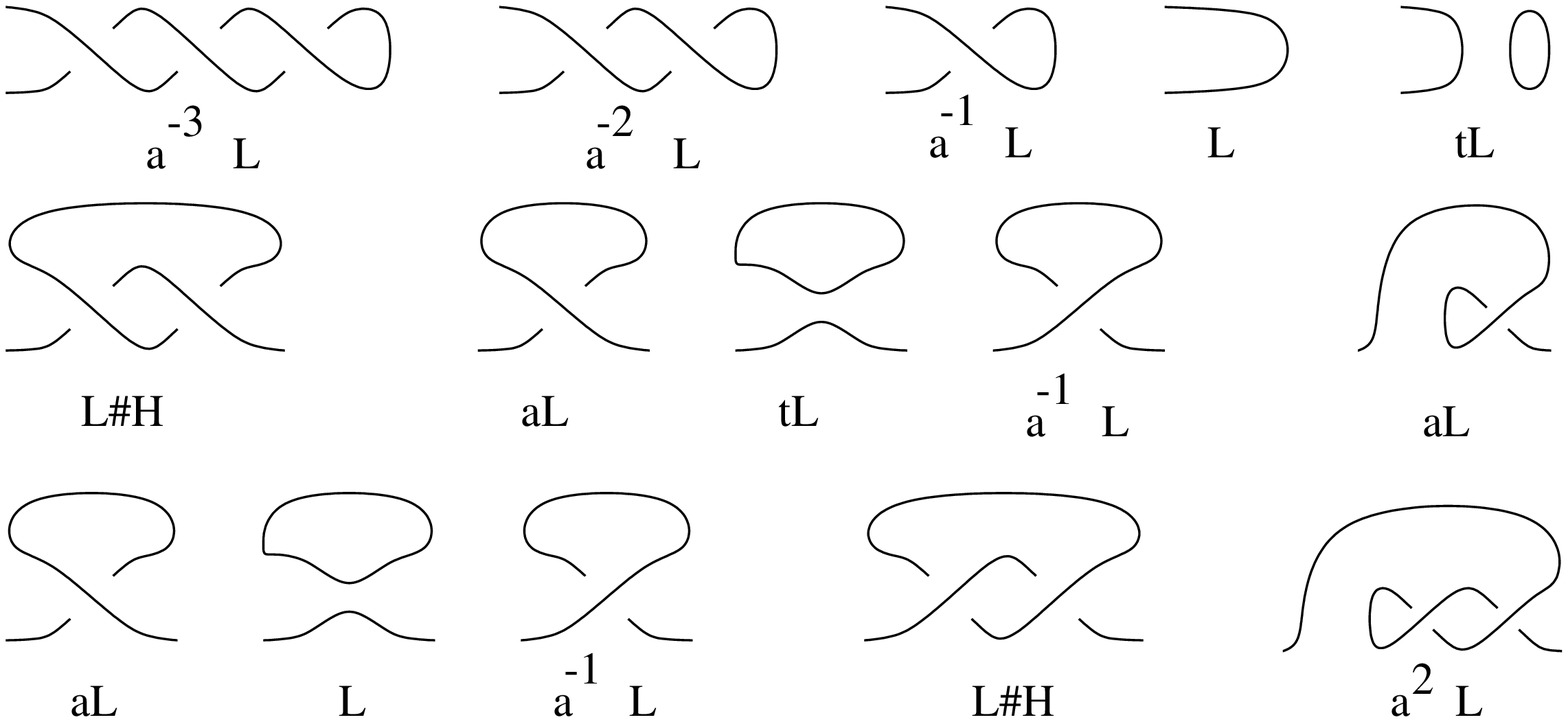,height=5.6cm}}
\centerline{Fig. 3.4}
\ \\
It is possible that ($*$) and ($**$) are the only relations in the module.
More precisely, we ask whether ${\cal S}_{4,\infty}(S^3;R)$
is the quotient ring $R[t]/({\cal I})$ where $t^i$ represents
the trivial link of $i$ components and ${\cal I}$ is the
ideal generated by ($*$) and ($**$) for $L=t$.\
The interesting substitution which satisfies the relations is 
$b_0=b_3=a=1$, $b_1=b_2=x$, $b_{\infty}=y$. This may lead to 
a new polynomial invariant (in $\mathbb{Z}[x,y]$) of unoriented links 
in $S^3$ satisfying the skein relation 
$L_3 +xL_2 + xL_1 + L_0 + yL_{\infty} = 0$.\footnote{This speculation 
should be modified keeping in mind the fact that the Montesinos-Nakanishi 
3-move conjecture does not hold \cite{D-P-1}.}

What about the relations to the Fox colorings?
One such a relation, that was already mentioned,
is the use of 3-colorings to estimate the number of
basic n-tangles (by $\prod_{i=1}^{n-1} (3^i +1)$) for
the skein module ${\cal S}_{4,\infty}$. I am also convinced that
${\cal S}_{4,\infty}(S^3;R)$ contains full information about the space of
Fox 7-colorings. It would be a generalization of the fact that
the Kauffman bracket polynomial contains information about 3-colorings
and the Kauffman polynomial contains information about 
5-colorings. In fact, Fran{\c c}ois Jaeger told me that he 
knew how to form a short  skein relation 
(of the type  $(\frac{p+1}{2},\infty)$) involving 
 spaces of $p$-colorings. 
Unfortunately,
Fran{\c c}ois died prematurely in 1997 and I do not know 
how to prove his statement\footnote{If 
$col_p(L)=|Col_p(L)|$ denotes the order of the space of Fox 
$p$-colorings  of the link $L$, then
among $p+1$ links $L_0,L_1,...,L_{p-1}$, and $L_{\infty}$,
$p$ of them has the same order  $col_p(L)$ 
and one has its order $p$ times larger
\cite{P-3}. This leads to the relation of type $(p,\infty)$.
The relation between Jones polynomial (or the Kauffman bracket)
and $col_3(L)$ has the form: $col_3(L)= 3|V(e^{\pi i/3})|^2$ and
the formula relating the Kauffman polynomial and $col_5(L)$ 
has the form: $col_5(L)= 5|F(1,e^{2\pi i/5} + e^{-2\pi i/5})|^2$.
This seems to suggest that the formula discovered by Jaeger 
involved Gaussian sums.}.

Finally, let me shortly  describe the skein module
related to the deformation of $(2,2)$-moves.
Because a $(2,2)$-move is equivalent to the rational
$\frac{5}{2}$-move, I will denote the skein module
by ${\cal S}_{\frac{5}{2}}(M;R)$.
\begin{definition}\label{3.3}
Let $M$ be an oriented 3-manifold. Let  ${\cal L}_{fr}$ be
 the set of unoriented framed links in $M$ (including the 
empty link, $\emptyset$) and let $R$ be any commutative 
ring with identity.  We define the $\frac{5}{2}$-skein module as 
 ${\cal S}_{\frac{5}{2}}(M;R) = R{\cal L}_{fr}/(I_{\frac{5}{2}})$
where $I_{\frac{5}{2}}$ is the submodule of $R{\cal L}_{fr}$
generated by the skein relation:\\
\textup{(i)} \ \ $b_2L_2 + b_1L_1 + b_0L_0  + b_{\infty}L_{\infty} +
b_{-1}L_{-1} + b_{-\frac{1}{2}}L_{-\frac{1}{2}} = 0$,\\
 its mirror image:\\
$(\bar{i})$ \ \ $b'_2L_2 + b'_1L_1 + b'_0L_0  + b'_{\infty}L_{\infty} +
b'_{-1}L_{-1} + b'_{-\frac{1}{2}}L_{-\frac{1}{2}} = 0$\\
 and the framing relation:\\
$L^{(1)} = a L$, where $a,b_2,b_2',b_{-\frac{1}{2}},
b'_{-\frac{1}{2}}$ are invertible elements in $R$ and 
$b_1,b'_1,b_0,b'_0$, $b_{-1},b'_{-1},b_{\infty}$, 
and $b'_{\infty}$ are any fixed elements of $R$. 
The links $L_2,L_1,L_0,L_{\infty},L_{-1},$ $L_{\frac{1}{2}}$ 
and $L_{-\frac{1}{2}}$
are illustrated in Fig.3.5.\footnote{Our notation is based on
Conway's notation for rational tangles. However, it differs
from it by a sign change. The reason is that the Conway
convention for a positive crossing is generally not
used in the setting of skein relations.}
\end{definition}
\centerline{\psfig{figure=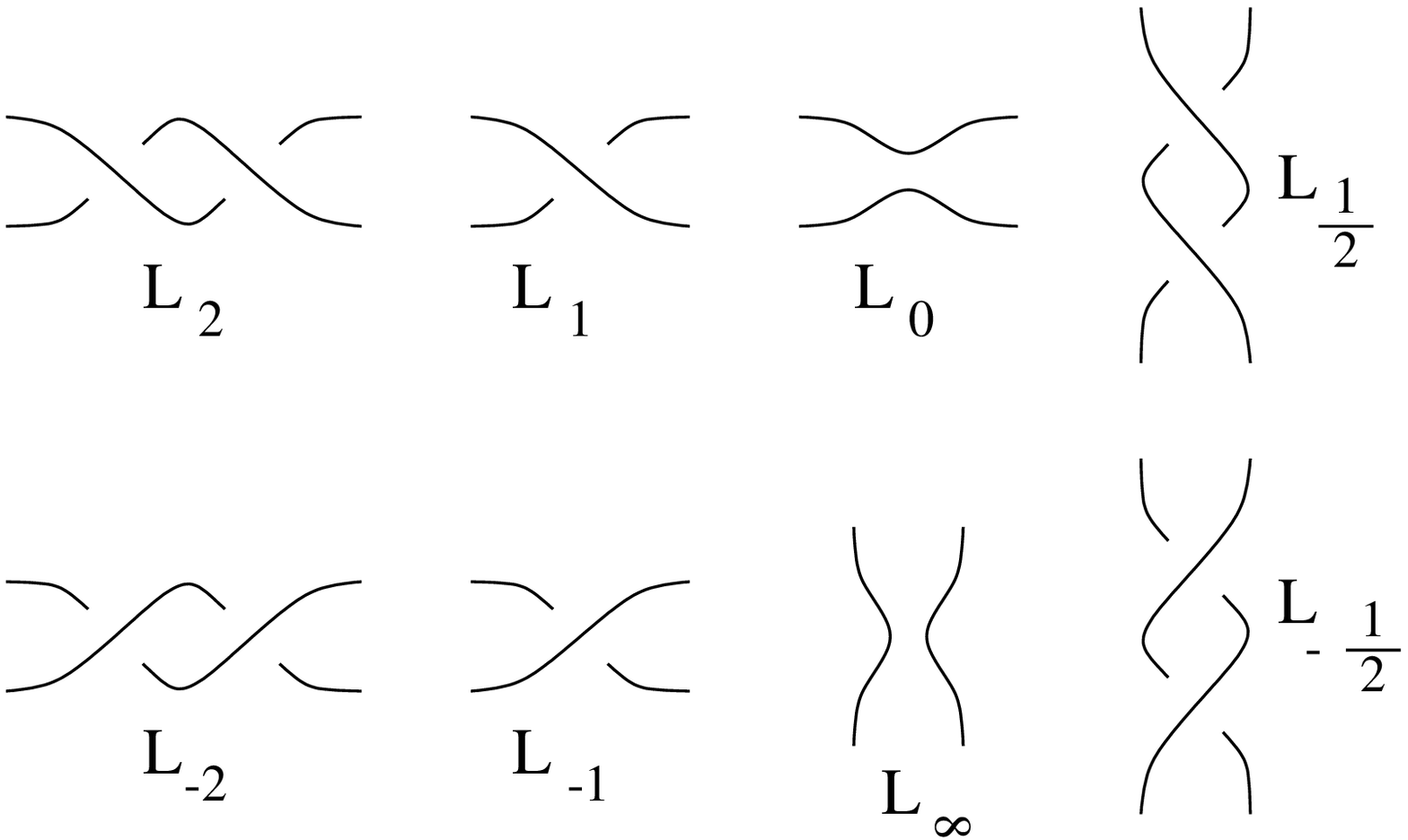,height=5.7cm}}
\centerline{Fig. 3.5}
\ \\

If we rotate the figure from the relation (i) we obtain:\\
(i')  $b_{-\frac{1}{2}}L_2 + b_{-1}L_1 + b_{\infty}L_0  + b_{0}L_{\infty} +
b_{1}L_{-1} + b_{2}L_{-\frac{1}{2}} = 0$\\
One can use (i) and (i') to eliminate $L_{-\frac{1}{2}}$ 
and to get the relation:\\
$(b_2^2- b^2_{-\frac{1}{2}}) L_2 +  
(b_1b_2 - b_{-1}b_{-\frac{1}{2}})L_1 +
((b_0b_2 - b_{\infty}b_{-\frac{1}{2}})L_0 +
(b_{-1}b_2 - b_{1}b_{-\frac{1}{2}})L_{-1} +
(b_{\infty}b_2 - b_{0}b_{-\frac{1}{2}})L_{\infty}=0$.\\
Thus, either we deal with the shorter relation (essentially the one in
the fourth skein module described before) or all
coefficients are equal to 0 and therefore (assuming 
that there are no zero divisors in $R$)
$b_2 = \varepsilon b_{-\frac{1}{2}}$, $b_1 = \varepsilon b_{-1}$, 
and $b_0 = \varepsilon b_{\infty}$. Similarly, we would get:
$b'_2 = \varepsilon b'_{-\frac{1}{2}}$, $b'_1 = \varepsilon b'_{-1}$, 
and $b'_0 = \varepsilon b'_{\infty}$, where $\varepsilon = \pm 1$.
Assume, for simplicity, that $\varepsilon =1$. Further relations
among coefficients follow from the computation of the 
Hopf link component using the amphicheirality of the
unoriented Hopf link. 
Namely, by comparing diagrams in Figure 3.6 and their mirror images we get
$$L\#H= -b_2^{-1}(b_1(a + a^{-1}) + a^{-2}b_2 + b_0(1+T_1))L$$
$$L\#H= -{b'}_2^{-1}(b'_1(a + a^{-1}) + a^{2}b'_2 + b'_0(1+T_1))L.$$
Possibly, the above equalities give the only other relations among
coefficients (in the case of $S^3$). I would present below
the simpler question (assuming $a=1, b_x=b'_x$ and writing
$t^n$ for $T_n$).\\
\ \\
\centerline{\psfig{figure=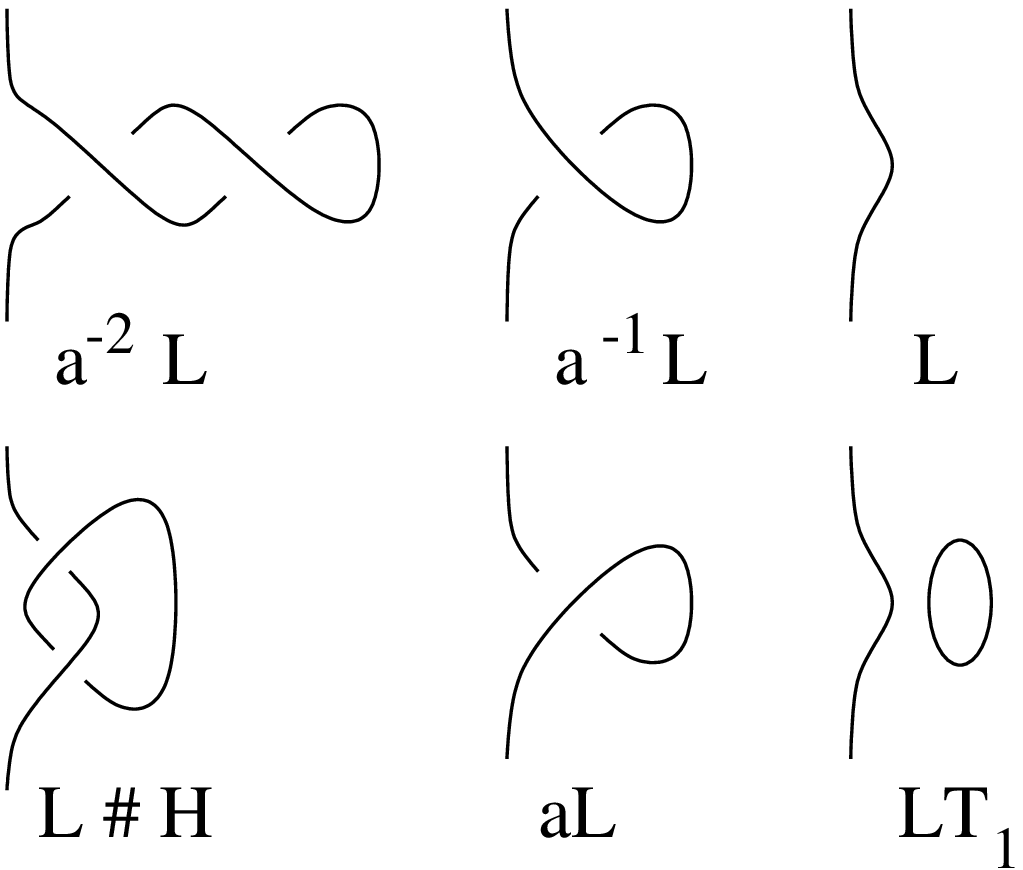,height=6.4cm}}
\centerline{Fig. 3.6}

\begin{question}\label{3.4}
Is there a polynomial invariant of unoriented links in $S^3$,
$P_{\frac{5}{2}}(L) \in {\mathbb Z}[b_0,b_1,t]$, which satisfies the
following conditions?
\begin{enumerate}
\item [\textup{(i)}] 
Initial conditions: $P_{\frac{5}{2}}(T_n) = t^n$, where $T_n$ is a trivial
link of $n$ components.
\item [\textup{(ii)}] 
Skein relations
$$P_{\frac{5}{2}}(L_2) + b_1P_{\frac{5}{2}}(L_1) +  b_0P_{\frac{5}{2}}(L_0) + 
b_0P_{\frac{5}{2}}(L_{\infty}) + b_1P_{\frac{5}{2}}(L_{-1}) +  
P_{\frac{5}{2}}(L_{-\frac{1}{2}})
=0.$$
$$P_{\frac{5}{2}}(L_{-2}) + b_1P_{\frac{5}{2}}(L_{-1}) +  
b_0P_{\frac{5}{2}}(L_0) + 
b_0P_{\frac{5}{2}}(L_{\infty}) + b_1P_{\frac{5}{2}}(L_{1}) +  
P_{\frac{5}{2}}(L_{\frac{1}{2}})
=0.$$ 
\end{enumerate}
\end{question}
Notice that by taking the difference of our skein relations one gets
the interesting identity:
$$P_{\frac{5}{2}}(L_2)  - P_{\frac{5}{2}}(L_{-2}) =
P_{\frac{5}{2}}(L_{\frac{1}{2}}) - P_{\frac{5}{2}}(L_{-\frac{1}{2}}).$$
Nobody has yet studied the skein module 
${\cal S}_{\frac{5}{2}}(M;R)$ seriously so everything that you can find will
be a new research, even a table of the polynomial $P_{\frac{5}{2}}(L)$
for small links, $L$.\\
I wish you luck!


\section{Added in proof -- the Montesinos-Nakanishi 3-move conjecture}
A preliminary calculation performed by my student Mietek D{\c a}bkowski
(February 21, 2002) 
shows that the Montesinos-Nakanishi 3-move conjecture does not hold
for the Chen link (Fig.1.14). Below is the text of the abstract we have sent
for the Knots in Montreal conference organized by Steve Boyer and
Adam Sikora in April 2002.\\
\ \\
Authors: Mieczys{\l}aw D{\c a}bkowski, J\'ozef H. Przytycki (GWU)\\
Title: Obstructions to the Montesinos-Nakanishi 3-move conjecture.\\
Yasutaka Nakanishi asked in 1981 whether a 3-move is an unknotting
operation. This question is called, in the Kirby's problem list,
{\it the Montesinos-Nakanishi Conjecture}.
Various partial results have been obtained by Q.Chen, Y.Nakanishi,
J.Przytycki and T.Tsukamoto. Nakanishi and Chen presented examples
which they couldn't reduce (the Borromean rings and the closure
of the square of the center of the fifth braid group, $\bar\gamma$,
respectively). The only tool, to analyze 3-move equivalence, till
1999, was the Fox 3-coloring (the number of Fox 3-colorings is unchanged
by a 3-move). It allowed to distinguish different
trivial links but didn't separate Nakanishi and Chen examples
from trivial links.
The group of 3-colorings of a link $L$ corresponds
to the first homology group with $Z_3$ coefficients
of the double branched cover of a link $L$,
$M_L^{(2)}$, i.e. $$Tri(L) = H_1(M_L^{(2)},Z_3)\oplus Z_3$$
We find more delicate invariants of 3-moves
using homotopy in place homology and we consider the fundamental
group of $M_L^{(2)}$.

We define an $n$th Burnside group of a link as the quotient of the
fundamental group of the double branched cover of the link divided
by all relations of the form $a^n=1$. For $n=2,3,4,6$ the
quotient group is 
finite\footnote{Burnside groups of links are instances of groups of
finite exponents. Our method of analysis of tangle moves rely on the well
developed theory of classical Burnside groups and the associated Lie
rings. 
A group $G$ is of a finite exponent if there is a finite integer
$n$ such that $g^{n}=e$ for all $g\in G$. If, in addition, there is no
positive integer $m<n$ such that $g^{m}=e$ for all $g\in G$, then we say
that $G$ has an exponent $n$.
Groups of finite exponents were considered for the first time by 
Burnside in 1902 \cite{Bur}. In particular, Burnside himself 
was interested in the case when $G$ is
a finitely generated group of a fixed exponent. He asked the question, known
as the Burnside Problem, whether there exist infinite and finitely generated
groups $G$ of finite exponents. \newline
Let $F_{r}=\langle x_{1},\,x_{2},\,\dots ,\,x_{r}\,|\,-\rangle $ be the free
group of rank $r$ and let $B(r,n)=F_{r}/N$, where $N$ is the normal subgroup
of $F_{r}$ generated by $\{g^{n}\,|\,g\in F_{r}\}$. The group $B(r,n)$ is
known as the $r$th generator Burnside group of exponent $n$. In this
notation, Burnside's question can be rephrased as follows.
For what values of $r$ and $n$ is the Burnside group $B(r,n)$ finite?
$B(1,n)$ is a cyclic group $Z_n$. $B(r,2) = Z_2^n$. Burnside proved that $%
B(r,3)$ is finite for all $r$ and that $B(2,4)$ is finite. In 1940 Sanov
proved that $B(r,4)$ is finite for all $r$, and in 1958 M.Hall
proved that $B(r,6)$ finite for all $r$.  
However, it was proved by Novikov and Adjan in 1968 
that $B(r,n)$ is infinite whenever $r > 1$ and $n$ is an odd and $n
\geq 4381$ (this result was later improved by Adjan, who showed
that $B(r,n)$ is infinite if $r > 1$ and $n$ odd and $n \geq 665$).
Sergei Ivanov
proved that for $k \geq 48$ the group $B(2,2^k)$ is infinite.
Lys\"enok found that $B(2,2^k)$ is infinite for $k \geq 13$.
 It is still an open problem though whether, for example, $B(2,5)$,
$B(2,7)$ or $B(2,8)$ are infinite or finite \cite{VL,D-P-3}.}.

The third Burnside group of a link is unchanged by 
3-moves\footnote{$p$th Burnside group is preserved by 
$\frac{ps}{q}$-moves. This fact allows us to disprove Conjecture 2.2
 \cite{D-P-3}.}.

In the proof we use the "core" presentation of the group from the diagram;
that is arcs are generators and each crossing gives a relation
$c=ab^{-1}a$ where $a$ corresponds to the overcrossing and $b$ and $c$ to
undercrossings.

The Montesinos-Nakanishi 3-move conjecture
does not hold for Chen's example $\hat \gamma$.

To show that $\hat \gamma$ has different third Burnside group
 than any trivial link it suffices to show
that the following element, $P$, of the Burnside free group
 $B(4,3)=\{x,y,z,t: (a)^3\}$ is nontrivial:
$P=uwtu^{-1}w^{-1}t^{-1}$
where \\
 $u=xy^{-1}zt^{-1}$ and $w=x^{-1}yz^{-1}t$.

With the help of GAP it has been achieved!! (Feb. 21, 2002).

We have confirmed our calculation using also computer algebra system 
Magnus.

\ \\
Department of Mathematics\\
George Washington University  \\
Washington, DC 20052 \\
USA\\
e-mail: przytyck@gwu.edu
\end{document}